\documentclass[reqno,11pt]{amsart}

\usepackage{amsmath,amssymb}
\usepackage{amsfonts}
\usepackage{verbatim}
\usepackage[frenchb]{babel}
\usepackage{times}
\usepackage{enumerate}

\usepackage{graphicx}
\usepackage{epsfig}
\input{amssym.def}

\theoremstyle{plain} 

\newtheorem{theoreme}{Th{\'e}or{\`e}me}
\newtheorem{proposition}{Proposition}[section]
\newtheorem{lemme}[proposition]{Lemme}
\newtheorem{corollaire}[proposition]{Corollaire}

\theoremstyle{definition}

\newtheorem{definition}[proposition]{D{\'e}finition}  

\newtheorem{remarque}[proposition]{Remarque}

\def\C{\mathbb C}

\def\H{\mathbb H}
\def\N{\mathbb N}
\def\P{\mathbb P}
\def\Q{\mathbb Q}
\def\R{\mathbb R}

\def\cC{{\mathcal C}}
\def\cD{{\mathcal D}}

\def\cM{{\mathcal M}}
\def\cO{{\mathcal O}}
\def\cP{{\mathcal P}}

\def\cS{{\mathcal S}}
\def\cT{{\mathcal T}}

\def\heta{\widehat{\eta}}

\def\hg{\widehat{g}}

\newcommand{\bz}{{\bar{z}}}
\newcommand{\bpartial}{{\bar{\partial}}}

\renewcommand{\a}{\alpha}
\renewcommand{\b}{\beta}

\newcommand{\e}{\varepsilon}

\renewcommand{\=}{ : = }

\newcommand{\diam}{\mathrm{diam}}

\newcommand{\can}{\mathrm{can}}

\newcommand{\diag}{\mathrm{Diag}}
\newcommand{\Gal}{\mbox{\rm Gal}}

\newcommand{\blim}{\overline{\lim}}
\newcommand{\lip}{\mathrm{Lip}\,}

\def\pc{{{\mathbb P}^1(\C)}}

\def\pp{{{\mathbb P}^1(\C_p)}}
\def\pv{{{\mathbb P}^1(\C_v)}}

\def\berp{{{\mathsf P}^1({\mathbb C}_p)}}
\def\berv{{{\mathsf P}^1({\mathbb C}_v)}}
\def\berCp{{\mathsf A}^1({\mathbb C}_p)}

\newcommand{\lpar}{(\!(\,}
\newcommand{\rpar}{\, )\!)}

\def\ov{\overline}
\def\pie{{\pi_{\varepsilon}}}

\def\He{{H_{\varepsilon}}}

\def\hn{{h_{\rm nv}}}

\newcommand{\ie}{i.e.\ }

\begin{document}
%
%%%%%%%%%%%%%%%%%%%%%%%%%%%%%%%%%%%%%%%%%%%%%%%%%%%%%%%%%%%%%%%%%%
%

\setcounter{tocdepth}{1}

\title[Points de petite hauteur]{Equidistribution quantitative des
  points de petite hauteur sur la droite projective}
\date{\today}

\author{Charles Favre \and Juan Rivera-Letelier}

\address{CNRS et Institut de Math{\'e}matiques de Jussieu\\
         Case 7012, 2  place Jussieu\\
         F-75251 Paris Cedex 05\\
         France}
\email{favre@math.jussieu.fr}

\address{Departamento de Matem{\'a}tica\\
         Universidad Cat{\'o}lica del Norte\\
         Casilla 1280\\
         Antofagasta - Chile}
\email{juanrive@ucn.cl}

\thanks{Le premier auteur tient {\`a} remercier chaleureusement le project
  MECESUP UCN0202, ainsi que l'ACI ``Syst{\`e}mes Dynamiques Polynomiaux''
  qui ont permis son s{\'e}jour {\`a} l'Universit{\'e} Catholique d'Antofagasta.
  Le deuxi{\`e}me auteur est partiellement soutenu par le projet FONDECYT
  N~1040683. Enfin, nous remercions le rapporteur pour sa lecture
  d{\'e}taill{\'e}e de l'article.}


%
%%%%%%%%%%%%%%%%%%%%%%%%%%%%%%%%%%%%%%%%%%%%%%%%%%%%%%%%%%%%%%%%%%
%
~
\maketitle

%
%%%%%%%%%%%%%%%%%%%%%%%%%%%%%%%%%%%%%%%%%%%%%%%%%%%%%%%%%%%%%%%%%%
%

\begin{abstract}
  Nous introduisons une classe de hauteurs ad{\'e}liques sur la droite
  projective dont nous donnons une estimation du minimum essentiel, et
  pour lesquelles nous d{\'e}montrons un r{\'e}sultat d'{\'e}quidistribution
  des points de petites hauteurs en toutes les places (finies et
  infinies), avec estimation pr{\'e}cise de la vitesse de convergence.  A
  toute fraction rationnelle $R$ en une variable et d{\'e}finie sur un
  corps de nombres $K$, est associ{\'e}e une hauteur normalis{\'e}e sur sa
  clot{\^u}re alg{\'e}brique $\ov{K}$. Nous montrons que ces hauteurs
  dynamiques sont ad{\'e}liques en notre sens, et en d{\'e}duisons des
  r{\'e}sultats d'{\'e}quidistribution de pr{\'e}images par $R$ en toutes
  les places. Notre approche suit celle de Bilu, et s'appuie sur la
  th{\'e}orie du potentiel dans $\C$, ainsi que dans l'espace de Berkovich
  associ{\'e} {\`a} la droite projective de $\C_p$, pour chaque nombre premier
  $p$.
\end{abstract}

\renewcommand{\abstractname}{Abstract}
\begin{abstract}
  We introduce a new class of adelic heights on the projective line.
  We estimate their essential minimum and prove a result of
  equidistribution (at every place) for points of small height with
  estimates on the speed of convergence.  To each rational function $R$
  in one variable and defined over a number field $K$, is associated a
  normalized height on the algebraic closure of $K$. We show that
  these dynamically defined heights are adelic in our sense, and deduce
  from this equidistribution results for preimages of points under $R$
  at every place of $K$.  Our approach follows that of Bilu, and relies on
  potential theory in the complex plane, as well as in the Berkovich
  space associated to the projective line over $\C_p$, for each prime
  $p$.
\end{abstract}

\tableofcontents

%
%%%%%%%%%%%%%%%%%%%%%%%%%%%%%%%%%%%%%%%%%%%%%%%%%%%%%%%%%%%%%%%%%%
%

\section{Introduction}
L'objet de cet article est de d{\'e}finir une large classe de hauteurs sur
la droite projective sur un corps de nombre, et de montrer de mani{\`e}re
quantitative, c'est-{\`a}-dire avec estimation des restes, que toute suite
de points alg{\'e}briques dont la hauteur tend vers z{\'e}ro s'{\'e}quidistribue
en toutes les places finies et infinies, selon une mesure ne
d{\'e}pendant que de la hauteur. Nous appliquons de plus ces r{\'e}sultats {\`a}
l'{\'e}tude dynamique de fractions rationnelles {\`a} coefficients
alg{\'e}briques.

Le premier r{\'e}sultat d'{\'e}qui\-distri\-bution des points de petite
hauteur a {\'e}t{\'e} obtenu par Szpiro-Ullmo-Zhang~\cite{SUZ} en lien avec la
conjecture de Bogomolov, et concerne l'{\'e}qui\-distri\-bution des points
de petite hauteur dans les vari{\'e}t{\'e}s ab{\'e}liennes par rapport {\`a} la mesure
de Haar.  Cet article fondateur a inspir{\'e} depuis lors de nombreux
travaux. Tout d'abord par Bilu qui s'est int{\'e}ress{\'e} dans~\cite{bilu} au
cas de la hauteur standard dans les espaces projectifs; puis {\`a} Rumely
qui a {\'e}tendu l'approche de Bilu pour une classe de hauteurs sur la
droite projective issues de la th{\'e}orie du potentiel complexe.
Autissier~\cite{Autissier1} a ensuite d{\'e}montr{\'e} une vaste
g{\'e}n{\'e}ralisation des th{\'e}or{\`e}mes de Bilu et Szpiro-Ullmo-Zhang dans le cas
des courbes d{\'e}finies sur un corps de nombre et en dimension
sup{\'e}rieure. Enfin plus r{\'e}cemment, Baker et Hsia dans~\cite{Baker-Hsia}
ont d{\'e}montr{\'e} des r{\'e}sultats d'{\'e}quidistribution aux places finies dans
un contexte dynamique pour une classe particuli{\`e}re de polyn{\^o}mes. Nous
renvoyons {\`a}~\cite{Ullmo} pour des r{\'e}f{\'e}rences plus compl{\`e}tes
concernant d'autres r{\'e}sultats d'{\'e}quidistribution en
arithm{\'e}tique.

Deux approches parall{\`e}les ont jusqu'{\`a} pr{\'e}sent {\'e}t{\'e} privil{\'e}gi{\'e} dans ces
probl{\`e}mes d'{\'e}quidistribution arithm{\'e}tique. La premi{\`e}re d{\'e}velopp{\'e}e par
Szpiro-Ullmo-Zhang et poursuivie par Autissier est d'inspiration
g{\'e}om{\'e}trique.  Les propri{\'e}t{\'e}s d'{\'e}qui\-distribution r{\'e}sulte dans ce
cadre d'un th{\'e}or{\`e}me d'Hilbert-Samuel arithm{\'e}tique convenablement
{\'e}nonc{\'e} en th{\'e}orie d'Arakelov.

L'autre approche initi{\'e}e par Bilu, et adapt{\'e}e par Rumely et
Baker-Hsia part d'une interpr{\'e}tation des hauteurs en termes de la
th{\'e}orie du potentiel. C'est celle-ci que nous allons adopter dans la
suite.

Un ingr{\'e}dient nouveau et important dans notre approche est
l'utilisation intensive d'une th{\'e}orie du potentiel convenable sur la
droite projective d{\'e}finie sur un corps $p$-adique. La topologie
$p$-adique {\'e}tant totalement discontinue, elle se pr{\^e}te de fait mal {\`a}
l'analyse et tout particuli{\`e}rement {\`a} la th{\'e}orie de la mesure et du
potentiel. On est donc naturellement amen{\'e} {\`a} travailler dans un espace
plus gros dont la topologie est plus maniable: la droite projective au
sens de Berkovich. Cet espace est un arbre r{\'e}el muni d'une topologie
compacte, et on peut y d{\'e}velopper une th{\'e}orie du potentiel
compl{\`e}tement analogue au cas complexe. Il a ainsi {\'e}t{\'e} construit par
M.~Jonsson et le premier auteur dans~\cite{FJ} un op{\'e}rateur de
Laplace, et ses propri{\'e}t{\'e}s caract{\'e}ristiques ont d{\'e}j{\`a} {\'e}t{\'e} utilis{\'e} avec
succ{\`e}s pour la construction de mesures invariantes en dynamique
$p$-adique, voir~\cite{FR}.  Nous donnons ici d'autres {\'e}l{\'e}ments de
th{\'e}orie du potentiel $p$-adique, mais nous nous sommes restreints {\`a}
ceux n{\'e}cessaires {\`a} l'{\'e}nonc{\'e} et {\`a} la preuve de nos r{\'e}sultats
principaux.  Nous renvoyons aux travaux ind{\'e}pendants de
Baker-Rumely~\cite{BR2} pour une {\'e}tude plus approfondie, ainsi qu'aux
travaux de th{\`e}se r{\'e}cents d'Amaury Thuillier~\cite{thuillier}. Enfin
mentionnons que Chambert-Loir~\cite{ACL2} a donn{\'e} quelques {\'e}l{\'e}ments
pour {\'e}tendre cette th{\'e}orie en dimension sup{\'e}rieure.

%%%%%%%%%%%%%%%%%%%%%%%%%%%%%%%%%%%%%%%%%%%%%%%%%%%%%%%%%%%%%%%%%%%%%%%
\subsection{Places, compl{\'e}tions et hauteurs ad{\'e}liques.}\label{S-completions}
Commen\c cons tout d'abord par mettre en place quelques notations
avant de d{\'e}finir les hauteurs que nous {\'e}tudierons par la suite.

Soit $M_\Q$ l'ensemble constitu{\'e} de tous les nombres entiers premiers
auquel on ajoute~$\infty$.  On d{\'e}signe par $| \cdot |_\infty$ la norme
usuelle sur $\Q$ et pour un nombre premier $p$ on d{\'e}signe par $| \cdot
|_p$ la norme $p$-adique, normalis{\'e}e par $|p|_p = p^{-1}$.  Pour toute
extension finie $K$ de $\Q$, on d{\'e}signera par $M_K$ la collection de
toutes les normes sur $K$ qui {\'e}tendent l'une des normes $| \cdot |_v$,
avec $v \in M_\Q$.  Un {\'e}l{\'e}ment de $M_K$ est appel{\'e} place.  Pour toute
place $v$, on notera encore $| \cdot |_v$ la norme sur $K$
correspondante.  On dira que $v$ est {\it infinie} lorsque sa
restriction {\`a} $\Q$ co{\"\i}ncide avec $ |\cdot|_\infty$ et que $v$ est {\it
  finie} sinon.

Pour tout  $v\in M_K$, on note $K_v$ le compl{\'e}t{\'e} du corps
valu{\'e} $(K, | \cdot |_v)$, et $\Q_v$ l'adh{\'e}rence de $\Q$ dans
$K_v$. On pose $N_v = [K_v : \Q_v] / [K : \Q]$, ainsi que $\| \cdot
\|_v = | \cdot |_v^{N_v}$.  La norme $| \cdot |_v$ s'{\'e}tend de
fa{\c c}on unique {\`a} la clot{\^u}re alg{\'e}brique de $K_v$.  On
d{\'e}signera par $\C_v$ le compl{\'e}t{\'e} de ce corps. Le corps $\C_v$
est alors complet et alg{\'e}briquement clos. En toute place infinie,
il est isomorphe au corps des nombres complexes~$\C$.

Aux places finies, $(\C_v, |\cdot|_v)$ est {\`a} la fois totalement
discontinu et non localement compact, ce qui rend d{\'e}licat toute
analyse sur cet espace. Pour contourner ces difficult{\'e}s,
suivant~\cite{Ber} on d{\'e}finit la droite projective au sens de
Berkovich $\berv$ comme la compl{\'e}tion (pour une m{\'e}trique convenable)
de l'ensemble des boules de rayon fini ou infini $B(z,r) = \{w ; \;
|z-w|_v \le r\}$ dans $\C_v$.  Cet espace est naturellement un arbre
r{\'e}el m{\'e}trique dans lequel les ensembles de la forme $\{ B(z,r) ; r\in
[0, +\infty]\}$ sont des segments, et dont le bord {\`a} l'infini
s'identifie canoniquement {\`a} la la droite projective standard. On
peut de plus le munir d'une topologie qui le rend localement connexe
et localement compact. Il est alors possible de d{\'e}finir une classe
$\cP$ de fonctions sur $\berv$ {\`a} valeurs r{\'e}elles, et un op{\'e}rateur
$\Delta$ d{\'e}fini sur $\cP$ et {\`a} valeurs dans les mesures sur $\berv$
qui joue le r{\^o}le analogue de l'op{\'e}rateur de Laplace sur la droite
projective complexe.

\smallskip

Apr{\`e}s ces pr{\'e}liminaires, rappelons bri{\`e}vement la d{\'e}finition de la
hauteur de Weil standard, ce qui permettra de motiver la d{\'e}finition de
hauteur ad{\'e}lique. Soit $\ov{K}$ une clot{\^u}re alg{\'e}brique de $K$.  La
\emph{hauteur de Weil} (ou \emph{hauteur na{\"\i}ve}) d'un sous-ensemble
fini $F$ de $\ov{K}$, et invariant par l'action du groupe de Galois
$\Gal(\ov{K} / K)$, est par d{\'e}finition\footnote{dans tout l'article
  $\log$ d{\'e}signe le logarithme n{\'e}p{\'e}rien}
\begin{equation}\label{hauteur naive K}
\hn(F) \= |F|^{-1} \sum_{\alpha \in F} \sum_{v \in M_K} \log^+
\|\alpha\|_v ~.
\end{equation}
Ici $|F|$ d{\'e}signe la cardinalit{\'e} de $F$ et $\log^+\|\cdot\|_v = \log
\max \{1, \| \cdot \|_v \}$.  On d{\'e}finit aussi la hauteur na{\"\i}ve d'un
{\'e}l{\'e}ment $\alpha$ de $\ov{K}$ par $\hn(\alpha) = \hn(F)$, o{\`u} $F$ est
l'orbite de $\alpha$ sous l'action du groupe de Galois $\Gal(\ov{K} /
K)$.  On {\'e}tend $\hn$ {\`a} $\P^1(\ov{K}) = \ov{K} \cup \{ \infty \}$ en
posant $\hn(\infty) = 0$.  Notons que par construction cette hauteur
est invariante sous l'action de $\Gal(\ov{K} / K)$.

L'{\'e}quation~\eqref{hauteur naive K} admet une interpr{\'e}tation
simple en terme de th{\'e}orie du potentiel. Pour expliquer cel{\`a},
notons $\diag_v = \{ (z, z)\ ; \ z \in \C_v \}$ la diagonale de $\C_v
\times \C_v$.  Lorsque $v$ est infinie, pour chaque paire de mesures
bor{\'e}liennes $\rho, \rho'$ support{\'e}es dans $\pv$, et telles que la
fonction $\log^+ \| z - w \|_v$ soit int{\'e}grable par rapport {\`a} $\rho
\otimes \rho'$ sur $\C_v \times \C_v \setminus \diag_v$, on pose
$$ \lpar \rho, \rho' \rpar_v \= - \int_{\C_v \times \C_v \setminus
\diag_v} \log \| z - w \|_v \ d\rho(z) \otimes d\rho'(w) ~.
$$
Pour simplifier les notations, pour tout ensemble fini $F$, on
notera $[F]$ la mesure de probabilit{\'e} atomique {\'e}quidistribu{\'e}e sur les
points de $F$.  Notons maintenant que $\lpar \lambda_v, \lambda_v
\rpar_v =0$ et que
\begin{equation}\label{cercle unite}
\lpar [\alpha], \lambda_v \rpar_v = - \log^+ \| \alpha \|_v ~,
\end{equation}
pour tout $\alpha \in \C_v$, o{\`u} $\lambda_v$ est la mesure de
probabilit{\'e} proportionnelle {\`a} la mesure de Lebesgue sur le cercle
unit{\'e} $S^1 \subset \C_v$. Lorsque $v$ est finie, un accouplement
$\lpar \cdot, \cdot \rpar_v$ peut {\^e}tre d{\'e}fini de fa{\c c}on analogue.  Dans
ce cas, la formule~\eqref{cercle unite} reste valide pour tout $\alpha
\in \C_v$ si l'on remplace $\lambda_v$ par la masse de Dirac situ{\'e}e au
point de $\berv$ associ{\'e} {\`a} $B(0,1)$.  On v{\'e}rifie que $\lpar
\lambda_v, \lambda_v \rpar =0$ pour tout $v$ dans $M_K$, et que la formule du
produit donne $\sum_{M_K} \lpar [F], [F] \rpar =0$.  De~\eqref{hauteur
  naive K} et~\eqref{cercle unite}, on tire enfin
\begin{equation}\label{e-local}
\hn(F) = \frac12 \sum_{v \in M_K} \lpar [F] -\lambda_v, [F] -
\lambda_v \rpar_v ~,
\end{equation}
pour tout sous-ensemble fini $F$ de $\mathbb{P}^1(\ov{K})$ invariant
par l'action du groupe de Galois $\Gal(\ov{K}/K)$.

Nous proposons une d{\'e}finition de hauteurs bas{\'e}e sur cette
{\'e}galit{\'e}.
\begin{definition}
  Une \emph{mesure ad{\'e}lique} $\rho$ est la donn{\'e}e en chaque place
  d'une mesure de probabilit{\'e} $\rho_v$ support{\'e}e dans $\berv$, telle
  que $\rho_v = \lambda_v$ hors d'un nombre fini de places, et telle
  que $\rho_v$ admette un potentiel continu en toutes les places, \ie
  $\rho_v = \lambda_v + \Delta g$ avec $g$ continue.
\end{definition}
\begin{definition}
  La \emph{hauteur ad{\'e}lique} $h_\rho$ associ{\'e}e {\`a} la mesure
  ad{\'e}lique $\rho$ est par d{\'e}finition donn{\'e}e par
\begin{equation}\label{e-hauteur adel}
h_\rho(F) \= \frac12 \sum_{v \in M_K} \lpar [F] -\rho_v, [F] -
\rho_v \rpar_v ~,
\end{equation}
pour tout ensemble fini $F\subset\ov{K}$ invariant par
$\Gal(\ov{K}/K)$.  Pour tout $\a\in\mathbb{P}^1(\ov{K})$, on posera
$h_\rho(\alpha) = h_\rho(F)$, o{\`u} $F$ est l'orbite de $\alpha$ sous l'action
du groupe de Galois $\Gal(\ov{K} / K)$. 
\end{definition}
Les hauteurs ad{\'e}liques peuvent toutes {\^e}tre d{\'e}finies de mani{\`e}re
{\'e}quivalente par une formule du type Mahler.
\begin{proposition}\label{P-standard}
Soit $\alpha \in \ov{K}$ et soit $P \in K[T]$ le polyn{\^o}me (unitaire)
minimal de $\alpha$ sur $K$.  Alors on~a
$$
h_\rho(\alpha) = h_\rho(\infty) + \frac1{\deg(\alpha)}\sum_{v \in
  M_K} \int_{\berv} \log \| P \|_v \ d\rho_{v} ~.
$$
\end{proposition}
Aux places finies telles que $\rho_v = \lambda_v$, l'int{\'e}grale
$\int_{\berv} \log \| P \|_v \ d\rho_{v}$ est {\'e}gale {\`a} la norme de
Gauss du polyn{\^o}me $P= \sum a_i T^i$ \ie $\max \{|a_i|_v\}$. Ceci
montre que le membre de droite est en r{\'e}alit{\'e} une somme finie. Dans le
cas de la hauteur na{\"\i}ve, l'{\'e}galit{\'e} ci-dessus se r{\'e}duit {\`a} la formule de
Mahler habituelle.

%%%%%%%%%%%%%%%%%%%%%%%%%%%%%%%%%%%%%%%%%%%%%%%%%%%%%%%%%%%%%%%%%%%%%%%

\subsection{R{\'e}sultats principaux}
Les {\'e}nonc{\'e}s ci-dessous r{\'e}sument les propri{\'e}t{\'e}s g{\'e}n{\'e}\-rales
des hauteurs ad{\'e}liques.
\begin{theoreme}\label{T-main1}
  Pour toute mesure ad{\'e}lique $\rho$, la fonction $h_\rho$ est une
  hauteur de Weil dont le minimum essentiel est non-n{\'e}gatif.  
  
  En d'autres termes, la diff{\'e}rence $h_\rho -\hn$ est uniform{\'e}ment
  born{\'e}e sur $\mathbb{P}^1(\ov{K})$, et pour tout $\e>0$, l'ensemble
  $\{ \a \in \P^1(\ov{K});\, h_\rho(\a) < - \varepsilon \}$ est fini.
\end{theoreme}
\begin{theoreme}\label{T-main2}  
  Soit $\{ F_n \}_{n \ge 0}$ une suite d'ensembles finis distincts
  deux {\`a} deux et $\Gal(\ov{K}/K)$-invariants telle que $ \lim_{n \to
    \infty} h_\rho(F_n) = 0 $. Alors pour toute place $v$ de $M_K$ on
  a convergence faible au sens des mesures $[F_n] \to \rho_v$ lorsque
  $n \to \infty$.
\end{theoreme}
%Il nous semble int{\'e}ressant de donner une interpr{\'e}tation en terme
%g{\'e}om{\'e}trique de ces deux r{\'e}sultats. Soit tout d'abord $z \in
%\mathbb{P}^1(\ov{K})$. Pour toute place $v\in M_K$, notons $Z_v$ la
%mesure de probabilit{\'e} sur $\C_v$ {\'e}quidistribu{\'e}e sur les
%conjugu{\'e}s de $z$.  On dira qu'une mesure ad{\'e}lique $\rho$ est
%\emph{globale} si et seulement si il existe une suite de points alg{\'e}briques $z_n$
%tels que pour toute place $v\in M_K$, $Z_{n,v} \to \rho_v$.
%L'estimation sur le minimum essentiel et le th{\'e}or{\`e}me
%d'{\'e}quidistribution ci-dessus  impliquent alors l'{\'e}nonc{\'e} suivant.
%\begin{corollaire}
%  Une hauteur ad{\'e}lique a un minimum essentiel nul si et seulement si elle provient
%  d'une mesure ad{\'e}lique globale.
%\end{corollaire}
Nous donnerons aussi des estimations quantitatives pr{\'e}cises de la
vitesse de convergence $[F_n]\to\rho_v$ du Th{\'e}or{\`e}me~\ref{T-main2}.
Afin d'{\'e}viter d'introduire trop de terminologie, nous ne mentionnons
ici qu'un {\'e}nonc{\'e} aux places infinies, et nous renvoyons au
Th{\'e}or{\`e}me~\ref{T-estimprecise} en \S\ref{S-vitesse} pour un {\'e}nonc{\'e}
analogue aux places finies.
\begin{theoreme}\label{T-estim}
  Soit $\rho$ une mesure ad{\'e}lique admettant un potentiel H{\"o}lder en
  toutes les places.  Alors il existe une constante $C>0$ telle que,
  pour tout ensemble fini $F \subset \ov{K}$, invariant par l'action
  du groupe de Galois et de cardinalit{\'e} $|F|$, pour toute place
  infinie $v \in M_K$, et pour toute fonction $\varphi$ de classe
  $\cC^1$ sur $\pv$, on~a
  $$
  \left|\frac1{|F|} \sum_{\a \in F} \varphi (\a) - \int \varphi \,
    d\rho_{v} \right| \le \left( h_{\rho}(F) + C \, \frac{\log
      |F|}{|F|} \right) \times \lip (\varphi) ~,$$
  o{\`u} $\lip(\varphi) =
  \sup_{x\neq y} |\varphi(x)-\varphi(y)|/\mathsf{d}(x,y)$ et
  $\mathsf{d}$ est la m{\'e}trique sph{\'e}rique sur $\pv$.
\end{theoreme}
Afin de comprendre la force de cette {\'e}nonc{\'e}, mentionnons le
corollaire nouveau suivant dans le cas de la hauteur na{\"\i}ve.  Ici
$\lambda_{S^1}$ est la mesure de Haar sur le cercle unit{\'e}.
\begin{corollaire} Il existe une constante $C>0$ telle que
  pour tout ensemble fini $F \subset \ov{\Q}$, invariant par l'action
  du groupe de Galois et de cardinalit{\'e} $|F|$, et pour toute fonction
  $\varphi$ de classe $\cC^1$ sur $\mathbb{P}^1(\C)$, on~a
  $$
  \left|\frac1{|F|} \sum_{\a \in F} \varphi (\a) - \int \varphi \,
    d\lambda_{S^1} \right| \le \left( \hn (F) + C\, \frac{\log
      |F|}{|F|} \right) \times \lip (\varphi) ~.$$
\end{corollaire}
Dans un travail r{\'e}cent, Petsche~\cite{petsche} a obtenu une
estimation moins forte, mais pour une classe de fonctions plus
g{\'e}n{\'e}rales.

%%%%%%%%%%%%%%%%%%%%%%%%%%%%%%%%%%%%%%%%%%%%%%%%%%%%%%%%%%%%%%%%%%%
\subsection{Exemples et applications dynamiques}
Les hauteurs ad{\'e}liques que nous avons construites recouvrent
essentiellement toutes les constructions de hauteurs d{\'e}j{\`a}
consid{\'e}r{\'e}es dans la lit\-t{\'e}rature (sur la droite projective).
Elles co{\"\i}ncident avec les hauteurs issues de m{\'e}triques (sur le
fibr{\'e} $O(1)$) dites \emph{ad{\'e}liques int{\'e}grables}, dont la
construction a {\'e}t{\'e} r{\'e}alis{\'e}e par Zhang et {\'e}tendue par
Chambert-Loir~\cite{ACL2} dans un travail ind{\'e}pendant du n{\^o}tre.  On
voit ainsi que notre Th{\'e}or{\`e}me~\ref{T-main2} est {\'e}quivalent
{\`a}~\cite[Th{\'e}or{\`e}me~4.2]{ACL2} dans le cas de la droite projective.
Nous nous contenterons ici de d{\'e}crire deux types de situations dans
lesquelles les Th{\'e}or{\`e}\-mes~\ref{T-main1},~\ref{T-main2},
et~\ref{T-estim} permettent d'{\'e}tendre et de pr{\'e}ciser certains
r{\'e}sultats d{\'e}j{\`a} connus.

\smallskip

Supposons donn{\'e} en chaque place infinie, un compact du plan complexe
$E_v$ de capacit{\'e} logarithmique strictement positive. Pour toute place
finie, posons $\rho_v = \lambda_v$, et pour toute place infinie notons
$\rho_v$ la mesure d'{\'e}quilibre (au sens de la th{\'e}orie du potentiel) de
$E_v$.  Toutes ces mesures sont {\`a} potentiel localement born{\'e}e et sous
l'hypoth{\`e}se suppl{\'e}mentaire que ces potentiels sont continus, on peut
donc leur associer une hauteur ad{\'e}lique $h_\rho$. Dans ce cadre, le
Th{\'e}or{\`e}me~\ref{T-main1} implique la partie ais{\'e}e du th{\'e}or{\`e}me de
Fekete-Szeg{\"o} et montre la finitude du nombre de points entiers dont
tous les conjugu{\'e}s sont dans un voisinage fixe de $E_v$ (en toutes les
places infinies) sous une hypoth{\`e}se convenable sur la capacit{\'e} des
$E_v$.  Le Th{\'e}or{\`e}me~\ref{T-estim} quant {\`a} lui nous donne une version
quantitative de~\cite[Th{\'e}or{\`e}me~1]{Rumely-bilu}.

\smallskip

%Une deuxi{\`e}me cat{\'e}gorie de hauteur est aussi de nature ad{\'e}lique.
%Supposons donn{\'e} un morphisme alg{\'e}brique $\phi$ d{\'e}fini sur $K$ de
%$\mathbb{P}^1(K)$ dans $\mathbb{P}^n(K)$, pour $n\ge1$. On peut alors
%munir $\mathbb{P}^n(\ov{K})$ de la hauteur de Weil standard $\hn$ et
%la ramener sur $\mathbb{P}^1(\ov{K})$ par $\phi$. Il est facile de
%voir que cette hauteur est la somme d'une constante et d'une hauteur
%de type ad{\'e}lique. Le Th{\'e}or{\`e}me~\ref{T-main1} donne dans ce cas une
%estimation sur le minimum essentiel de $\hn\circ \phi$ sur
%$\mathbb{P}^1(\ov{K})$. Cette estimation peut {\^e}tre {\`a} nouveau utilis{\'e}
%pour montrer la finitude des points de petites hauteurs sur
%$\mathbb{P}^1(\ov{K})$.  En appliquant le Th{\'e}or{\`e}me~\ref{T-estim} {\`a}
%l'application $t\in\ov{\Q} \mapsto (t, t-1)\in\ov{\Q}^2$, on montre
%ainsi que le minimum essentiel de la hauteur de Weil sur la courbe
%$X+Y =1 $ dans le plan projectif est minor{\'e} par $\int_0^{\pi/3} \log
%|\sin \theta|^{-1} \, \frac{d\theta}{2\pi}$ qui est strictement
%positif. Ceci donne une version faible de~\cite[Th{\'e}or{\`e}me~1]{Za}.
%Il serait int{\'e}ressant de voir si le Th{\'e}or{\`e}me~\ref{T-estim}
%permet d'obtenir des estimations de la hauteur pour des courbes plus
%g{\'e}n{\'e}rales.
%
%
%\smallskip

De notre point de vue cependant, les applications les plus
significatives concer\-nent une classe de hauteurs issues des syst{\`e}mes
dynamiques.  Soit donc $R$ une fraction rationnelle {\`a} coefficients
dans un corps de nombres $K$ et de degr{\'e} $D \ge 2$. On peut montrer
que la limite $\lim_{n\to\infty} D^{-n} \hn \circ R^n$ existe et
d{\'e}finit une hauteur de Weil $h_R$ qui v{\'e}rifie $h_R
\circ R = D\, h_R$. Nous allons voir que $h_R$ est une hauteur
ad{\'e}lique, mais pour ce faire, il nous faut tout d'abord d{\'e}crire
quelques r{\'e}sultats de nature dynamique.

En toute place infinie $v\in M_K$, les it{\'e}r{\'e}s $\{R^n\}_{n\ge0}$ de $R$
induisent un syst{\`e}me dynamique sur la sph{\`e}re de Riemann. Bien que la
nature des suites $\{R^n(z)\}_{n\ge0}$ d{\'e}pendent tr{\`e}s fortement du
choix du point $z$, l'action par images inverses de $R$ pr{\'e}sente des
propri{\'e}t{\'e}s d'uniformit{\'e} tout {\`a} fait remarquable. On d{\'e}montre
en effet qu'il existe une mesure de probabilit{\'e} $\rho_{R,v}$, dite
\emph{mesure d'{\'e}quilibre}, telle que pour tout $z_0 \in \pc$ non
exceptionnel pour $R$, on~a
\begin{equation}\label{E-equicomp}
\lim_{n\to \infty} D^{-n} [R^{-n}\{z_0\}] = \rho_{R,v}~.
\end{equation}
Rappelons qu'un point est dit {\it exceptionnel} pour $R$ si son
orbite inverse est finie, et qu'une fraction rationnelle donn{\'e}e admet
au plus deux points exceptionnels. Ce r{\'e}sultat est d{\^u} {\`a}
Brolin~\cite{Brolin} dans le cas des polyn{\^o}mes, et ind{\'e}pendemment {\`a}
Lyubich~\cite{Lyubich} et {\`a} Freire-Lopez-Ma{\~n}{\'e}~\cite{FLM}, dans le
cas des fractions rationnelles.  La mesure $\rho_{R,v}$ est support{\'e}e
sur son ensemble de Julia\footnote{la partie de la sph{\`e}re de Riemann o{\`u}
  la dynamique est chaotique.}, et permet d'obtenir de nombreuses
informations sur le syst{\`e}me dynamique induit par $R$. C'est de plus
une mesure {\`a} potentiel continue (et m{\^e}me H{\"o}lder), voir~\S~\ref{S-holder}.

En toute place finie $v$, une mesure $\rho_{R,v}$ satisfaisant {\`a} une
propri{\'e}t{\'e} analogue a {\'e}t{\'e} construite dans~\cite{FR}.  Celle-ci n'est
pas support{\'e}e en g{\'e}n{\'e}ral dans l'espace projectif standard $\pv$, mais
dans l'espace de Berkovich associ{\'e}, et elle est aussi {\`a} potentiel
continu.

On peut maintenant {\'e}noncer le
\begin{theoreme}\label{T-dyn}
  Soit $R$ une fraction rationnelle de degr{\'e} au moins~$2$ et {\`a}
  coefficients dans un corps de nombre $K$.  Pour chaque place $v$ de
  $K$ soit $\rho_{R, v}$ la mesure d'{\'e}quilibre de $R$ correspondante.
  Alors $\rho_R = \{ \rho_{R, v} \}_{v \in M_K}$ est une mesure
  ad{\'e}lique et la hauteur normalis{\'e}e $h_R$ co{\"\i}ncide avec la hauteur
  $h_{\rho_R}$ d{\'e}finie par la mesure ad{\'e}lique $\rho_R$.
\end{theoreme}
La Proposition~\ref{P-standard} s'applique donc, et nous donne ainsi
une formule de Mahler pour toutes ces hauteurs. Dans ce cadre,
celle-ci avait {\'e}t{\'e} {\'e}nonc{\'e}e et d{\'e}montr{\'e}e dans~\cite{PST}.

La hauteur $h_R$ est aussi redevable des
Th{\'e}or{\`e}mes~\ref{T-main1},~\ref{T-main2} et~\ref{T-estim} ci-dessus.
Notons cependant que le Th{\'e}or{\`e}me~\ref{T-main1} ne nous donne aucune
information.  En effet, on~a par construction $h_R\ge0$, et le minimum
essentiel de $h_R$ est nul, car $\{ h_R =0 \}$ qui est constitu{\'e} des
points pr{\'e}p{\'e}riodiques de $R$, est toujours infini.  Par contre, le
Th{\'e}or{\`e}me~\ref{T-main2} permet d'obtenir des r{\'e}sultats
d'{\'e}quidistribution remarquables.  On obtient ainsi le corollaire
suivant, qui donne une preuve arith\-m{\'e}\-tique d'un r{\'e}sultat d{\^u}
{\`a}~\cite{Lyubich} dans le cas complexe.
\begin{corollaire}\label{C-periodique}
  Soient $R$ et $S$ deux fractions rationnelles {\`a} coefficients dans un
  corps de nombres $K$, avec $\deg(R)\ge2$.  Pour $n\ge1$, notons
  $F_n$ le sous-ensemble des solutions dans $\ov{K}$ {\`a} l'{\'e}quation $R^n
  = S$.  Si pour $n$ grand les ensembles $F_n$ sont distincts deux {\`a}
  deux, alors pour toute place $v\in M_K$, on~a
$$
\lim_{n\to \infty} [F_n] = \rho_{R,v}~.
$$
\end{corollaire}
Lorsque $S(z) = z$, l'ensemble $F_n$ est {\'e}gal {\`a} l'ensemble des points
p{\'e}riodiques de $R$ dans $\ov{K}$, de p{\'e}riode $n$.  Dans ce cas les
ensembles $F_n$ sont distincts deux {\`a} deux (pour $n$ grand) et on
obtient l'{\'e}quidistribution des points p{\'e}riodiques de $R$ selon la
mesure d'{\'e}quilibre.  C'est un r{\'e}sultat nouveau pour toute place finie.

Lorsque $z_0 \in K$ n'est pas exceptionnel pour $R$ et lorsque la
fraction rationnelle $S$ est constante {\'e}gale {\`a} $z_0$, on sait que les
ensembles $F_n$ sont distincts deux {\`a} deux.  On obtient alors le
r{\'e}sultat d'{\'e}quidistribution des pr{\'e}images it{\'e}r{\'e}es de $z_0$, mentionn{\'e}
ci-dessus en~\eqref{E-equicomp}.

C'est dans ce contexte que plusieurs cas particuliers du
Corollaire~\ref{C-periodique} ont {\'e}t{\'e} obtenus pr{\'e}c{\'e}demment.
Bilu~\cite{bilu} l'a tout d'abord d{\'e}montr{\'e} pour les morphismes de
puissance.  Autissier~\cite{Autissier1} a ensuite obtenu ce r{\'e}sultat {\`a}
la place infinie pour toutes les fractions rationnelles comme un cas
particulier d'un r{\'e}sultat concernant les courbes arithm{\'e}tiques.  Le
Corollaire~\ref{C-periodique} a r{\'e}cemment {\'e}t{\'e} d{\'e}montr{\'e} pour les
polyn{\^o}mes par Baker-Hsia~\cite{Baker-Hsia} {\`a} la place infinie et, sous
certaines hypoth{\`e}ses, aux places finies.

Finalement, dans des travaux ind{\'e}pen\-dants des n{\^o}tres,
Baker-Rumely~\cite{BR} d'une part et Chambert-Loir~\cite{ACL2} d'autre
part ont d{\'e}\-montr{\'e} le Corollaire~\ref{C-periodique} pour $S(z) =z$ en
toutes les places (finies et infinies).

\medskip La mesure ad{\'e}lique $\rho_R$ est {\`a} potentiel H{\"o}lder
(voir~\S~\ref{S-holder}) et par cons{\'e}quent la hauteur normalis{\'e}e $h_R$
v{\'e}rifie toutes les conditions du Th{\'e}or{\`e}me~\ref{T-estim}.  Ceci permet
d'obtenir imm{\'e}diatement la version quantitative suivante du
Corollaire~\ref{C-periodique}.
\begin{corollaire}
  Soit $R$ une fraction rationnelle de degr{\'e} $D\ge2$ a coefficients
  dans $\ov{\Q}$, et notons $\rho_R$ sa mesure d'{\'e}quilibre sur
  $\mathbb{P}^1(\C)$. Alors il existe une constante $C>0$ telle que
  pour toute fonction $\varphi$ de classe $\cC^1$ sur
  $\mathbb{P}^1(\C)$, pour tout point $z\in\mathbb{P}^1(\ov{\Q})$ non
  exceptionnel, et pour tout $n\ge0$, on ait
  $$
  \left|\frac1{D^n} \sum_{\a \in R^{-n}\{z\}} \varphi (\a) - \int
\varphi \, d\rho_{R} \right|
\le
\left( \frac{h_R(z) + C n}{D^n} \right) \times \lip (\varphi) ~.$$
\end{corollaire}
Ce r{\'e}sultat est tout {\`a} fait surprenant dans la mesure o{\`u} pour une
fraction rationnelle {\`a} coefficients complexes quelconques,
l'estimation en $nD^{-n}$ ne semble pas connue.  Des estimations en
$\exp(-\sqrt{n})$ ont {\'e}t{\'e} obtenues dans~\cite{denker} (et pour des
fonctions $f$ H{\"o}lder), et raffin{\'e}es en $\sigma^n$ pour un $\sigma
<1$ proche $1$ dans~\cite{haydn}, voir aussi~\cite{polsharp} pour des
r{\'e}sultats plus faibles.

\smallskip

Il nous semble int{\'e}ressant de mentionner aussi le corollaire suivant
dont une preuve directe par des m{\'e}thodes complexes semble d{\'e}licate. Ce
corollaire nous a {\'e}t{\'e} inspir{\'e} par~\cite[Theorem~8.13]{Baker-Hsia} dont
le th{\'e}or{\`e}me ci-dessous en est une version quantitative.  Fixons un
entier $D\ge2$, et regardons l'ensemble des polyn{\^o}mes de la forme
$P_c(z) = z^D + c$ pour $c \in \C$. Il est int{\'e}ressant de regarder
l'ensemble dit de Mandelbrot et not{\'e} $\cM_D$, constitu{\'e} des param{\`e}tres
$c$ pour lesquels l'orbite de $z = 0$ pour $P_c$ est born{\'e}e. On montre
que $\cM_D$ est un ensemble compact. On peut donc consid{\'e}rer sa mesure
harmonique $\mu_D$, qui est caract{\'e}ris{\'e}e de mani{\`e}re dynamique par la
formule $\mu_D = \lim_{n\to\infty} D^{-n} \Delta \log^+|P^n_c(0)|$.

On dit qu'un param{\`e}tre $c \in \C$ est {\it critiquement fini}, s'il
existe des entiers distincts $n$ et $m$ tels que $P_c^n(0) =
P_c^m(0)$. Notons que pour de tel param\`etres, $c \in \bar{\Q}$.
\begin{theoreme}\label{thm:mandel}
  Il existe une constante $C > 0$ telle que pour tout ensemble fini et
  $\Gal(\bar{\Q}/\Q)$-invariant $F \subset \C$ de param{\`e}tres
  critiquement finis et toute fonction $\varphi$ de classe $\cC^1$ sur
  $\P^1(\C)$, on~a
  $$
  \left| |F|^{-1} \sum_{\alpha \in F} \varphi(\alpha) - \int
    \varphi \, d\mu_D \right| \le \mathrm{C}\, \frac{\log |F|}{|F|}
  \times \lip (\varphi) ~.$$
  En particulier pour toute suite
  d'ensembles finis $\{ F_n \}_{n \ge 1}$ distincts deux {\`a} deux
  v{\'e}rifiant les propri{\'e}t{\'e}s ci-dessus, on~a $[F_n]\to\mu_D$.
\end{theoreme}
La preuve du Th{\'e}or{\`e}me~\ref{thm:mandel} est donn{\'e}
en~\S~\ref{S-mandel}.  Nous indiquons maintenant rapidement la preuve
du Corollaire~\ref{C-periodique}.
\begin{proof}[D{\'e}monstration du Corollaire~\ref{C-periodique}]
  On fixe tout d'abord des constantes $B,C>0$ telles que $h_R(S(z))
  \le B \cdot h_R(z) + C$ (on peut en fait prendre $B = \deg (S)$).
  Pour tout $z\in F_n$, on~a alors $D^n h_R(z) = h_R(R^n(z)) =
  h_R(S(z)) \le B \cdot h_R(z) + C$, donc
$$
h_R(z) \le \frac{C}{D^n - B} \text{ pour tout } z \in F_n~.
$$ Comme $R$ et $S$ sont {\`a} coefficients dans $K$, les ensembles
$F_n$ sont $\Gal(\ov{K}/K)$-invariant, et l'estimation pr{\'e}c{\'e}dente
donne $\lim_{n\to\infty} h_R(F_n) =0$. Sous l'hypoth{\`e}se suppl{\'e}mentaire
que les ensembles $F_n$ sont distincts deux {\`a} deux, le
Th{\'e}or{\`e}me~\ref{T-main2} s'applique et donne alors $\lim_{n \to \infty} [F_n] =
\rho_{R,v}$.
\end{proof}

%%%%%%%%%%%%%%%%%%%%%%%%%%%%%%%%%%%%%%%%%%%%%%%%%%%%%%%%%%%%%%%%%%

\subsection{Strat{\'e}gie de la preuve.}\label{S-esquisse}
Nous indiquons maintenant comment d{\'e}montrer nos r{\'e}sultats
principaux, les Th{\'e}or{\`e}mes~\ref{T-main1} et~\ref{T-main2}.

Le fait que pour toute mesure ad{\'e}lique $\rho$, la fonction $h_\rho$
soit une hauteur de Weil r{\'e}sulte directement de notre hypoth{\`e}se de
continuit{\'e} sur les potentiels des mesures $\rho_v$. Les deux autres
{\'e}nonc{\'e}s, le fait que le minimum essentiel de $h_\rho$ soit positif et
le r{\'e}sultat d'{\'e}quidistribution, sont en r{\'e}alit{\'e} la cons{\'e}quence d'une
m{\^e}me estimation de positivit{\'e} de chacun des termes $\lpar [F] -\rho_v,
[F] - \rho_v \rpar_v$ intervenant dans la d{\'e}finition de $h_\rho$, que
nous expliquons dans le cas des places infinies.

Etant diff{\'e}rence de deux mesures de probabilit{\'e}, chaque mesure
$[F]-\rho_v$ s'{\'e}crit $\Delta g$ pour une fonction $g$ d{\'e}finie
globalement sur la sph{\`e}re de Riemann. Une int{\'e}gra\-tion par partie
montre alors que\footnote{{\`a} une constante multiplicative pr{\`e}s tenant
  compte du fait que $v$ est r{\'e}el ou non.} $\lpar [F] -\rho_v, [F] -
\rho_v \rpar_v = \int_\pv dg \wedge d^cg $ d{\`e}s que cette int{\'e}grale est
bien d{\'e}finie.  C'est le cas lorsque $g$ est lisse, ce que nous allons
supposer un instant pour les besoins de la discussion.  Dans ce cas
$\int dg \wedge d^cg = \int |\partial g/ \partial x|^2 + |\partial
g/\partial y|^2 \, dxdy$ est positif et s'annule si et seulement si
$g$ est constante, ou bien de mani{\`e}re {\'e}quivalente si et seulement si
$[F]-\rho_v =\Delta g =0$.  On expliquera en \S\ref{S-positive} que
tout ceci reste vrai sous l'hypoth{\`e}se plus faible que $g$ est
continue.

Cependant la mesure $[F]$ est atomique et donc $g$ n'est m{\^e}me pas
localement born{\'e}e. L'id{\'e}e consiste alors {\`a} approcher $[F]$ par
une famille de mesures lisses $[F]_\e$ (en convolant par un noyau
lisse), et le point cl{\'e} est d'estimer pr{\'e}cis{\'e}ment la diff{\'e}rence $\lpar
[F] -\rho_v, [F] - \rho_v \rpar_v - \lpar [F]_\e -\rho_v, [F]_\e -
\rho_v \rpar_v$.  C'est le contenu des Lemmes~\ref{L-100}
et~\ref{L-120}, qui permettent de contr{\^o}ler cette diff{\'e}rence en termes
du param{\`e}tre $\e$ et de la cardinalit{\'e} de $F$.

La m{\^e}me analyse est possible aux places finies, si l'on remplace
l'espace projectif standard par la droite projective au sens de
Berkovich. On utilise dans ce cas un proc{\'e}d{\'e} de r{\'e}gularisation de
nature {\'e}l{\'e}mentaire (bas{\'e} sur la structure d'arbre de $\berv$) pour
estimer la positivit{\'e} des termes $\lpar [F] -\rho_v, [F] -
\rho_v\rpar$. Ceci aboutit aux Lemmes~\ref{L-300} et~\ref{L-320}.

Une fois ces estimations faites, la preuve du
Th{\'e}or{\`e}me~\ref{T-main1} est une application simple du th{\'e}or{\`e}me de
Northcott sur la finitude des points de hauteur (na{\"\i}ve) et de
degr{\'e} born{\'e}s.  Concernant le Th{\'e}or{\`e}me~\ref{T-main2}, si $F_n$ est une
suite d'ensembles finis telle que $h_\rho(F_n)\to 0$, nos estimations de
positivit{\'e} et le fait que $|F_n| \to \infty$ impliquent en toutes les
places $ \lim_{n \to 0} \lpar [F_n]- \rho_v , [F_n] - \rho_v \rpar_v
=0 $, dont on d{\'e}duit  $[F_n]\to \rho_v$.

Notons finalement que les estimations de positivit{\'e} que nous donnons
sont essentiellement optimales.  Nous d{\'e}crivons en~\S\ref{S-exemple}
un exemple qui montre qu'en tous les cas ces estimations sont
n{\'e}cessaires. On construit en effet une hauteur ad{\'e}lique (de type
dynamique) $h_\rho$, une suite d'ensemble fini $F_n$ invariant sous
$\Gal(\ov{K}/K)$ et de cardinalit{\'e} tendant vers l'infini, pour
lesquels il existe une place (finie) telle que $\lpar [F_n]-\rho_v ,
[F_n]-\rho_v \rpar <0$ pour tout $n$.

%
%%%%%%%%%%%%%%%%%%%%%%%%%%%%%%%%%%%%%%%%%%%%%%%%%%%%%%%%%%%%%%%%%%
%
\subsection{Plan de l'article.}
Cet article est divis{\'e} en cinq parties.  Dans la
Section~\ref{S-pot_complexe}, nous rappelons les {\'e}l{\'e}ments de th{\'e}orie
du potentiel sur $\C$ n{\'e}cessaires {\`a} notre analyse.  Bien que le
contenu de cette partie soit classique, le traitement que nous donnons
ici est adapt{\'e}e pr{\'e}cis{\'e}ment {\`a} nos besoins et sert de base au
traitement de la th{\'e}orie du potentiel que nous d{\'e}veloppons aux places
finies par la suite.  Dans la Section~\ref{S-berkovich}, nous faisons
quelques rappels sur la g{\'e}om{\'e}trie de la droite projective au sens de
Berkovich sur $\C_p$.  Dans la Section~\ref{S-energieCp}, nous
d{\'e}crivons les r{\'e}sultats de th{\'e}orie du potentiel sur $\C_p$ analogues {\`a}
ceux de la Section~\ref{S-pot_complexe}.  La Section~\ref{S-hauteur
  Arakelov} est d{\'e}di{\'e}e aux preuves des
Th{\'e}or{\`e}mes~\ref{T-main1},~\ref{T-main2}, et~\ref{T-estim}. Enfin nous
montrons dans la derni{\`e}re Section~\ref{S-dyn} le Th{\'e}or{\`e}me~\ref{T-dyn}
{\'e}tablissant que les hauteurs dynamiques sont des hauteurs ad{\'e}liques,
ainsi que le Th{\'e}or{\`e}me~\ref{thm:mandel}.

\newpage 
%
%
%%%%%%%%%%%%%%%%%%%%%%%%%%%%%%%%%%%%%%%%%%%%%%%%%%%%%%%%%%%%%%%%%%%%%%
%
%

\section{Th{\'e}orie du potentiel dans le cas complexe}\label{S-pot_complexe}
%
%%%%%%%%%%%%%%%%%%%%%%%%%%%%%%%%%%%%%%%%%%%%%%%%%%%%%%%%%%%%%%%%%%%%%
%
%
\subsection{Forme de Dirichlet}\label{S-dirichlet}
On identifie $\C$ {\`a} $\R \oplus i \R$ et pour $z = x + i y \in \C$ on
pose $\bz = x - iy$.  Toute application lin{\'e}aire r{\'e}elle d{\'e}finie
sur $\C$ et {\`a} valeurs dans $\C$, s'{\'e}crit de fa{\c c}on unique
comme une combinaison lin{\'e}aire complexe des applications $z \mapsto
z$ et $z \mapsto \bz$.  Toute $1$-forme $\omega$ sur $\C$ se
d{\'e}compose donc uniquement sous la forme $\omega = \omega_{1,0} dz +
\omega_{0,1} d\bz$; la forme $\omega_{1,0} dz$ (resp.
$\omega_{0,1} d\bz$) est la composante dite de type $(1,0)$
(resp. $(0,1)$) de la forme $\omega$. Si $f: \C \to \C$ est une
fonction de classe $\cC^1$, on note $df = \partial f + \bpartial f$,
o{\`u} $\partial f$ est la composante de type $(1,0)$ et $\bpartial f$
celle de type $(0,1)$ de $df$. On pose $d^c f \= \frac 1{2\pi i}
(\partial f - \bpartial f)$. C'est un op{\'e}rateur r{\'e}el au sens o{\`u}
$\overline{d^c f} = d ^c \overline{f}$. Si $f$ et $g$ sont deux
fonctions de classe $\cC^1$ {\`a} valeurs \emph{r{\'e}elles}, on v{\'e}rifie
que
\begin{equation}\label{e-dirichlet}
df \wedge d^c g = d g \wedge d^c f = \left( \frac{\partial f}{\partial x}
\frac{\partial g}{\partial x}+\frac{\partial f}{\partial
y}\frac{\partial g}{\partial y}\right)~\frac{dx\wedge dy}{2\pi}~.
\end{equation}
Fixons $f,g$ deux fonctions r{\'e}elles de classe $\cC^1$. Pour
tout ouvert connexe $D \subset \C$, on notera
$\langle f, g \rangle _D \= \int_D df \wedge d^c g \in \R$.
Cet accouplement est appel{\'e} classiquement \emph{forme de Dirichlet}.
Il est clair que $\langle \cdot , \cdot \rangle_D$ d{\'e}finit une forme
bilin{\'e}aire positive sur l'espace des fonctions de classe $\cC^2$ et que
l'on a $ \langle f , g \rangle_D^2 \le \langle f , f \rangle_D \cdot
\langle g , g \rangle_D $ avec {\'e}galit{\'e} si et seulement si il
existe une constante $c \in \R$ telle que la fonction $f-c g$ soit
constante sur $D$.

Enfin, on notera que pour toute fonction de classe $\cC^2$, on~a
$dd^cf =(\Delta f) dx \wedge dy$, o{\`u} $\Delta f = (2\pi)^{-1} \left(
\frac{\partial^2 f}{\partial x^2}+ \frac{\partial^2 f}{\partial y^2}
\right)$ est le Laplacien standard de $f$ sur $\C$. On verra {\`a} la section
suivante la justification du choix de normalisation par $2\pi$.  Enfin
on fera souvent l'abus de notations $\Delta f = dd^cf$.

%
%%%%%%%%%%%%%%%%%%%%%%%%%%%%%%%%%%%%%%%%%%%%%%%%%%%%%%%%%%%%%%%%%%%%%
%
%
\subsection{Fonctions sous-harmoniques}\label{S-sh}
Dans la preuve du th{\'e}or{\`e}me principal, nous aurons besoin de
travailler avec la forme de Dirichlet appliqu{\'e}e {\`a} des fonctions de
classe de r{\'e}gularit{\'e} plus faible que $\cC^1$. Ceci est possible si
les fonctions poss{\'e}dent des propri{\'e}t{\'e}s de convexit{\'e} compatible
avec la structure complexe.  Commen{\c c}ons par une d{\'e}finition.
\begin{definition}
Une fonction $u: \C \to \R \cup \{ - \infty\}$ non identiquement
$-\infty$ est dite sous-harmonique si elle est semi-continue
sup{\'e}rieure, et qu'elle v{\'e}rifie en tout point et pour tout rayon
$r>0$, l'in{\'e}galit{\'e} dite de sous-moyenne $u(z)\le \int_{[0,2\pi]} u ( z
+ r e^{it}) \frac{dt}{2\pi}$.
\end{definition}
On v{\'e}rifie que toute fonction sous-harmonique est localement
int{\'e}grable. Gr{\^a}ce aux in{\'e}galit{\'e}s de sous-moyenne, on montre
que par convolution toute fonction sous-harmo\-nique est limite
\emph{d{\'e}croissante} de fonctions sous-harmoniques \emph{lisses}.  On
en d{\'e}duit alors que pour toute fonction sous-harmonique $u$, la
distribution $dd^c u $ est une mesure
\emph{positive}. R{\'e}ciproquement, on montre que toute fonction
$L^1_{\mathrm{loc}}$ dont le $dd^c$ au sens des distributions est une
mesure positive est {\'e}gale presque partout {\`a} une fonction
sous-harmonique.

Pour tout $z_0 \in \C$, la fonction $z \mapsto \log |z-z_0|$ est
sous-harmonique et on~a $dd^c \log|\cdot - z_0| = [z_0]$ (c'est pour
que cette formule soit valide que l'on a normalis{\'e} le laplacien en
divisant par $2\pi$) . Par int{\'e}gration, on en d{\'e}duit que pour
toute mesure de probabilit{\'e} $\rho$ sur $\C$ et {\`a} support born{\'e},
la fonction $g_\rho(z) \= \int_\C \log|z-w|\, d\rho(w)$ est
sous-harmonique dans $\C$, et l'on a $\Delta g_\rho = \rho$.  Lorsque
$\rho \= \lambda_1$ est la mesure de probabilit{\'e} proportionnelle {\`a}
la mesure de Lebesgue sur le cercle unit{\'e}, on obtient $g_{\lambda_1} (z) =
\log^+|z|$.

Dans toute la suite, on d{\'e}signera par potentiel d'une mesure $\rho$
d{\'e}finie sur un domaine $\Omega$, toute fonction sous-harmonique $u$
d{\'e}finie sur $\Omega$ et telle que $\Delta u = \rho$. Un potentiel
est d{\'e}fini {\`a} une fonction harmonique pr{\`e}s. On dira que $\rho$
est {\`a} potentiel continu (resp. born{\'e}) si $\rho = \Delta u$ avec $u$
continue (resp. born{\'e}e). Ceci ne d{\'e}pend pas du choix du
potentiel car toute fonction harmonique est lisse.
%

%
%%%%%%%%%%%%%%%%%%%%%%%%%%%%%%%%%%%%%%%%%%%%%%%%%%%%%%%%%%%%%%%%%%%%%
%
%
\subsection{R{\'e}gularit{\'e} des potentiels.}

Comme la fonction $z \mapsto \log|z|$ est un noyau fondamental pour le
laplacien, pour toute fonction sous-harmonique $u$ d{\'e}finie dans le
disque unit{\'e}, la fonction $ u(z) - \int_{|z| \le 1} \log |z-w| \,
\Delta u (w)~, $ est harmonique, donc lisse. On d{\'e}duit de ce fait
qu'une fonction sous-harmonique est dans $L^p_{\mathrm{loc}}$ pour
tout $1\le p < \infty$ (car c'est le cas pour $\log|z|$); et que ses
d{\'e}riv{\'e}es partielles $\frac{\partial u}{\partial x}$ et $\frac{\partial
  u}{\partial y}$ sont dans $L^{2-\e}_{\mathrm{loc}}$ pour tout $\e>0$
(car c'est le cas pour la d{\'e}riv{\'e}e $| \cdot |^{-1}$ de la fonction
$\log| \cdot |$).  En particulier, $du$ est une $1$-forme {\`a}
coefficients $L^{2-\e}_{\mathrm{loc}}$.

Notons cependant que, {\'e}tant donn{\'e}e une partie ouverte et connexe $D$
de~$\C$, ces propri{\'e}t{\'e}s ne suffisent pas pour pouvoir d{\'e}finir $\langle
u, v \rangle_D$ pour un couple arbitraire de fonctions
sous-harmoniques $u$ et $v$.  La condition la plus faible sous
laquelle ce produit est d{\'e}fini et pour laquelle l'in{\'e}galit{\'e} de
Cauchy-Schwartz est v{\'e}rifi{\'e}e est naturellement $u \in L^1 ( dd^c v)$.
Cependant, dans ce cas la $2$-forme $du \wedge d^c v$ n'est que
mesurable en g{\'e}n{\'e}ral, ce qui complique nettement l'exposition.  Nous
nous contenterons donc de r{\'e}sultats plus faibles.

\begin{lemme}\label{L-intsh}
Soit $u$ une fonction sous-harmonique \emph{localement born{\'e}e}.
Alors la $1$-forme $du$ est {\`a}
coefficients $L^2_{\mathrm{loc}}$.
\end{lemme}
Il est en fait vrai que $u \in L^1_{\mathrm{loc}} ( dd^c u)$ 
{\'e}quivaut au fait que  $du$ est {\`a} coefficients $L^2_{\mathrm{loc}}$,
mais nous n'aurons pas besoin de ce r{\'e}sultat plus fort.
\medskip

\begin{proof}
On fixe donc une suite d{\'e}croissante de fonctions sous-harmoni\-ques
lis\-ses $u_n$ convergeant vers $u$.  On peut construire $u_n$ par
convolution, et on montre alors que $u_n \to u$ dans
$L^p_{\mathrm{loc}}$ pour tout $1\le p<\infty$, et que $du_n \to du$
dans $L^{2-\e}$ pour tout $\e \in (0, 1)$.  Soit $\chi$ une fonction
lisse positive dans $D$ {\`a} support compact.  Stokes donne:
\begin{eqnarray*}
\int_D d \left( \chi u_n d^c u_n\right) 
&=&
\int_{\partial D}\chi u_n\, d^c u_n =0
\\
&=&
\int_D \chi\, d u_n \wedge d^c u_n
+
\int_D u_n\, d \chi  \wedge d^c u_n
+ 
\int_D \chi u_n\, dd^c u_n~.
\end{eqnarray*}
On a donc $0 \le \int_D \chi\, d u_n \wedge d^c u_n = - \int_D u_n\, d
\chi \wedge d^c u_n - \int_D \chi u_n\, dd^c u_n$.  On va montrer que
les deux derniers termes sont born{\'e}s uniform{\'e}ment en $n$.  On peut
donc extraire une sous-suite de la suite de $1$-formes $du_n$
convergeant faiblement dans $L^2_\mathrm{loc}$. Comme $du_n \to du$ presque
partout, le th{\'e}or{\`e}me de convergence domin{\'e} implique que la
limite est n{\'e}cessairement $du$, ce qui montre que $du$ est {\`a}
coefficients $L^2_\mathrm{loc}$.

On a tout d'abord $d^c u_n \to d^c u$ dans $L^{2-\e}_{\mathrm{loc}}$
pour tout $\e>0$, et $u_n \to u$ dans $L^p_{\mathrm{loc}}$ pour
tout $p>0$ donc $u_n d^c u_n \to u d^c u$ dans $L^1_{\mathrm{loc}}$.
Comme $\chi$ est lisse {\`a} support compact, la suite $\int_D u_n\, d
\chi \wedge d^c u_n$ converge donc vers $\int_D u\, d \chi \wedge d^c
u \in \R$.

Ensuite, $u_n\ge u$ donc $\int_D \chi u_n\, dd^c u_n \ge \int_D \chi
u\, dd^c u_n$.  La fonction $u$ est born{\'e}e inf{\'e}rieurement sur le
support de $\chi$ par une constante $C$.  On~a donc l'in{\'e}galit{\'e}
$\int_D \chi u_n\, dd^c u_n \ge C \int_D \chi dd^c u_n$.  Ce dernier
terme est uniform{\'e}ment born{\'e}. Ceci termine la preuve du lemme.
\end{proof}
Nous aurons aussi besoin du lemme suivant.
\begin{lemme}\label{L-critint}
Soient $u$ et $v$ deux fonctions sous-harmoniques. Si $v$ est
localement born{\'e}e, alors $u$ est localement int{\'e}grable par rapport
{\`a} la mesure $dd^c v$.
\end{lemme}
\begin{proof}
Soient $u_n \to u$, $v_m \to v$ deux suites r{\'e}gularisantes et $\chi$
une fonction test d{\'e}finie sur un domaine $D$.
On a vu {\`a} la preuve du Lemme~\ref{L-intsh} que pour tout $n$,
\begin{eqnarray*}
\int_D \chi\, d u_n \wedge d^c v_m
+
\int_D u_n\, d \chi  \wedge d^c v_m
+ 
\int_D \chi u_n\, dd^c v_m
 =0 ;
\\
\int_D \chi\, d v_m \wedge d^c u_n
+
\int_D v_m\, d \chi  \wedge d^c u_n
+ 
\int_D \chi v_m\, dd^c u_n
 =0 ~.
\end{eqnarray*}
Comme pr{\'e}c{\'e}demment, les termes $ \int_D u_n\, d \chi \wedge
d^c v_m$ et $\int_D v_m\, d \chi \wedge d^c u_n$ convergent
respectivement vers $ \int_D u\, d \chi \wedge d^c v$ et $\int_D v\, d
\chi \wedge d^c u$ et sont finis.  Quand $n\to \infty$, $dd^c u_n \to
dd^cu$ et comme $v$ est localement born{\'e}e, $C \ge v_m \ge -C$ pour
une constante $C>0$ et pour tout $m$. La suite $\int_D \chi v_m\,
dd^c u_n$ est donc uniform{\'e}ment born{\'e}e.  Enfin on~a la sym{\'e}trie $d
u_n \wedge d^c v_m = d v_m \wedge d^c u_n$. On en d{\'e}duit donc
que $\int_D \chi u_n\, dd^c v_m$ est uniform{\'e}ment born{\'e}e en
$n,m$. Comme $dd^c v_m \to dd^cv$ et $u_n$ d{\'e}croit vers $u$, il
s'ensuit que $u \in L_\mathrm{loc}^1( dd^cv )$.
\end{proof}

%
%%%%%%%%%%%%%%%%%%%%%%%%%%%%%%%%%%%%%%%%%%%%%%%%%%%%%%%%%%%%%%%%%%%%%
%
%

\subsection{Energie}\label{S-energie}
Notons $\diag = \{ (z,z)\, ; \, z \in \C\}$ la diagonale de $\C \times \C$.
Soient $\rho$ et $\rho'$ deux mesures sign{\'e}es sur
$\pc$. Notons $|\rho|$ et $|\rho'|$ leur mesure trace et supposons que
$\log |z-w| \in L^1( |\rho| \otimes |\rho'|)$ dans $\C^2\setminus
\diag$.  On d{\'e}finit alors l'{\'e}nergie mutuelle de $\rho$ et $\rho'$
par l'int{\'e}grale:
\begin{equation}\label{e-energie-complexe}
(\rho, \rho') \= - \int_{\C \times \C \setminus \diag} \log | z-w|\, d\rho(z)
\otimes d\rho'(w)~.
\end{equation}
Lorsque $\rho = \sum m_i [z_i]$, $\rho' = \sum m_j' [z_j']$ sont deux
mesures {\`a} support fini, l'hypoth{\`e}se
d'int{\'e}\-gra\-bi\-lit{\'e} est imm{\'e}diatement satisfaite, et on~a 
$$ (\rho, \rho') = - \sum_S m_i m_j' \log |z_i - z_j'| \text{ avec } S =
\{ (i,j) \, ; i\ne j , \, z_i \ne \infty, \, z_j' \ne \infty \} ~.$$
Dans la suite, on utilisera aussi le crit{\`e}re d'int{\'e}grabilit{\'e}
suivant.  Pour toute mesure positive $\rho$ d{\'e}finie sur $\pc$, on
utilisera la locution $\rho$ est {\`a} \emph{potentiel continu} pour dire
que localement en tout point de $\pc$, $\rho = \Delta u$ avec $u$
continu. 
\begin{lemme}\label{L-integr}
Soit $\rho$ une mesure sign{\'e}e dont la mesure trace est {\`a}
potentiel continu.  Soit de plus $\rho'$
une mesure satisfaisant {\`a} l'une des propri{\'e}t{\'e}s suivantes~:
\begin{itemize}
\item[$\bullet$]
$\rho'$ est une mesure {\`a} support fini ne chargeant pas l'infini;
\item[$\bullet$]
$|\rho'|$ est {\`a} potentiel continu.
\end{itemize}
Alors $\log | z - w | \in
L^1( |\rho| \otimes |\rho'|)$ dans $\pc \times \pc$.
En particulier, $(\rho, \rho')$ est bien d{\'e}fini.
\end{lemme}
\begin{proof}
Dans le premier cas, il suffit par lin{\'e}arit{\'e} de montrer que $\log
| z- w_0| \in L^1(|\rho|)$ pour tout $w_0\in\C$ fix{\'e}.  Localement en
un point $z\in \C$, cel{\`a} r{\'e}sulte du Lemme~\ref{L-critint} car
$|\rho|$ admet un potentiel continu. Au point infini, on peut faire le
changement de variable $Z = 1/z$ et on~a alors $\log | z- w_0| = \log
| 1- Zw_0| - \log |Z|$ qui est une diff{\'e}rence de deux fonctions
sous-harmoniques. Donc $\log | z- w_0|$ est aussi localement
int{\'e}grable par rapport {\`a} $|\rho|$ au voisinage de l'infini. 

Pour d{\'e}montrer l'int{\'e}grabilit{\'e} de $\log|z-w|$ dans le second
cas, on remarque que par lin{\'e}arit{\'e}, il suffit de le v{\'e}rifier
lorsque $\rho$ et $\rho'$ sont des mesures positives {\`a} potentiel
 continu. C'est compl{\`e}tement clair dans $\C \times \C$ priv{\'e}
de la diagonale. Localement en un point de la diagonale, cel{\`a}
r{\'e}sulte comme pr{\'e}c{\'e}demment du Lemme~\ref{L-critint}, combin{\'e}
maintenant au th{\'e}or{\`e}me de Fubini. En un point de la forme
$(\infty, w)$ avec $w \in \C$, le changement de variables $Z = 1/z$
donne $\log |z -w| = \log | 1- Zw| - \log |Z|$, qui est localement une
diff{\'e}rence de deux fonctions sous-harmoniques. Le raisonnement
pr{\'e}c{\'e}dent s'applique donc. Au point $(\infty,\infty)$, on pose $Z
= 1/z$ et $W= 1/w$. On~a alors $\log |z -w| = \log | Z- W| - \log |Z|
- \log |W|$. On conclut de m{\^e}me que $\log|z-w|$ est int{\'e}grable par
rapport {\`a} $|\rho| \otimes |\rho'|$ en ce point. On~a donc prouv{\'e}
que $\log|z-w|\in L^1(|\rho| \otimes |\rho'|)$. 
\end{proof}
\begin{lemme}\label{L-cal}
Soient $\rho,\rho'$ deux mesures sign{\'e}es telles que 
$\log|z-w|\in L^1(|\rho| \otimes |\rho'|)$ dans $\pc \times \pc$.
Alors la fonction $g_\rho ( z) \= \int_\C \log |z-w|\, d\rho(w)$
est int{\'e}grable par rapport {\`a} $\rho'$ et on~a
\begin{equation}\label{e-cal}
(\rho , \rho') = - \int_\C g _\rho d \rho' ~.
\end{equation}
\end{lemme}
\begin{proof}
L'hypoth{\`e}se implique tout d'abord que $\rho \otimes \rho'$ ne charge
ni $\diag$, ni $\{ \infty\} \times \pc$, ni $\pc \times \{
\infty\}$. De plus, Fubini implique que $g_\rho ( z) \= \int_\C \log
|z-w| d\rho(w)$ est bien d{\'e}fini pour presque tout $z$ et que l'on a
$g_\rho \in L^1(|\rho'|)$.  On peut donc {\'e}crire:
\begin{multline*}
(\rho,\rho')
=
-\int_{\C \times \C \setminus \diag} \log |z-w|\, d\rho(z) \otimes d\rho'(w)
=\\=
-\int_{\C \times \C} \log |z-w|\, d\rho(z) \otimes d\rho'(w)
=\\=
-\int_{w\in \C} \left[ \int_{z\in \C} \log |z-w|\, d\rho(z)\right]   d\rho'(w)
=
-\int_\C g_\rho d\rho'~.
\end{multline*}
Ce qui conclut la preuve.
\end{proof}

%
%%%%%%%%%%%%%%%%%%%%%%%%%%%%%%%%%%%%%%%%%%%%%%%%%%%%%%%%%%%%%%%%%%%%%
%
%

\subsection{Positivit{\'e}.}\label{S-positive}
On va maintenant montrer que l'{\'e}nergie d'une mesure $\rho$ poss{\`e}de
des propri{\'e}t{\'e}s de positivit{\'e}, au moins lorsque $\rho(\pc)=0$ et
lorsque $\rho$ admet un potentiel suffisamment r{\'e}gulier.  Pour
m{\'e}moire, notons qu'une mesure sign{\'e}e $\rho$ d{\'e}finie sur $\pc$
s'{\'e}crit $\rho = \Delta g$ avec $g: \pc \to \R$ localement
int{\'e}grable si et seulement si $\rho ( \pc) = 0$.
\begin{proposition}\label{P-casdeg}
  Soit $\rho$ une mesure sign{\'e}e sur $\pc$ telle que $\rho (\pc) =0$,
  et dont la mesure trace est {\`a} potentiel continu. On peut alors
  {\'e}crire $\rho = \Delta g$ avec $g$ continue, et $(\rho, \rho)$ est
  bien d{\'e}finie. De plus, $dg$ est une $1$-forme {\`a} coefficients $L^2$
  et on~a
\begin{equation}\label{e-140}
(\rho, \rho ) =  \int_\pc dg \wedge d^c g \ge 0~.
\end{equation}
De plus $(\rho, \rho ) = 0$ si et seulement si $\rho = 0$.
\end{proposition}
\begin{proof}
  Ecrivons la d{\'e}com\-position de Jordan de $\Delta g$ sous la forme
  $\Delta g = \rho_1 - \rho_2$, avec $\rho_1$ et $\rho_2$ deux mesures
  positives de m{\^e}me masse.  Par hypoth{\`e}se, localement en tout point,
  $\rho_1$ et $\rho_2$ admettent des potentiels continus. Au
  voisinage d'un point $p$, on peut donc bien {\'e}crire $\rho = \Delta
  h$ avec $h$ continue.  Mais $g-h$ est harmonique, donc $g$ est aussi
  continue.  Le fait que $(\rho,\rho)$ est bien d{\'e}finie r{\'e}sulte du
  Lemme~\ref{L-integr}. Pour montrer que $dg$ est {\`a} coefficients
  $L^2$, on applique le Lemme~\ref{L-intsh} aux potentiels locaux de
  $\rho_1$ et $\rho_2$. Les formes $dg$ et $d^cg$ sont donc toutes
  deux localement {\`a} coefficients $L^2$, donc globalement $L^2$ par
  compacit{\'e} de $\pc$.  On fixe alors $R>0$ tr{\`e}s grand. Stokes donne: $
  \int_{D(0,R)} dg \wedge d^c g = -\int_{D(0,R)} g \, dd^c g +
  \int_{\{|z|=R\}} g d^c g $. Comme $g$ est continue, et quitte {\`a}
  remplacer $g$ par $g - g(\infty)$, on peut supposer que $g\to0$ {\`a}
  l'infini, donc le dernier terme tend vers $0$ pour une sous-suite
  $R_n$ croissant vers l'infini, convenablement choisie. En passant {\`a}
  la limite, on en d{\'e}duit que
\begin{equation}\label{e-130}
\int_\pc dg \wedge d^c g= -
\int_\pc g\, dd^c g~. 
\end{equation}
Le Lemme~\ref{L-cal} implique que $g_\rho (z) = \int_{w\in \C} \log
|z-w|\, d\rho(w)$ est bien d{\'e}finie et int{\'e}grable par rapport {\`a}
$\rho$.  Mais $\Delta g_\rho = \rho = \Delta g$ donc $g_\rho - g$ est
constante.  Comme $\rho$ a un potentiel local continu en l'infini, on
a $\rho \{\infty\} = 0$. Donc $\rho(\pc)=0$ implique $\lim_{z\to \infty}
g_\rho (z)=0 = g(\infty)$. On conclut que $g = g_\rho$.  
Finalement~\eqref{e-cal} implique
\begin{equation*}
(\rho,\rho)
=
- \int_\C g_\rho\, d\rho
= 
-\int_\C g \, dd^c g = \int_\pc dg \wedge d^c g~.
\end{equation*}
Ceci d{\'e}montre~\eqref{e-140}.  Enfin $(\rho,\rho) =0$ implique $dg
=0$ presque partout (rappelons que $dg$ est une $1$-forme {\`a}
coefficients $L^2$). La continuit{\'e} de $g$ et le fait que $g(\infty)
=0$ implique que $g\equiv 0$. En particulier, $\rho =0$.
\end{proof}
%
%%%%%%%%%%%%%%%%%%%%%%%%%%%%%%%%%%%%%%%%%%%%%%%%%%%%%%%%%%%%%%%%%%%%%
%
%
\subsection{Energie et r{\'e}gularisation}
Dans ce paragraphe, on d{\'e}montre la Proposi\-tion~\ref{P-130}. Ce
r{\'e}sultat sera fondamental dans le reste de l'article.

On fixe dans toute la suite une fonction lisse d{\'e}croissante $\varphi
: [0, \infty) \to [0, 1]$, telle que $\varphi\equiv 0$ hors du segment
$[0,1]$ et $\int_0^1\varphi=1$. Pour tout $\e>0$, on note
$\varphi_{\e}(r) = \e^{-1} \varphi(\frac{r}{\e})$.  Pour toute
fonction continue $\chi$ sur $\P^1(\C)$, on d{\'e}finit
\begin{equation}\label{e-160}
\chi_\e(z) \= \int_0^\e \left[ \int_0^{2\pi} \chi
(z + r e^{it}) \frac{dt}{2\pi} \right] \varphi_\e(r) dr~,
\end{equation}
avec la convention $\chi_\e(+\infty) = \chi(\infty)$.
C'est une fonction lisse sur $\C$, continue sur $\P^1(\C)$, et
$\chi_\e$ tend uniform{\'e}ment vers $\chi$ lorsque $\e\to0$.  Pour toute
mesure sign{\'e}e $\rho$ sur $\pc$, on d{\'e}finit sa convol{\'e}e $\rho_\e \=
\varphi_\e \ast \rho$ en posant $\int \chi d\rho_\e \= \int \chi_\e
d\rho$ pour toute fonction continue $\chi$.  Si $\rho$ est une mesure
de probabilit{\'e}, $\rho_\e$ est encore une mesure de probabilit{\'e}. On
v{\'e}rifie facilement le
\begin{lemme}\label{l:pourtefaire plaisir}
  Pour toute mesure sign{\'e}e $\rho$, on~a $\rho_\e \to \rho$ lorsque
  $\e\to0$. Pour toute suite de mesures sign{\'e}es telle que
  $\rho_n\to\rho$, on~a $\rho_{n,\e} \to \rho_\e$. Enfin si $\rho$ est
  {\`a} support compact dans $\C$, alors $\rho_\e$ est {\`a} potentiel
  continu.
\end{lemme}
En particulier, si $z\in \C$, la mesure de probabilit{\'e} $[z]_\e$ est
une mesure {\`a} potentiel continu.  Dans la suite, on utilisera la
terminologie suivante. Si $F\subset \C$ est un ensemble fini de
cardinal $|F|$, on notera
$$
[F] = |F|^{-1}\sum_{z\in F} [z]
\ \mbox{ et } \
[F]_\e= |F|^{-1} \sum_{z\in F} [z]_\e~.
$$
Dans la proposition suivante on note par $\lambda_1$ la mesure de
probabilit{\'e} proportionnelle {\`a} la mesure de Lebesgue sur le cercle
unit{\'e}.

Rappelons qu'un \emph{module de continuit{\'e}} pour une fonction continue
$h$ sur $\P^1(\C)$ est une fonction $\eta: \R_+ \to \R_+$, telle que
pour tous points $z,w \in \P^1(\C)$ tels que $\mathsf{d}(z,w)\le \e$, on ait
$|h(z) - h(w)|\le \eta(\e)$. Ici $\mathsf{d}$ d{\'e}note la m{\'e}trique sph{\'e}rique sur
la sph{\`e}re de Riemann. Notons que la m{\'e}trique euclidienne $|\cdot|$ sur $\C$
est plus grande que la m{\'e}trique sph{\'e}rique, donc tout module de
continuit{\'e} pour $\mathsf{d}$ est un module de continuit{\'e} pour $|\cdot|$.
\begin{proposition}\label{P-130}
  Soit $\rho$ une mesure de probabilit{\'e} {\`a} potentiel continu, soit
  $\heta$ un module de continuit{\'e} d'un potentiel de $\rho - \lambda_1$
  et posons $\eta(\e) \= \heta(\e) + \e$.  Alors il existe une
  constante $C > 0$ telle que pour tout $\e>0$ et tout sous ensemble
  fini $F$ de $\C$, on ait
\begin{align}
([F] - \rho , [F] -\rho) 
& \ge  ([F]_\e - \rho , [F]_\e -\rho) 
- 2\, \eta(\e) - |F|^{-1} (C + \log \e^{-1}) \label{e-190}\\
& \ge  - 2\, \eta(\e) - |F|^{-1} (C + \log \e^{-1}) \label{e-191}
\end{align}
\end{proposition}
La preuve, qui sera donn{\'e}e ci-dessous, repose sur
l'{\'e}tude du comportement de l'{\'e}nergie apr{\`e}s r{\'e}gularisation des
mesures.
\begin{lemme}\label{L-100}
Soient $\rho$ et $\eta$ comme dans la proposition.
Alors pour tout ensemble fini de points $F\subset \C$, on~a
\begin{equation}\label{e-150}
|([F], \rho ) - ([F]_\e , \rho )|\le \eta(\e)~.
\end{equation}
\end{lemme}
\begin{lemme}\label{L-120}
Il existe une constante $C$ telle que pour tout $\e>0$ et tout
ensemble fini $F$ de points, on ait
\begin{equation}\label{e-170}
([F]_\e , [F]_\e ) \le ([F], [F]) + |F|^{-1} (C + \log \e^{-1}) ~.
\end{equation}
\end{lemme}
\begin{proof}[Preuve de la proposition~\ref{P-130}]
Comme $\rho$ et $[F]_\e$ sont des mesures de probabilit{\'e} {\`a}
potentiel continu, il existe une
fonction $g$ d{\'e}finie et continue sur $\pc$ telle que $\Delta g =
\rho - [F]_\e$.  L'{\'e}quation~\eqref{e-140} implique alors
$([F]_\e-\rho, [F]_\e-\rho) \ge0$.  En particulier,~\eqref{e-190}
implique~\eqref{e-191}.

Le r{\'e}sultat est alors une cons{\'e}quence imm{\'e}diate des
Lemmes~\ref{L-100} et~\ref{L-120}.
\end{proof}

\begin{proof}[Preuve du Lemme~\ref{L-100}]
On traite tout d'abord le cas o{\`u} $\rho=\lambda_1$.
On a alors $g_{\lambda_1}(z) = \int_\C
\log|z-w|\, d\lambda_1(w) = \log^+|w|$ sur $\C$.  Pour tout $z\in \C$,
les {\'e}quations~\eqref{e-cal} et~\eqref{e-160} donnent alors:
\begin{multline*}
([z]_\e, \lambda_1)
- 
([z], \lambda_1)
=
\int_\C \log^+|w| (\varphi_\e \ast [z] - [z])
=\\=
\int_0^\e \left[ \int_0^{2\pi}\left( 
\log^+\left|z + r e^{it}\right| - \log^+|z|\right)
 \frac{dt}{2\pi} \right] \varphi_\e(r) dr~.
\end{multline*}
On v{\'e}rifie facilement que $\left| \log^+\left|z + r e^{it}\right| -
\log^+|z|\, \right| \le \e$ pour tout $z,t$ et tout $r\le \e$.  On a
donc $|([z]_\e, \lambda_1) - ([z], \lambda_1)|\le \e$ et par suite
$|([F]_\e, \lambda_1) - ([F], \lambda_1)|\le \e$.

Dans le cas g{\'e}n{\'e}ral on {\'e}crit $\rho - \lambda_1 = \Delta h$ sur $\pc$.
Par hypoth{\`e}se $h$ est continue et la fonction $\heta$ est un module de
continuit{\'e} de $h$.  Sur $\C$, on~a alors $\rho = \Delta g $ avec
$g\=\log^+|w| + h$.  De~\eqref{e-cal}, on d{\'e}duit
\begin{multline*}
([F], \rho) -  ([F]_\e,\rho)
= \int_\C g \, d ([F]_\e - [F])
=\\=
 \int_\C h \, d ([F]_\e - [F]) + ([F]_\e, \lambda_1) -  ([F],\lambda_1)~.
\end{multline*}
Ce dernier terme est born{\'e} par $\eta (\e) \= \widehat{\eta} (\e) + \e$, ce
qui termine la d{\'e}monstration du lemme.
\end{proof}

\begin{proof}[Preuve du Lemme~\ref{L-120}]
Fixons $r, r' > 0$ et $z, z' \in \C$ distincts.
Alors
\begin{multline*}
\int_{[0,2\pi]^2} \log 
\left| (z+ re^{it}) - (z'+r'e^{it'}) \right|
\frac{dt\otimes dt'}{(2\pi)^2}
= \\ =
\int_{[0,2\pi]}
\max \left\{ \log 
| z - (z'+r'e^{it'})| , \log r
\right\} \frac{dt'}{2\pi}
\ge \\ \ge
\max \left\{ \int_{[0,2\pi]}
\log 
| z - (z'+r'e^{it'})| 
\,  \frac{dt'}{2\pi} ,\,  \log r\right\}
\ge \\ \ge 
\max \{ \log|z-z'|, \log r', \log r \} 
\ge \log |z-z'| 
~.
\end{multline*}
Pour tous $z, z' \in \C$ distincts, on obtient $([z]_\e, [z']_\e) \le
- \log |z-z'| = ([z], [z'])$ en int{\'e}\-grant cette suite
d'in{\'e}galit{\'e}s.  Lorsque $z=z'$, les in{\'e}galit{\'e}s pr{\'e}c{\'e}dentes
se r{\'e}{\'e}crivent
\begin{multline*}
\int_{[0,2\pi]^2} \log 
\left| (z+re^{it}) - (z+r'e^{it'}) \right|
\frac{dt\otimes dt'}{(2\pi)^2}
= \\ =
\int_{[0,2\pi]^2} \log 
\left| re^{it} - r'e^{it'} \right|
\frac{dt\otimes dt'}{(2\pi)^2}
= \max \{ \log r , \log r' \}~,
\end{multline*}
et donc
\begin{multline*}
([z]_\e, [z]_\e) \le - \int_{[0,\e]^2} \max \{ \log r , \log r' \}
\varphi_\e(r) \varphi_\e(r')\, dr \otimes dr' 
\le \\ \le -
\int_0^\e \log r \cdot \varphi_\e(r)  dr 
=
C + \log \e^{-1}~,
\end{multline*}
 pour une certaine constante $C$.  On d{\'e}duit de tout cel{\`a}:
\begin{multline*}
([F]_\e,[F]_\e)
=
\frac{1}{|F|^2}\sum_{z\ne z'\in F} ([z]_\e,[z']_\e) 
+
\frac1{|F|^2} \sum_{z\in F} ([z]_\e,[z]_\e) 
\le \\ \le
\frac1{|F|^2} \sum_{z\ne z'\in F} ([z],[z'])
 +
\frac1{|F|^2} \sum_{z\in F} ([z]_\e,[z]_\e)
=\\=
([F], [F]) 
+
\frac1{|F|^2} \sum_{z\in F} ([z]_\e,[z]_\e)~.
\end{multline*}
Le dernier terme se majore par $(C + \log \e^{-1})/|F|$, ce qui termine
la d{\'e}monstration.
\end{proof}

%
%%%%%%%%%%%%%%%%%%%%%%%%%%%%%%%%%%%%%%%%%%%%%%%%%%%%%%%%%%%%%%%%%%%%%
%
%
\subsection{Le r{\'e}sultat clef}
Nous pouvons maintenant d{\'e}montrer le r{\'e}sultat clef qui nous servira
dans la preuve du th{\'e}or{\`e}me principal de l'article.
\begin{proposition}\label{P-casdeg2}
Soit $\rho$ une mesure de probabilit{\'e} sur $\pc$, telle que
localement en tout point, il existe une fonction $g$ continue
 v{\'e}rifiant $\rho = \Delta g$.  Soit $F_n \subset \P^1(\C)$ une
suite d'ensembles finis telle que $|F_n| \to \infty$. 
Alors
$$
\ov{\lim}_{n\to \infty} ([F_n] - \rho, [F_n] - \rho ) \le 0 
\, \text{ implique } \,
\lim_{n\to \infty} [F_n] = \rho~.
$$
\end{proposition}

\begin{proof}
  Tout d'abord notons que l'on peut toujours supposer que $F_n\subset
  \C$. En effet, par construction $([F_n] - \rho, [F_n] - \rho )\to0$
  et $|F_n|\to\infty$ impliquent $ ([F_n\setminus\{\infty\}] - \rho,
  [F_n\setminus\{\infty\}] - \rho )\to0$.

Pour d{\'e}montrer la proposition, on se ram{\`e}ne {\`a} la
Proposition~\ref{P-casdeg} en r{\'e}gu\-la\-ri\-sant les mesures
$[F_n]$.  Les mesures $[F_n]_\e$ sont absolument continues par rapport
{\`a} la mesure de Lebesgue sur $\C$, on peut donc {\'e}crire $[F_n]_\e -
\rho = \Delta g_{n,\e}$ avec $g_{n,\e}$ continues, que l'on
normalisera par $g_{n,\e} (\infty) = 0$.  On va montrer que $\Delta
g_{n, \e}$ converge faiblement vers la mesure nulle lorsque l'on fait tendre
$n \to \infty$ puis $\e \to 0$.  Comme les mesures $\Delta g_{n,
\e}$ sont de masse totale {\'e}gale {\`a}~$0$, il suffit de montrer que
pour toute fonction lisse $\chi$ {\`a} support compact dans $\C$ on~a
\begin{equation}\label{e:presque-gagne}
\lim_{\e \to 0} \lim_{n \to \infty} \int_\C \chi dd^c g_{n, \e} = 0.
\end{equation}
En effet, prenons un point d'adh{\'e}rence $\tilde{\rho}$ de la suite de
mesures de probabilit{\'e} $[F_n]$. On peut donc trouver une sous-suite
$[F_{n_k}] \to \tilde{\rho}$. L'op{\'e}rateur de r{\'e}gularisation est
continue dans l'espace des mesures, donc pour tout $\e>0$, on~a
$[F_{n_k}]_\e \to \tilde{\rho}_\e$, voir Lemme~\ref{l:pourtefaire plaisir}.
L'{\'e}quation~\eqref{e:presque-gagne} se traduit par l'{\'e}galit{\'e}
$\lim_{\e\to0}\tilde{\rho}_\e =\rho$, et on en d{\'e}duit donc
$\tilde{\rho}=\rho$. Tous les points d'adh{\'e}rence de $[F_n]$ {\'e}tant
{\'e}gaux {\`a} $\rho$, on conclut $[F_n]\to\rho$.

Pour montrer~\eqref{e:presque-gagne}, on proc{\`e}de comme suit.
De~\eqref{e-190}, on tire
\begin{equation*}
([F_n]_\e - \rho, [F_n]_\e -\rho ) 
\le
([F_n] - \rho , [F_n] - \rho ) + |F_n|^{-1} (C + \log \e^{-1}) 
 + 2 \eta ( \e) ~.
\end{equation*}
Comme $|F_n|\to \infty$ et que par hypoth{\`e}se
$\ov{\lim}_{n \to \infty} ([F_n]-\rho,[F_n]-\rho) \le 0$, on en d{\'e}duit que
\begin{equation}\label{e-310}
 \blim_{\e \to 0} \blim_{n\to \infty} ([F_n]_\e - \rho, [F_n]_\e -\rho
 ) \le 0~.
\end{equation}

De~\eqref{e-140}, on tire $ ([F_n]_\e - \rho, [F_n]_\e -\rho ) =
\int_\C dg_{n,\e} \wedge d^c g_{n,\e}\ge 0$.  Au vu de~\eqref{e-310},
on en d{\'e}duit $ \lim_{\e \to 0} \lim_{n\to \infty} \int_\C
dg_{n,\e} \wedge d^cg_{n,\e} =0$, puis $ \lim_{\e \to 0} \lim_{n\to
\infty} d^c g_{n,\e} =0$ dans $L^2$, donc dans $L^1$.  Si $\chi$ est
une fonction {\`a} support compact dans $\C$ quelconque, le
Th{\'e}or{\`e}me de Stokes implique
$$
\lim_{\e \to 0} \lim_{n \to \infty} \int_\C \chi dd^c g_{n, \e}
=
\lim_{\e \to 0} \lim_{n \to \infty} - \int_\C d \chi \wedge d^c g_{n, \e}
=
0.
$$
Ceci termine la preuve de la proposition.
\end{proof}

%
%%%%%%%%%%%%%%%%%%%%%%%%%%%%%%%%%%%%%%%%%%%%%%%%%%%%%%%%%%%%%%%%%%
%
%%%%%%%%%%%%%%%%%%%%%%%%%%%%%%%%%%%%%%%%%%%%%%%%%%%%%%%%%%%

\section{L'espace de Berkovich de $\C_p$.}\label{S-berkovich}

L'espace $\C_p$ muni de sa norme $p$-adique est un espace totalement
discontinu et non localement compact et de ce fait se pr{\`e}te mal {\`a}
l'analyse ou {\`a} la th{\'e}orie de la mesure.  Pour contourner cette
difficult{\'e}, on ``connexifie'' $\C_p$ en construisant un arbre $\berCp$
dans lequel $\C_p$ s'identifie {\`a} un sous-espace de ses bouts.  Cette
construction d{\^u}e {\`a} Berkovich s'av{\`e}re tout {\`a} fait fondamentale.  Nous
verrons au paragraphe suivant qu'il est ainsi possible de construire
un op{\'e}rateur de Laplace convenable sur $\berCp$.  Dans ce paragraphe,
nous d{\'e}crivons les propri{\'e}t{\'e}s topologiques essentielles de $\berCp$.

%%%%%%%%%%%%%%%%%%%%%%%%%%%%%%%%%%%%%%%%%%%%%%%%%%%%%%%%%%%

\subsection{Le corps $\C_p$.}
Fixons une cl{\^o}ture alg{\'e}brique $\ov{\Q}$ du corps des nombres rationnels
$\Q$ et un nombre premier~$p$.  On d{\'e}signe par $| \cdot |$ la norme
$p$-adique sur $\Q$, normalis{\'e}e par $|p| = p^{-1}$.  Cette norme
s'{\'e}tend de fa{\c c}on unique en une norme d{\'e}finie sur la compl{\'e}tion
$\Q_p$ du corps valu{\'e} $(\Q, | \cdot |)$, puis sur une clot{\^u}re
alg{\'e}brique $\ov{\Q_p}$. On d{\'e}signera toutes ces normes par $| \cdot
|$.  On notera enfin $\C_p$ la compl{\'e}tion de  $(\ov{\Q_p}, |
\cdot |)$.  Le groupe
$$
| \C_p^* | = \{ |z| ; \,  z \in \C_p^* \}
$$ est appel{\'e} le {\it groupe des valeurs} et il est {\'e}gal {\`a} $\{
p^r ; \,  r \in \Q \}$.  Enfin, on notera $\pp$ la droite projective de
$\C_p$, que l'on peut identifier naturellement {\`a} $\C_p \cup \{
\infty \}$.

%%%%%%%%%%%%%%%%%%%%%%%%%%%%%%%%%%%%%%%%%%%%%%%%%%%%%%%%%%%

\subsection{L'espace des semi-normes.}\label{espace des semi-normes}
Soit $\berCp$ l'espace de toutes les semi-normes multiplicatives
d{\'e}finies sur $\C_p[T]$, dont la restriction {\`a} $\C_p$ est {\'e}gale {\`a}
$|\cdot|$.  On note de plus $\cS_\infty$ la fonction d{\'e}finie sur
$\C_p[T]$, qui est constante {\'e}gale {\`a}~$\infty$ sur tous les polyn{\^o}mes
non constants de $\C_p$ et telle que pour chaque polyn{\^o}me constant $P
\equiv a$ on ait $\cS_\infty(P) = |a|$.  On pose $\berp = \berCp \sqcup
\{ \cS_\infty \}$ et on munit $\berp$ de la topologie la moins fine
telle que pour chaque $P \in \C_p[T]$ la fonction $\cS \mapsto \cS(P)$
soit continue.  L'espace $\berp$ est alors compact et sa topologie
admet une base d{\'e}nombrable.  On l'appelle \emph{espace analytique de
  Berkovich} associ{\'e} {\`a}~$\pp$.

Tout $z \in \C_p$ induit une semi-norme $\cS_z$ d{\'e}finie par
$\cS_z(P) = |P(z)|$.
L'application $z \mapsto \cS_z$ est un
hom{\'e}omorphisme de $\pp = \C_p \cup \{ \infty \}$ sur son image.
Dans la suite, on identifiera $\pp$ avec son image dans $\berp$.

A chaque boule $B = \{ |z - z_0| \le r \}$ correspond la semi-norme
$\cS_B$ dans $\berCp$, d{\'e}finie par $\cS_B(P) = \sup_B |P(z)|$.  Plus
g{\'e}n{\'e}ralement, toute suite d{\'e}croissante $\{ B_i \}_{i \ge 0}$ de
boules de $\C_p$ induit une semi-norme $P \mapsto \lim_{i \to \infty}
\cS_{B_i}(P)$. R{\'e}cipro\-que\-ment, toute semi-norme dans $\berCp$ est de
cette forme et les points de $\berp$ se rangent donc dans l'une des
quatre cat{\'e}gories suivantes (voir par exemple~\cite[p.18]{Ber}):
\begin{enumerate}
\item[{\rm i)}]
{\it les points de $\pp$};
\item[{\rm ii)}] {\it les points rationnels}, de la forme $\cS_B$,
avec $B = \{ |z - a| \le r \}$ et $r \in |\C_p^*|$;
\item[{\rm iii)}] {\it les points irrationnels}, de la forme $\cS_B$,
avec $B = \{ |z - a| \le r \}$ et $r \not \in |\C_p^*|$;
\item[{\rm iv)}] {\it les points singuliers}, associ{\'e}s {\`a} une suite
d{\'e}croissante de boules de $\C_p$ dont l'intersection est vide.
\end{enumerate}
Notons que tous les points de type (ii),~(iii) et~(iv) sont des normes
qui s'{\'e}tendent {\`a} $\C_p(T)$, alors que la semi-norme associ{\'e}e {\`a} un
point $z\in\pp$ v{\'e}rifie $\cS_z(T-z) =0$.

On appelle \emph{point canonique} la norme associ{\'e}e {\`a} la boule
unit{\'e} $\{ |z| \le 1 \}$ et on le note $\cS_\can$.  Etant donn{\'e} un
point rationnel ou irrationnel $\cS$, on d{\'e}signe par $B_\cS$ la boule
de $\C_p$ correspondante.  Lorsque $z \in \C_p$ on pose $B_z = \{ z
\}$.

Chaque fonction rationnelle $R \in \C_p(T)$ agit sur $\berp\setminus
\pp$, envoyant toute norme $\cS$ sur la norme $R_*(\cS)$ d{\'e}finie par
$R_*(\cS)(P) \= \cS ( P \circ R)$.  Cette action s'{\'e}tend contin{\^u}ment
en une action de $R_*$ sur $\berp$ qui co{\"\i}ncide avec l'action
naturelle de $R$ sur $\pp$.

\medskip

C'est un fait fondamental que $\berp$ poss{\`e}de une structure d'arbre,
que nous allons maintenant d{\'e}crire bri{\`e}vement. Consid{\'e}rons
l'ordre partiel $\le$ d{\'e}fini sur l'espace $\berp$ par~: $\cS \le \cS'$ si et
seulement si pour tout $P \in \C_p[T]$ on~a $\cS(P) \le \cS'(P)$.
Lorsque $\cS$ et $\cS'$ sont non singuliers, on~a $\cS \le \cS'$ si et
seulement si $B_\cS \subset B_{\cS'}$. On v{\'e}rifie que le point
$\cS_\infty$ est l'unique {\'e}l{\'e}ment maximal de $\berp$ et que
l'ensemble des {\'e}l{\'e}ments minimaux co{\"\i}ncide avec l'union de
$\C_p$ et des points singuliers.

Etant donn{\'e}s $\cS$ et $\cS'$ dans $\berp$,  on d{\'e}finit $\cS
\wedge \cS' \in \berp$ par
$$
(\cS \wedge \cS')(P) = \inf \{ \widehat{\cS}(P) ; \,  \widehat{\cS}
\in \berp, \ \cS \le \widehat{\cS}, \ \cS' \le \widehat{\cS} \}.
$$
On v{\'e}rifie qu'on a $\cS \wedge \cS' = \cS$ si et seulement si $\cS'
\le \cS$ et que $\cS \wedge \cS' = \cS_\infty$ si et seulement si
$\cS$ ou $\cS'$ est {\'e}gale {\`a} $\cS_\infty$.  Lorsque $\cS$ et $\cS'$
sont des points non singuliers dans $\berCp$, le point $\cS \wedge
\cS'$ est la semi-norme associ{\'e}e {\`a} la plus petite boule de $\C_p$
qui contient $B_\cS$ et $B_{\cS'}$.

L'ordre partiel $\le$ d{\'e}finit alors une structure d'arbre dans
$\berCp$ (resp. $\berp$) au sens suivant. Pour chaque paire de points
distincts $\cS$ et $\cS'$, l'ensemble
$$
[\cS, \cS']
=
\{ \widetilde{\cS} ; \,  \cS \le \widetilde{\cS} 
\le 
\cS \wedge \cS' \ \mbox{ ou } \ \cS' \le \widetilde{\cS} 
\le \cS \wedge \cS' \}.
$$ est l'unique arc topologique dans $\berCp$ (resp. $\berp$) ayant
$\cS$ et $\cS'$ comme extr{\'e}mi\-t{\'e}s. Un ensemble de la forme $ [\cS,
\cS']$ est appel{\'e} \emph{segment}.  On dira qu'un point $\cS$ est
{\it entre} les points $\cS'$ et $\cS''$ lorsque $\cS \in [\cS',
\cS'']$.  Dans ce cas on~a $[\cS', \cS''] = [\cS', \cS] \cup [\cS,
\cS'']$.  Notons que pour chaque triplet de points $\cS$, $\cS'$ et
$\cS''$ il existe un unique point qui est entre $\cS$ et $\cS'$, entre
$\cS'$ et $\cS''$ et entre $\cS''$ et $\cS$.

%%%%%%%%%%%%%%%%%%%%%%%%%%%%%%%%%%%%%%%%%%%%%%%%%%%%%%%%

\subsection{La fonction $\sup \{ \cdot, \cdot \}$}\label{S-sup}
Diverses fonctions d{\'e}finies sur $\C_p$ s'{\'e}tendent de ma\-ni{\`e}\-re
naturelle {\`a} $\berCp$ et jouent un r{\^o}le fondamental dans la
suite. La fonction $\sup$ mentionn{\'e}e dans le titre {\'e}tend la norme
$|z-z'|$ et nous permettra (entre autre) de d{\'e}finir une m{\'e}trique
naturelle sur $\berp\setminus \pp$.

\medskip

Commen{\c c}ons par d{\'e}finir les fonctions $| \cdot |$ et $\diam :
\berCp \to [0, \infty)$ comme suit.
Pour $z \in \C_p$ on pose $P_z (T) = T - z \in \C_p[T]$.
Alors,
$$
|\cS| = \cS(P_0)
\ \mbox{ et } \
\diam(\cS) = \inf_{z \in \C_p} \cS(P_z)~.
$$
Lorsque $\cS$ est un point non singulier de $\berCp$, on~a
$$
| \cS | = \sup_{B_\cS}|z|
\ \mbox{ et } \
\diam(\cS) = \diam(B_\cS).
$$ En particulier, la restriction de $| \cdot |$ {\`a} $\C_p$
co{\"\i}ncide avec la norme de $\C_p$. La fonction $| \cdot |$
s'annule uniquement au point $0$.  Pour tout $\cS \in
\berCp$, on~a $|\cS| \ge \diam(\cS)$ et $\diam (\cS) =0$ si et
seulement si $\cS \in \C_p$.  Enfin, les deux fonctions $|\cdot|$ et
$\diam$ sont continues et s'{\'e}tendent contin{\^u}ment {\`a} $\berp$ en
posant $|\infty| = \diam(\infty) = \infty$.

\medskip
A l'aide des fonctions pr{\'e}c{\'e}dentes, on d{\'e}finit maintenant:
$$ \sup \{ \cS, \cS' \} = \diam(\cS \wedge \cS')~, \text{ pour }
\cS, \cS' \in\berp ~.
$$
Lorsque $\cS$ et $\cS'$ sont des points non singuliers de $\berCp$, on~a
$$
\sup \{ \cS, \cS' \} = \sup \{ |z - z'| ; \,  z \in B_\cS, z' \in B_{\cS'} \}~,
$$ et en particulier pour tout $z, z' \in \C_p$ on~a $\sup \{ z, z' \}
= |z - z'|$.  On v{\'e}rifie ais{\'e}ment que
$$
\sup \{ \cdot, 0 \} = | \cdot |
\mbox{ et }
\sup \{ \cdot, \cS_\can \} = \max \{ |\cdot|, 1 \}~.
$$ Introduisons maintenant quelques notations.  Une \emph{boule}
ouverte (resp. ferm{\'e}e) de $\berCp$ est un ensemble de la forme $ \{
\sup \{ \cS, z \} < r \}$ (resp. $\{ \sup \{ \cS, z \} \le r \}$),
o{\`u} $z \in \C_p$ et $r > 0$.  Une boule ouverte (resp. ferm{\'e}e) de
$\berp$ est une boule ouverte (resp. ferm{\'e}e) de $\berCp$ ou le
compl{\'e}mentaire dans $\berp$ d'une boule ferm{\'e}e (resp. ouverte) de
$\berCp$.  Il est facile de voir que toute boule de $\berp$ est
connexe.  Les boules ouvertes de $\berp$ forment une sous-base de la
topologie de $\berp$.

%%%%%%%%%%%%%%%%%%%%%%%%%%%%%%%%%%%%%%%%%%%%%%%%%%%%%%%%

\subsection{Produit de Gromov dans l'espace hyperbolique~$\H_p$.}
\label{section produit de Gromov}
Notons $\H_p$ l'ensemble $\berp\setminus\pp$. La fonction $d$
d{\'e}finie par
$$ d(\cS, \cS') = 2 \log\sup \{ \cS, \cS' \} - \log \diam(\cS) - \log
\diam(\cS'),
$$ est une distance sur $\H_p$.  On a $ d(\cS, \cS') = \log \left(
\diam(\cS)/ \diam(\cS')\right)$ lorsque $\cS \le \cS'$.  L'espace
m{\'e}trique $(\H_p, d)$ est complet et c'est un arbre r{\'e}el au sens de
J.~Tits~: pour chaque paire de points $\cS, \cS'$ dans $\H_p$, l'arc
$[\cS, \cS']$ est isom{\'e}trique {\`a} l'intervalle $[0, d(\cS, \cS')]$
de $\R$. Notons de plus que $d$ est invariante par l'action du groupe
des automorphismes de $\pp$.

\medskip
Fixons un point base $\cS_0 \in \H_p$.  Le produit de Gromov est la
fonction $\langle \cdot\, , \cdot \rangle_{\cS_0} : \berp \times \berp
\to [0, \infty]$ d{\'e}finie comme suit.  Etant donn{\'e}s $\cS, \cS' \in
\berp$, notons $\cS''$ l'unique point de $\berp$ qui est entre $\cS$
et $\cS'$, entre $\cS$ et $\cS_0$ et entre $\cS'$ et $\cS_0$.
%%% (voir Figure~\ref{fig1}).
On pose alors
$$
\langle \cS , \cS' \rangle_{\cS_0}
=
\begin{cases}
d(\cS'', \cS_0) & \text{si } \cS'' \in \H_p ~;\\
\infty & \text{si } \cS'' \in \pp~.
\end{cases}
$$
On   v{\'e}rifie    facilement   que    $\langle   \cS    ,   \cS'
\rangle_{\cS_0}=\infty$ si et  seulement si $\cS = \cS'  \in \pp$; et
que $\langle \cS , \cS' \rangle_{\cS_0} = 0$ si et seulement si $\cS_0
\in  [\cS,\cS']$.
En  particulier, pour  tout  $\cS \in  \berp$ on  a
$\langle \cS , \cS_0 \rangle_{\cS_0} = 0$.
En g{\'e}n{\'e}ral, on~a $\langle \cS , \cS' \rangle_{\cS_0} \le d(\cS,
\cS_0)$, avec {\'e}galit{\'e} si et seulement si $\cS \in [\cS', \cS_0]$.
Lorsque $\cS,\cS' \in \H_p$, on v{\'e}rifie la formule:
\begin{equation*}
2 \langle \cS , \cS' \rangle_{\cS_0} 
= 
d(\cS, \cS_0) + d(\cS', \cS_0) - d(\cS, \cS')~,
\end{equation*}
et par cons{\'e}quent, pour tout $\cS, \cS', \cS_0, \cS_1 \in \H_p$, on~a
\begin{equation}\label{changement point base}
\langle \cS , \cS' \rangle_{\cS_0} - \langle \cS , \cS_1 \rangle_{\cS_0}
=
\langle \cS , \cS' \rangle_{\cS_1} - \langle \cS' , \cS_0 \rangle_{\cS_1}~.
\end{equation}
Chacun des termes {\`a} gauche et {\`a} droite est d{\'e}fini et continu
pour $\cS, \cS' \in \berp$, cette {\'e}galit{\'e} est donc aussi valable lorsque
$\cS$ ou $\cS'$ appartiennent {\`a} $\pp$.

\begin{lemme}\label{L-infini comme point base}
Pour $\cS, \cS' \in \berCp$ et $\cS^0 \in \H_p$ on~a $\log\sup\{ \cS,
\cS^0 \} = \langle \cS , \infty \rangle_{\cS^0} + \log \diam(\cS^0)$
et
\begin{equation}\label{infini comme point base}
\log \sup \{ \cS, \cS'\}
=
\langle \cS , \infty \rangle_{\cS^0} + 
\langle \cS' , \infty \rangle_{\cS^0} - 
\langle \cS , \cS' \rangle_{\cS^0} + \log\diam(\cS^0) ~.
\end{equation}
\end{lemme}
\begin{proof}
On a
\begin{eqnarray*}
\langle \cS , \infty \rangle_{\cS^0}
& = &
d(\cS \wedge \cS^0, \cS^0) \\
& = &
\log\diam(\cS \wedge \cS^0) - \log\diam(\cS^0) \\
& = &
\log \sup\{ \cS, \cS^0 \} - \log\diam(\cS^0) ~.
\end{eqnarray*}
De la m{\^e}me fa{\c c}on on~a $\langle \cS' , \infty \rangle_{\cS^0} =
\log\sup\{ \cS', \cS^0 \} - \log\diam(\cS^0)$.

Lorsque $\cS^0 \le \cS \wedge \cS'$ on~a
$$
\sup\{ \cS, \cS' \}
=
\max \{ \sup \{ \cS, \cS^0 \}, \sup \{ \cS', \cS^0 \} \}
\ \mbox{ et}
$$
$$
\langle \cS , \cS' \rangle_{\cS^0}
=
\log \min \{ \sup \{ \cS, \cS^0 \}, \sup \{ \cS', \cS^0 \} \} 
- \log \diam(\cS^0) ~,
$$
d'o{\`u} on obtient l'{\'e}quation d{\'e}sir{\'e}e.
D'autre part, lorsque $\cS^0 \not\le \cS \wedge \cS'$, on~a
$$ \diam((\cS \wedge \cS') \wedge \cS^0) = \sup\{ \cS, \cS^0 \} = \sup
\{ \cS', \cS^0 \}
$$
et par la premi{\`e}re formule,
\begin{eqnarray*}
\langle \cS , \cS' \rangle_{\cS^0}
& = &
d(\cS^0, (\cS \wedge \cS') \wedge \cS^0) 
+ d((\cS \wedge \cS') \wedge \cS^0, \cS \wedge \cS') \\
& = &
2 \log \diam((\cS \wedge \cS') \wedge \cS^0) - \log \diam(\cS^0) 
- \log\diam(\cS \wedge \cS') \\
& = &
\langle \cS , \infty \rangle_{\cS^0} + \langle \cS' , \infty \rangle_{\cS^0}
 - \log\sup \{ \cS, \cS' \} + \log\diam(\cS^0) ~.
\end{eqnarray*}
Ceci conclut la preuve.
\end{proof}

%
%%%%%%%%%%%%%%%%%%%%%%%%%%%%%%%%%%%%%%%%%%%%%%%%%%%%%%%%%%%%%%%%%%
%

\section{Energie dans $\C_p$}\label{S-energieCp}

Comme dans le cas complexe, notre preuve du th{\'e}or{\`e}me principal
repose de mani{\`e}re essentielle sur une th{\'e}orie du potentiel
adapt{\'e}e sur $\berp$. Dans cette section, nous indiquons les
{\'e}l{\'e}ments les plus importants de cette th{\'e}orie, c'est-{\`a}-dire la
construction d'un op{\'e}rateur Laplacien $\Delta$, d{\'e}fini sur un
espace convenable de potentiels $\cP$ et {\`a}
valeurs dans les mesures sign{\'e}es sur $\berp$.
Cette th{\'e}orie repose de mani{\`e}re fondamentale sur~\cite[Chapitre~7]{FJ}.
Ensuite nous introduisons la notion d'{\'e}nergie mutuelle pour deux mesures
dans $\berp$ et nous d{\'e}montrons pour celle-ci les r{\'e}sultats analogues
{\`a} ceux du paragraphe~\ref{S-pot_complexe}.

%%%%%%%%%%%%%%%%%%%%%%%%%%%%%%%%%%%%%%%%%%%%%%%%%%%%%%%%

\subsection{Op{\'e}rateur de Laplace.}\label{S-poten}

Munissons $\berp$ de la tribu des bor{\'e}liens associ{\'e}e {\`a} sa topologie
faible. On note $\cM^+$ l'ensemble des mesures bor{\'e}liennes positives
et finies, support{\'e}es dans $\berp$. On d{\'e}signe par $\cM$ l'espace
vectoriel des mesures r{\'e}elles sign{\'e}es, diff{\'e}rences de mesures dans
$\cM^+$.  Comme $\berp$ est compact (et admet une base d{\'e}nombrable),
toute suite de mesures de probabilit{\'e} dans $\cM^+$ admet une
sous-suite convergente pour la topologie de la convergence vague.
Notons que toute mesure dans $\cM$ est de Radon et est donc
repr{\'e}sent{\'e}e par une forme lin{\'e}aire continue sur l'espace des fonctions
continues de $\berp$ (voir par exemple~\cite[Proposition~7.14]{FJ}).

\medskip

Nous allons maintenant d{\'e}finir l'espace $\cP$ et l'op{\'e}rateur
$\Delta$ mentionn{\'e}s ci-des\-sus.  Pour cel{\`a}, fixons un point base
$\cS_0 \in \H_p$.  Notons que pour tout $\cS \in \H_p$, la fonction
$\cS' \mapsto \langle \cS , \cS' \rangle_{\cS_0}$ est non-n{\'e}gative
et major{\'e}e par $d(\cS, \cS_0)$.  Etant donn{\'e}e une mesure
bor{\'e}lienne $\rho \in \cM$, on peut donc d{\'e}finir $\hg_\rho : \H_p
\to \R$ par
$$ \hg_\rho(\cS) \= - \rho(\berp) - \int_\berp \langle \cS , \cS'
\rangle_{\cS_0}\,  d \rho(\cS'),
$$ et on l'appelle le {\it potentiel de $\rho$ bas{\'e} en $\cS_0$}.
Notons qu'on a $\hg_\rho(\cS_0) = -\rho (\berp)$ et que $\hg_{[\cS_0]}$ est la
fonction constante {\'e}gale {\`a} $-1$ sur tout $\berp$. Plus
g{\'e}n{\'e}ralement, $\hg_{[\cS']}(\cS) = -1 - \langle \cS , \cS'
\rangle_{\cS_0}$.  On d{\'e}signe par $\cP$ l'ensemble de tous les
potentiels.  C'est un espace vectoriel qui contient toutes les
fonctions de la forme $\langle \cdot \, , \cS' \rangle_{\cS_0}$.

Il r{\'e}sulte de~\cite[Th{\'e}or{\`e}me~7.50]{FJ} que
l'application $\rho \mapsto \hg_\rho$ induit une bijection entre $\cM$ et
$\cP$.  On peut donc poser
$$
\Delta \hg_\rho = \rho - \rho(\berp) \cdot [\cS_0]~.
$$ Ceci d{\'e}finit une application lin{\'e}aire $\Delta : \cP \to \cM$
que l'on appelle {\it le Laplacien}. Par construction, pour tout $g\in
\cP$ on~a $\Delta g ( \berp) =0$. R{\'e}ciproquement, toute mesure
v{\'e}rifiant $\rho (\berp ) =0$ est le Laplacien d'une fonction de
$\cP$. Dans toute la suite, on appellera \emph{potentiel} d'une mesure
bor{\'e}lienne $\rho$ toute fonction $g \in \cP$ telle que $\rho =
\Delta g$.
Notons que l'on a 
$$
\Delta \langle \cdot \, , \cS \rangle_{\cS_0} = [\cS_0] - [\cS]~.
$$ 
Le Lemme~\ref{L-infini comme point base} implique donc
\begin{equation}\label{e-330}
\Delta \log \sup \{ \cdot, \cS \} = [\cS]- [\infty]~.
\end{equation}
Pour tout $g \in \cP$, la fonction $\hg_{\Delta g}- g$ est constante.
On en d{\'e}duit que pour deux potentiels $g$ et $g'$, on~a $\Delta g =
\Delta g'$ si et seulement si la fonction $g - g'$ est constante.
\begin{proposition}
L'espace des potentiels $\cP$ et le Laplacien $\Delta : \cP \to \cM$
ne d{\'e}pendent pas du choix du point base.  De plus, la diff{\'e}rence
des potentiels d'une mesure, pris par rapport {\`a} des points base
distincts, s'{\'e}tend en une fonction d{\'e}finie et continue sur
$\berp$.
\end{proposition}
\begin{proof}
Choisissons $\cS_0, \cS_1 \in \H_p$ deux points base.  Etant donn{\'e}
une mesure $\rho \in \cM$, on d{\'e}signe par $g_0$ et $g_1$ les
potentiels de $\rho$ bas{\'e}s en $\cS_0$ et $\cS_1$ respectivement.
Lorsqu'on int{\`e}gre~\eqref{changement point
base} contre la mesure $\rho$, on obtient
\begin{equation}\label{e-230}
g_0(\cS)  
- 
C \, \langle \cS , \cS_1 \rangle_{\cS_0}
=
g_1(\cS) - g_1(\cS_0)~,
\end{equation}
avec $ C =\rho (\berp)$.  Par cons{\'e}quent, la diff{\'e}rence $g_0 - g_1$
s'{\'e}tend en une fonction d{\'e}finie et continue sur $\berp$.  Pour
montrer que l'espace des potentiels ne d{\'e}pend pas du point base, il
suffit de montrer que la fonction $g_1$ ci-dessus est un potentiel
lorsque $\cS_0$ est choisi comme point base. Ceci r{\'e}sulte
imm{\'e}diatemment de~\eqref{e-230}, car les deux fonctions $g_0$ et
$\langle \cdot , \cS_1\rangle_{\cS_0}$ sont des potentiels.  Enfin si
$\Delta\!^0$ et $\Delta\!^1$ d{\'e}notent les op{\'e}rateurs de Laplace avec
$\cS_0$ et $\cS_1$ pour point base respectivement, on~a
\begin{multline*}
\Delta\!^1 g_0 = \Delta\!^1 g_1 +C\,  \Delta\!^1 \langle \cdot ,
\cS_1\rangle_{\cS_0}
=\\ =
\rho - C\, [\cS_1] +
C\, \Delta\!^1  \langle \cdot ,
\cS_0\rangle_{\cS_1}  
= 
\rho - C\, [\cS_0] = \Delta\!^0 g_0~,
\end{multline*}
ce qui termine la preuve.
\end{proof}
%%%%%%%%%%%%%%%%%%%%%%%%%%%%%%%%%%%%%%%%%%%%%%%%%%%%%%%%%%%%%%%%%%%%%
%
%
\subsection{Le cas des potentiels {\`a} support fini.}
Bien que cel{\`a} ne soit pas pas strictement n{\'e}cessaire pour la
suite, nous allons indiquer ici de mani{\`e}re plus concr{\`e}te le
fonctionnement de l'op{\'e}rateur de Laplace sur une classe
particuli{\`e}re de fonctions.

Fixons $\cT$ un sous-arbre de $\berp$ donn{\'e} comme enveloppe convexe
d'un nom\-bre fini de points de $\H_p$. Cet arbre poss{\`e}de un nombre
fini de points de branchements, et on l'{\'e}crit comme r{\'e}union de
segments ferm{\'e}s $\cT = I_1 \cup \cdots \cup I_n$ de telle sorte que
pour tout $k\neq l$, $I_k \cap I_l$ est soit vide, soit r{\'e}duit {\`a} un point. 

Soit maintenant $g$ une fonction sur $\berp$ satisfaisant aux
propri{\'e}t{\'e}s suivantes:
\begin{itemize}
\item
$g$ est localement constante hors de $\cT$;
\item
la restriction de $g$ sur $\cT$ est continue;
\item la restriction de $g$ sur chaque segment $I_k$ est de classe
  $\cC^2$.
\end{itemize}
Il n'est pas difficile de montrer que $g$ est un potentiel au sens
pr{\'e}c{\'e}dent. La mesure $\Delta g$ est alors une combinaison de deux
termes. La restriction sur l'int{\'e}rieur des segments est le laplacien
standard au sens r{\'e}el (c'est-{\`a}-dire qu'il est donn{\'e} par la
d{\'e}riv{\'e}e seconde); aux points de branchement la masse de $\Delta g$
se calcule de mani{\`e}re combinatoire.

De mani{\`e}re pr{\'e}cise, sur chaque segment $I_k$, fixons un point extr{\'e}mal
$\cS_k$, et notons $\a_k : I_k \to \R_+$ la fonction $\a_k(\cS) =
d(\cS, \cS_k)$.  Cette fonction induit une isom{\'e}trie de $I_k$ sur son
image, et on note $\a_k \mapsto \cS_{\a_k}$ son inverse.  Alors la
restriction de $\Delta g$ sur l'int{\'e}rieur du segment $I_k$ est donn{\'e}e
par $-\frac{d^2}{d\a_k^2} g (\cS_{\a_k})$.  Fixons maintenant un point
extr{\'e}mal $\cS_*$ d'un des segments recouvrant $\cT$. La d{\'e}riv{\'e}e
sortante de $g$ en $\cS_*$ le long d'un segment $I_k$ contenant
$\cS_*$, est par d{\'e}finition le nombre $ \left.\frac{d}{d\a_k}\right|_0
g(\cS_{\a_k})$ si $\a_k(\cS_*)=0$, et
$\left. -\frac{d}{d\a_k}\right|_{\a_k \to \a_k(\cS)} g(\cS_{\a_k})$ sinon. On
note ce nombre $D_{I_k}g (\cS_*)$. 
On a alors
$$
\Delta g \{ \cS_* \} = \sum_{k,\, \cS \in I_k} 
D_{I_k} \, g(\cS_*)~.
$$

%%%%%%%%%%%%%%%%%%%%%%%%%%%%%%%%%%%%%%%%%%%%%%%%%%%%%%%%

\subsection{R{\'e}gularit{\'e} des potentiels.}\label{S-reg-Cp}
Fixons un point base $\cS_0 \in \H_p$.  On dira qu'une mesure de
probabilit{\'e} $\rho$ est {\`a} potentiel born{\'e} (resp. continu) si $\rho
-\lambda_p =\Delta g$ pour une fonction $g$ uniform{\'e}ment born{\'e} sur
$\H_p$ (resp.  d{\'e}finie {\`a} valeurs dans $\R$ et continue sur $\berp$).
De mani{\`e}re {\'e}quivalente, $\rho$ est {\`a} potentiel born{\'e}
(resp. continue) si il existe un point $\cS_0$ tel que le potentiel de
$\rho$ bas{\'e} en $\cS_0$ est born{\'e} (resp. continue). 
D'apr{\`e}s la proposition pr{\'e}c{\'e}dente, ces propri{\'e}t{\'e}s ne d{\'e}pendent pas du
choix du point base (bien que le potentiel lui-m{\^e}me en d{\'e}pende).
Toute combinaison lin{\'e}aire de masses de Dirac situ{\'e}es en des points de
$\H_p$ est {\`a} potentiel continu (donc born{\'e}).  Par ailleurs, le
potentiel d'une masse de Dirac en un point $z \in \pp$ tend vers
l'infini en ce point.  On en d{\'e}duit qu'une mesure positive {\`a} potentiel
born{\'e} ne charge pas les points de $\pp$.

L'analogue du Lemme~\ref{L-critint} d{\'e}montr{\'e} dans le cas
complexe s'{\'e}nonce de la mani{\`e}\-re suivante.
\begin{lemme}\label{L-int-Cp}
Soit $\rho$ une mesure positive finie sur $\berp$ {\`a} potentiel
born{\'e}.
Alors toute fonction dans $\cP$ est int{\'e}grable par rapport
{\`a} $\rho$.
\end{lemme}
\begin{proof}
Notons $g= -\hg_\rho$.

On d{\'e}montre tout d'abord que pour tout $\cS\in \berp$, on~a $g_\cS\=
\langle \cdot\, , \cS\rangle_{\cS_0} \in L^1(\rho)$. Lorsque $\cS \in
\H_p$, ceci est clair car $g_\cS$ est born{\'e}e et  par
d{\'e}finition on~a  $\int_\berp g_\cS \, d\rho = g(\cS) +C$ avec $C = - \rho
(\berp)$.
On peut donc supposer que $\cS = z \in \pp$. Pour tout
$t \ge 0$, on note $\cS_t$ l'unique point du segment $(z, \cS_0]$ {\`a}
distance $t$ de $\cS_0$.  Pour tout $t \ge 0$, on~a $\langle \cdot \,
, \cS_t \rangle_{\cS_0} = \min \{ t, \langle \cdot \, , z
\rangle_{\cS_0} \}$, donc la suite de fonctions $g_{\cS_t}$ converge
en croissant vers $g_z$ lorsque $t$ tend vers l'infini.  Par
convergence monotone, on obtient
\begin{equation*} 
0 \le \int g_z d \rho
=
\lim_{t \to \infty} \int g_{\cS_t} d\rho
=
\lim_{t \to \infty}  g(\cS_t) +C~.
\end{equation*}
Mais le potentiel $g$ de $\rho$ est born{\'e} uniform{\'e}ment dans $\H_p$ par
hypoth{\`e}se, donc $g_z$ est bien int{\'e}grable par rapport {\`a} $\rho$, et
on a $0\le \int g_z d \rho \le \sup | g| + |C|$ quantit{\'e} ind{\'e}pendante
de $z$. Il est alors facile d'en d{\'e}duire par int{\'e}gration que $g' \in
L^1(\rho)$ pour tout potentiel $g' \in \cP$.
\end{proof}

%%%%%%%%%%%%%%%%%%%%%%%%%%%%%%%%%%%%%%%%%%%%%%%%%%%%%%%%

\subsection{Energie.}\label{capacite p-adique}

Notons $\diag = \{ (z, z) , \; z \in \C_p \}$ la diagonale de $\C_p$
dans $\berCp \times \berCp$.  Pour chaque paire de mesures $\rho,
\rho' \in \cM$, telle que la fonction $\log \sup \{ \cdot, \cdot \}$
soit int{\'e}grable par rapport {\`a} $\rho \otimes \rho'$ dans
$\berCp\times \berCp \setminus \diag$, on pose
$$
(\rho, \rho')
=
- \int_{\berCp \times \berCp \setminus \diag} \log \sup \{
\cS, \cS' \} \ d\rho(\cS) \otimes d\rho'(\cS')~.
$$
Comme dans le cas complexe, lorsque les mesure sont atomiques $\rho
= \sum m_i [\cS_i]$, $\rho' = \sum m_j' [\cS'_j]$, l'hypoth{\`e}se
d'int{\'e}\-gra\-bi\-lit{\'e} est automatiquement satisfaite et on~a
$(\rho, \rho')= -\sum m_i m_j' \log \sup \{ \cS_i , \cS'_j\}$, la somme
{\'e}tant prise sur les indices tels que $\cS_i\ne \cS_j$ et
$\cS_i,\cS_j \ne \infty$. De m{\^e}me pour tout $\cS \in \berCp$, on~a
$$
([\cS], [\cS_\can]) = - \log \sup \{ \cS, \cS_\can \} = - \log^+|\cS|~.
$$ Gr{\^a}ce au Lemme~\ref{L-int-Cp} ci-dessus, la preuve du
Lemme~\ref{L-integr} se recopie mot {\`a} mot dans le cas $p$-adique et
on a donc:
\begin{lemme}\label{L-integr-Cp}
Soit $\rho$ une mesure sign{\'e}e dont la mesure trace est {\`a} potentiel
continu.  De plus soit $\rho'$ une mesure satisfaisant l'une des
propri{\'e}t{\'e}s suivantes~:
\begin{itemize}
\item[$\bullet$]
$\rho'$ est une mesure {\`a} support fini ne chargeant pas l'infini~;
\item[$\bullet$]
$|\rho'|$ est {\`a} potentiel continu.
\end{itemize}
Alors $\log \sup \{ \cdot , \, \cdot \} \in L^1( |\rho| \otimes
|\rho'|)$ dans $\berp \times \berp$.  En particulier, $(\rho, \rho')$
est bien d{\'e}fini.
\end{lemme}
De m{\^e}me la preuve du Lemme~\ref{L-cal} s'adapte au cas $p$-adique
pour donner le lemme suivant.
\begin{lemme}\label{L-calculCp}
  Soient $\rho,\rho'$ deux mesures sign{\'e}es telles $\log
  \sup \{ \cdot, \cdot \} \in L^1( |\rho| \otimes |\rho'| )$ dans
  $\berp \times \berp$.  
  
  Alors la fonction $g_\rho (\cS) \= \int_{\berCp} \log \sup \{ \cS,
  \cS'\} \, d\rho (\cS')$ est int{\'e}grable par rapport {\`a} $\rho'$ et on~a
\begin{equation}\label{e-500}
(\rho,\rho') = - \int_{\berCp} g_\rho \, d\rho' ~.
\end{equation}
\end{lemme}

%%%%%%%%%%%%%%%%%%%%%%%%%%%%%%%%%%%%%%%%%%%%%%%%%%%%%%%%

\subsection{Positivit{\'e}.}
On cherche {\`a} d{\'e}montrer un analogue de la
Proposition~\ref{P-casdeg}.

Fixons un point base $\cS_0\in \H_p$.
On d{\'e}finit la relation d'ordre $\preccurlyeq$ par,
$\cS \preccurlyeq \cS'$ si et seulement si $\cS \in [\cS_0, \cS']$.
Etant donn{\'e}e une mesure bor{\'e}lienne $\rho$, on d{\'e}signe par $f_\rho :
\berp \to \R$ la fonction nulle sur $\pp$ et telle que pour $\cS \in
\H_p$ on ait $f_\rho(\cS) = \rho \{ \cS \preccurlyeq \cdot \}$.  On
introduit aussi $\lambda$ la mesure sur $\berp$ qui est nulle sur
$\pp$ et qui co{\"\i}ncide avec la mesure de Hausdorff de dimension~1
sur $\H_p$, par rapport {\`a} la distance~$d$. En particulier, $\lambda
[\cS, \cS'] = d( \cS, \cS')$ pour tout $\cS, \cS' \in \H_p$.
\begin{proposition}\label{P-casdeg-Cp}
Soit $\rho$ une mesure sign{\'e}e sur $\pp$ telle que $\rho(\berp) =0$,
et dont la mesure trace est {\`a} potentiel continu.  Alors $(\rho,
\rho)$ est bien d{\'e}finie, la fonction $f_\rho$ d{\'e}finie ci-dessus
est dans $L^2(\lambda)$ et on~a
\begin{equation}\label{e-440}
(\rho, \rho ) =  \int_{\H_p} f_\rho^2 \, d\lambda \ge 0~.
\end{equation}
De plus $(\rho, \rho ) = 0$ si et seulement si $\rho = 0$.
\end{proposition}
Remarquons que cette {\'e}quation indique que l'{\'e}nergie d'une mesure
de masse totale nulle, est invariante par le groupe des automorphismes
de $\pp$.
Ainsi, lorsque $\rho = [\cS] - [\cS']$ avec $\cS, \cS' \in
\H_p$ on v{\'e}rifie que $(\rho,\rho) = d(\cS, \cS')$ (en utilisant
directement la d{\'e}finition ou en appliquant~\eqref{e-440}).

La d{\'e}monstration de cette proposition d{\'e}pend du lemme suivant.

\begin{lemme}\label{L-270}
Soit $\rho$ une mesure positive dont la mesure trace est {\`a} potentiel
continu.  Notons $f_\rho(\cS) = \rho \{ \cS \preccurlyeq \cdot \}$ et
$\tilde{g}_\rho(\cS) = \int_{[\cS_0, \cS]} f_\rho \, d\lambda$.  Alors
$\tilde{g}_\rho$ appartient {\`a} $\cP$, on~a $\Delta \tilde{g}_\rho = -
\rho + \rho(\berp) [\cS_0]$ et
\begin{equation}\label{e-270}
\int_{\H_p} f^2_\rho \, d\lambda = \int_\berp \tilde{g}_{\rho} \, d\rho ~.
\end{equation}
Plus g{\'e}n{\'e}ralement, cette {\'e}quation est valide pour toute mesure
bor{\'e}lienne sign{\'e}e d{\`e}s que $f_\rho \in L^2(\rho)$.
\end{lemme}

\begin{proof}[Preuve de la Proposition~\ref{P-casdeg-Cp}]
Par hypoth{\`e}se, $|\rho|$ est {\`a} potentiel continu.  Par le lem\-me
pr{\'e}c{\'e}dent, on d{\'e}duit que $\tilde{g}_{|\rho|}$ est un potentiel
pour $-|\rho|$ et qu'il est donc d{\'e}fini et continu dans $\berp$.
En particulier, l'int{\'e}grale $\int_\berp \tilde{g}_{|\rho|} d|\rho|$ est
finie.
Donc $|f_\rho| \le f_{|\rho|} \in L^2(\lambda)$
par~\eqref{e-270}. Le lemme pr{\'e}c{\'e}dent s'applique alors {\`a} la
mesure sign{\'e}e $\rho$ et~\eqref{e-270} donne $\int \tilde{g}_\rho \,
d\rho = \int f_\rho^2 d\lambda \ge0$.

Pour d{\'e}montrer~\eqref{e-440}, on va relier $\tilde{g}_\rho$ {\`a}
la fonction $g_\rho$, d{\'e}finie dans le Lemme~\ref{L-calculCp}.
En int{\'e}grant~\eqref{e-330} et utilisant l'hypoth{\`e}se
$\rho(\berp ) = 0$, on obtient $\Delta g_\rho = \int_\berp ( [\cS] - [\infty]
) \, d\rho(\cS) =  \rho$. Donc $g_\rho$ diff{\`e}re de la
fonction $- \tilde{g}_\rho$ d'une constante et de~\eqref{e-500} on tire
finalement
$$ (\rho, \rho) = - \int g_\rho \, d\rho = \int \tilde{g}_\rho \,
d\rho= \int f_\rho^2 d\lambda \ge0~.$$ Enfin si $(\rho,\rho)=0$, alors
$f_\rho= 0$ pour $\lambda$-presque tout point. Pour tout $\cS \in
\H_p$, on peut trouver une suite $\cS_k \to \cS$ croissante pour
$\preccurlyeq$ et telle que $f_\rho(\cS_k) =0$.  Mais $f_\rho(\cS_k) =
\rho \{ \cS_k \preccurlyeq \cdot \} \to \rho \{ \cS \preccurlyeq \cdot
\} = f_\rho(\cS)$, donc $f_\rho = 0$ pour tout point de $\H_p$.  On
remarque maintenant que tout ensemble de la forme $\{ \cS \preccurlyeq
\cdot \}$ est une boule ferm{\'e}e de $\berp$.  R{\'e}ciproquement, toute
boule de $\berp$ est de cette forme quitte {\`a} faire varier le point
base.  On en d{\'e}duit finalement que $\rho(B)=0$ pour toute boule
ferm{\'e}e de $\berp$.  L'ensemble des boules forme une base de la
topologie de $\berp$ et $\rho$ est r{\'e}guli{\`e}re, donc $\rho =0$.
\end{proof}
\begin{proof}[Preuve du Lemme~\ref{L-270}]
Pour $\cS\in \H_p$, notons $T \= d(\cS_0, \cS)$ et soit $\cS_t$
l'unique point du segment $[\cS_0, \cS]$ {\`a} distance $t$ de $\cS_0$.
Notons que dans ce cas $\langle \cS , \cS'\rangle_{\cS_0} \le T$ pour
tout $\cS' \in \berp$. On a alors l'{\'e}quation suivante:
\begin{multline*}
\int_\berp \langle \cS , \cS'\rangle_{\cS_0} \, d\rho (\cS')
=
\int_0^T \rho \{  \langle \cS , \cdot \, \rangle_{\cS_0} \ge t \} \, dt
=
\int_0^T \rho \{ \cS_t \preccurlyeq \cdot \, \} \, dt
=\\=
\int_0^T f_\rho (\cS_t)\, dt 
= 
\int_{[\cS_0, \cS]} f_\rho\, d\lambda
=
\tilde{g}_\rho(\cS)~.
\end{multline*}
Par d{\'e}finition du laplacien on~a $$\Delta \tilde{g}_\rho (\cS) =
\Delta \left( \int_\berp \langle \cS , \cdot\rangle \, d\rho \right)=
- \rho + \rho(\berp)[\cS_0]~.$$

Pour d{\'e}montrer~\eqref{e-270}, on {\'e}crit les {\'e}galit{\'e}s suivantes:
\begin{multline*}
\int_\berp \tilde{g}_{\rho} \, d\rho
=
\int_\berp \left[ \int_{[\cS_0, \cS]} f_\rho \, d\lambda \right] d\rho(\cS)
=\\=
\int_{\cS' \preccurlyeq \cS, \cS' \in \H_p}
f_\rho(\cS') \, d\lambda(\cS')\otimes d\rho(\cS)
=\\=
\int_{\cS'\in \H_p}
f_\rho(\cS') \times \rho \{ \cS' \preccurlyeq \cS\} d\lambda(\cS')
=
\int_{\H_p} f_\rho^2 d\lambda~.
\end{multline*}
Lorsque $\rho$ est sign{\'e}e et que $f_\rho\in L^2(\lambda)$, on
v{\'e}rifie que tous les termes ont un sens et que donc
l'{\'e}galit{\'e}~\eqref{e-270} reste valide. Ceci termine la preuve du
lemme.
\end{proof}

%%%%%%%%%%%%%%%%%%%%%%%%%%%%%%%%%%%%%%%%%%%%%%%%%%%%%%%%

\subsection{R{\'e}gularisation des mesures}\label{S-regularisation-Cp}
Comme dans le cas complexe, la r{\'e}gularisation des mesures est
fondamentale dans notre approche. Dans le cas pr{\'e}sent, l'op{\'e}ration
de convolution est remplac{\'e}e par une proc{\'e}dure de projection que
nous {\'e}tudions tout d'abord.

Etant donn{\'e} $\varepsilon \ge 0$, on d{\'e}signe par $\pi_\varepsilon :
\berp \to \berp$ l'application envoyant $\cS$ sur l'unique
point $\cS'\in [\cS, \infty]$ tel que $\diam(\cS') = \max\{
\diam(\cS), \varepsilon \}$. 
Sur l'ensemble
$$ \He = \{ \cS \in \berp , \,  \diam(\cS) \ge \varepsilon \}
$$ on~a $\pi_\varepsilon = \mathrm{id}$, tandis que $\pie(\berp) =
\He$.  Il est clair que $H_0 = \berp$ et que $\pi_0 = \mathrm{id}$ sur
$\berp$.  Enfin pour tout $\varepsilon, \varepsilon' \ge 0$, on~a
$\pi_\varepsilon \circ \pi_{\varepsilon'} = \pi_{\max \{ \varepsilon,
\varepsilon' \}}$.
\begin{lemme}\label{projectionepsilon}
Pour tout $\varepsilon \ge 0$ et toute boule $B$ de $\berp$, on~a
$$
\pi_\varepsilon^{-1}(B)
=
\begin{cases}
\emptyset &
\text{si } B \cap H_\varepsilon = \emptyset~; \\
\berp &
\text{si }H_\varepsilon \subset B~;\\
B &
\text{si }B \cap H_\varepsilon \neq \emptyset 
\text{ et }H_\varepsilon \not\subset B ~.
\end{cases}
$$ En particulier la fonction $\pie: \berp \to \berp$ est faiblement
continue.
\end{lemme}
\begin{proof}
Consid{\'e}rons la boule $B' = \berp \setminus B$ et notons que les ensembles
$\pie^{-1}(B)$ et $\pie^{-1}(B')$ forment une partition de $\berp$.
Lorsque $B$ (resp. $B'$) est disjointe de $\He =
\pie(\berp)$, on~a $\pie^{-1}(B) = \emptyset$ (resp. $\pie^{-1}(B') =
\emptyset$).  Il suffit donc de montrer que lorsque $B$ (resp. $B'$)
rencontre $\He$, on~a $B \subset \pie^{-1}(B)$ (resp. $B' \subset
\pie^{-1}(B')$).  On se ram{\`e}ne au cas de la boule $B$.

Pour tout $\cS \in \berp$ et tout $\cS_0 \in \He$, le point
$\pie(\cS)$ appartient {\`a} $[\cS, \cS_0]$.  Comme la boule $B$ est
connexe, lorsque $\cS$ appartient {\`a} $B$ et $\cS_0$ appartient {\`a} $B
\cap \He$, on~a $\pie(\cS) \in [\cS, \cS_0] \subset B$.  On a donc
$\cS \in \pie^{-1}(B)$.
\end{proof}

A l'aide de cette projection, on d{\'e}finit la r{\'e}gularis{\'e}e d'une
mesure bor{\'e}lienne $\rho$ sur $\berp$, en posant $\rho_\e =
\pi_{\varepsilon *}\rho$ et on v{\'e}rifie facilement le lemme suivant.
\begin{lemme}\label{L-reg-Cp}
  Pour toute mesure sign{\'e}e $\rho$ sur $\berp$, on~a $\rho_\e = \rho$
  lorsque $\varepsilon \to 0$. Pour toute suite de mesures sign{\'e}es
  $\rho_n\to\rho$, on~a $\rho_{n,\e}\to\ \rho_\e$. 
  
  Enfin, si $F\subset \berCp$ est un ensemble fini, notons $[F] =
  |F|^{-1}\sum_{z \in F} [z]$ la mesure {\'e}quidistri\-bu{\'e}e sur les
  points de $F$.  Alors pour tout $\e > 0$ la mesure de probabilit{\'e}
  $[F]_\e$ est {\`a} potentiel continu.
\end{lemme}
Plus g{\'e}n{\'e}ralement, on peut montrer que pour toute mesure $\rho$
dont le support {\'e}vite le point $\infty$, les mesures $\rho_\e$ sont
{\`a} potentiel born{\'e}.

%%%%%%%%%%%%%%%%%%%%%%%%%%%%%%%%%%%%%%%%%%%%%%%%%%%%%%%%

\subsection{Energie et r{\'e}gularisation.}
Ce paragraphe est consacr{\'e} {\`a} la preuve du r{\'e}sultat suivant.  Rappelons
qu'un \emph{module de continuit{\'e}} pour une fonction continue $h$ sur
$\berp$ est une fonction $\eta: \R_+ \to \R_+$, telle que pour tous
points $\cS, \cS' \in \berp$ tels que $\mathsf{d}(\cS,\cS')\le \e$, on
ait $|h(\cS) - h(\cS')|\le \eta(\e)$.  Ici $\mathsf{d}$ d{\'e}note la
m{\'e}trique chordale sur $\berp$, qui peut {\^e}tre d{\'e}finie par
$$
\mathsf{d}(\cS, \cS') 
\= 
\frac{\sup\{\cS, \cS'\}}{\max \{1, |\cS|\} \max \{1, |\cS'|\}} 
- \frac{\diam(\cS)}{2\max \{1, |\cS|\}^2} 
- \frac{\diam(\cS')}{2\max \{1, |\cS'|\}^2}~.
$$
Notons que pour tous $\cS, \cS'\in \berp$ on~a $\mathsf{d}(\cS,
\cS') \le \sup\{\cS, \cS' \}$.
\begin{proposition}\label{P-230}
  Soit $\rho$ une mesure de probabilit{\'e} {\`a} potentiel continu et soit
  $\eta$ un module de continuit{\'e} d'un potentiel de $\rho - [S_\can]$.
  Alors pour tout $\e \in (0, 1)$ et tout sous ensemble fini $F$ de
  $\berCp$ on~a
\begin{eqnarray}
([F] - \rho , [F] -\rho) 
& \ge & ([F]_\e - \rho , [F]_\e -\rho) 
- 2\eta(\e) - |F|^{-1} \log \e^{-1} \label{e-290}\\
& \ge &
- 2 \eta(\e) - |F|^{-1} \log \e^{-1} \label{e-291}
\end{eqnarray}
\end{proposition}
Comme dans le cas complexe, la preuve repose sur les deux lemmes suivants.
\begin{lemme}\label{L-300}
  Soit $\rho$ une mesure de probabilit{\'e} sur $\berv$ {\`a} potentiel
  continu et soit $\eta$ un module de continuit{\'e} d'un potentiel de
  $\rho - [\cS_\can]$.  Alors pour tout $\e \in (0, 1)$ et pour tout
  sous ensemble fini $F$ de $\berCp$ on~a
\begin{equation*}%\label{e-350}
|([F], \rho ) - ([F]_\e , \rho )|\le \eta(\e)~.
\end{equation*}
\end{lemme}
\begin{proof}
On peut supposer que $F$ est r{\'e}duit {\`a} un point $z\in \C_p$.  Pour
$\varepsilon \ge 0$, on d{\'e}signe par $z(\varepsilon)$ l'unique point
de $\berCp$ tel que $z(\varepsilon) \ge z$ et $\diam(z(\varepsilon)) =
\varepsilon$~; on~a par d{\'e}finition $z(0) = z$ et $[z]_\e = [z(\e)]$
(cf.~\S\ref{S-regularisation-Cp}).

Notons $g_\rho (\cS) = \int_\berp \log \sup \{ \cS, \cS'\} \, d\rho
(\cS')$.  Le Lemme~\ref{L-integr-Cp} implique que la fonction $\log
\sup \{ \cdot, \cdot \}$ appartient {\`a} $L^1(|\rho| \otimes [z])$ et
{\`a} $L^1(|\rho| \otimes [z(\varepsilon)])$ et le
Lemme~\ref{L-calculCp} implique que $([z],\rho) - ([z(\e)], \rho) = -
g_\rho(z) + g_\rho(z(\e))$.  On va maintenant relier $g_\rho$ {\`a} un
potentiel de $\rho$.

Pour ce faire, on int{\`e}gre l'{\'e}quation~\eqref{e-330}.  On obtient
$\Delta g_\rho = \rho - [\infty]$.  Fixons comme point base le point
$\cS_\can \in \H_p$ et remarquons que $g_{[\cS_\can]} (\cS) = \log^+ |
\cS|$. Comme $\Delta g_{[\cS_\can]} = [\cS_\can] - [\infty]$, on~a
finalement $\Delta ( g_\rho- g_{[\cS_\can]} ) = \rho - [\cS_\can]$. La
fonction $h \= g_\rho - g_{[\cS_\can]}$ d{\'e}finit donc un potentiel de
$\rho$.  Par hypoth{\`e}se c'est une fonction continue ayant $\eta$ comme
module de continuit{\'e}.  Alors pour tous $z, z' \in \berCp$ distincts
on a $ |h(z) - h(z')| \le \eta(\sup\{ z, z' \}) $.  Pour tout $\e>0$,
on a $\sup \{ z, z(\e) \} = z(\e)$ et lorsque $\e \in (0, 1)$ on~a
$\log^+|z| = \log^+|z(\e)|$, donc
\begin{multline*}
|([z],\rho) - ([z]_\e, \rho)|
=
|g_\rho(z) - g_\rho(z(\e))|
\le \\ \le
|h(z) - h(z(\e))|+ \left| \log^+ |z| - \log^+ |z(\e)| \right|
\le 
\eta(\e)~.
\end{multline*}
\end{proof}
\begin{lemme}\label{L-320}
Pour tout $\e>0$ et tout sous-ensemble fini $F$ de $\berCp$, on~a
\begin{equation*}
([F]_\e , [F]_\e ) \le ([F], [F]) + |F|^{-1} \log \e^{-1}~.
\end{equation*}
\end{lemme}
\begin{proof}
Pour tout $z \ne  z' \in \C_p$, on~a $\sup \{ z(\varepsilon),
z'(\varepsilon) \} \ge \sup \{ z, z' \}$.  Comme $\sup \{
z(\varepsilon), z(\varepsilon) \} = \varepsilon$ on obtient
\begin{eqnarray*}
([F]_\varepsilon, [F]_\varepsilon)
& = & 
- |F|^{-2} \sum_{z, z' \in F} 
\log (\sup \{ z(\varepsilon), z'(\varepsilon) \}) \\
& & \le
([F], [F]) + |F|^{-1} \log \e^{-1}~.
\end{eqnarray*}
\end{proof}
\begin{proof}[Preuve de la Proposition~\ref{P-230}]
Comme $[F]_\varepsilon$ est une mesure de probabilit{\'e} qui est une
combinaison lin{\'e}aire de masses de Dirac situ{\'e}es en des points de
$\H_p$, elle admet un potentiel d{\'e}fini et continu sur $\berp$.  Le
Lemme~\ref{P-casdeg-Cp} implique alors que l'{\'e}nergie
$([F]_\varepsilon -\rho, [F]_\varepsilon -\rho)$ est bien d{\'e}finie et
qu'on a $([F]_\varepsilon -\rho, [F]_\varepsilon -\rho) \ge 0$.

Le r{\'e}sultat est alors une cons{\'e}quence imm{\'e}diate des
Lemmes~\ref{L-300} et~\ref{L-320}.
\end{proof}

%%%%%%%%%%%%%%%%%%%%%%%%%%%%%%%%%%%%%%%%%%%%%%%%%%%%%%%%

\subsection{Le r{\'e}sultat clef.}
\begin{proposition}\label{critere de convergence}
Soit $\rho$ une mesure de probabilit{\'e} admettant un potentiel
d{\'e}fini et continu sur $\berp$.  Soit $\{ F_n \}_{n \ge 0}$ une suite
de sous-ensembles finis de $\pp$, tels que $|F_n| \to \infty$. Alors
$$
\ov{\lim}_{n\to \infty} ([F_n] - \rho, [F_n] - \rho ) \le 0 
\, \text{ implique } \,
\lim_{n\to \infty} [F_n] = \rho~.
$$
\end{proposition}
\begin{proof}
  Comme dans la preuve de la Proposition~\ref{P-casdeg2}, on peut
  supposer que $F_n\subset \C_p$ pour tout $n$.

Quitte {\`a} prendre une sous-suite, on suppose de plus que $[F_n]$
converge vaguement vers une mesure $\rho'$.  C'est une mesure de
probabilit{\'e} sur $\berp$.  Com\-me $\pi_\e$ est faiblement continue, on
a $\lim_{n \to \infty} [F_n]_\e = \rho'_\e$ et $\lim_{\e \to 0}
\lim_{n \to \infty} [F_n]_\e = \rho'$.

On va montrer que pour chaque boule $B$ ouverte de $\berp$ on~a
$\rho'(B) \ge \rho(B)$.  Comme $\rho$ et $\rho'$ sont des mesures de
probabilit{\'e}, ceci implique que $\rho' = \rho$.

Fixons alors une boule $B$ ouverte de $\berp$. Choisissons de plus
$\cS_0\in\H_p$ un point base et notons comme pr{\'e}c{\'e}demment
$\preccurlyeq$ la relation d'ordre telle que $\cS \preccurlyeq \cS'$
si et seulement si $\cS \in [\cS_0 ,\cS']$.  On choisit $\cS_0$ dans
le compl{\'e}mentaire de $\overline{B}$, de telle fa{\c c}on qu'il existe un
point $\cS$ distinct de $\cS_0$ tel que $B = \{ \cS \prec \cdot
\}$.  Pour tout $n,\e$, on introduit alors la fonction $f_{n,\e} (
\cS) \= ( [F_n]_\e - \rho ) \{ \cS \preccurlyeq \cdot \}$.

Les mesures $[F_n]_\e$ sont {\`a} support fini inclus dans $\H_p$, donc {\`a}
potentiel continu.  On peut donc appliquer la
Proposition~\ref{P-casdeg-Cp} et on obtient
$$
([F_n]_\e - \rho, [F_n]_\e -\rho ) = \int_{\H_p} f_{n,\e}^2\,
d\lambda \ge 0 ~.
$$
De~\eqref{e-290}, on tire
\begin{equation*}
\int_{\H_p} f_{n,\e}^2 \, d \lambda 
= 
([F_n]_\e - \rho, [F_n]_\e -\rho ) 
\le
([F_n] - \rho , [F_n] - \rho ) + \eta (\e) +  |F_n|^{-1} \log  \e^{-1} ~.
\end{equation*}
Comme $|F_n|\to \infty$ et que par hypoth{\`e}se
$\ov{\lim}_{n \to \infty} ([F_n]-\rho,[F_n]-\rho) \le 0$, on en d{\'e}duit que
\begin{equation}\label{e-318}
 \blim_{\e \to 0} \blim_{n\to \infty} \int_{\H_p} f_{n,\e}^2 \, d \lambda = 0
\end{equation}
En particulier sur tout segment de $\H_p$ on~a $f_{n,\e} \to 0$
lorsque $n \to \infty$ et $\e\to 0$, et ce, presque partout par
rapport {\`a} la mesure~$\lambda$.

Choisissons deux points
$\cS_0 \preccurlyeq \cS' \preccurlyeq \cS'' \preccurlyeq \cS$ tels que
$$
\lim_{\e \to 0} \lim_{n \to \infty} f_{n,\e} (\cS'')= 0 ~.
$$
On a alors
$$
[F_{n,\e}]\{\cS' \prec \cdot \}
\ge
[F_{n,\e}]\{\cS'' \preccurlyeq \cdot\}
=
\rho\{\cS'' \preccurlyeq \cdot\} + f_{n,\e}(\cS'')~,
$$ donc $\rho'\{\cS' \preccurlyeq \cdot\}\ge \blim_{\e \to 0} \blim_{n
\to \infty} [F_{n,\e}]\{\cS' \prec \cdot\} \ge \rho\{\cS'' \preccurlyeq
\cdot\}$.  On peut maintenant faire tendre $\cS'$ et $\cS''$ vers
$\cS$ et on obtient $\rho'(B) \ge \rho(B)$. Ce qui termine la preuve.
\end{proof}

%
%%%%%%%%%%%%%%%%%%%%%%%%%%%%%%%%%%%%%%%%%%%%%%%%%%%%%%%%%%%%%%%%%%
%

\section{Hauteurs ad{\'e}liques.}
\label{S-hauteur Arakelov}
Apr{\`e}s quelques g{\'e}n{\'e}ralit{\'e}s au~\S~\ref{S-generalites}, on introduit la
hauteur d{\'e}finie par une mesure ad{\'e}lique et on d{\'e}crit ses premi{\`e}res
propri{\'e}t{\'e}s au~\S~\ref{S-comp}.  On montre la formule de
Mahler au~\S~\ref{S-formule de Mahler} et l'{\'e}quidistribution des
points de petite hauteur au~\S~\ref{S-equidistribution}.  Enfin, on
montre la version quantitative de l'{\'e}quidistribution
au~\S~\ref{S-vitesse}, dans le cas o{\`u} la mesure ad{\'e}lique est {\`a}
potentiel H{\"o}lder.

%%%%%%%%%%%%%%%%%%%%%%%%%%%%%%%%%%%%%%%%%%%%%%%%%%%%%%%%%%%%

\subsection{G{\'e}n{\'e}ralit{\'e}s et notations.}\label{S-generalites}
Pour toute cette section, nous renvoyons
{\`a}~\cite[Part~B]{hindry-silverman} pour plus de pr{\'e}cisions.

Soit $M_\Q$ l'ensemble constitu{\'e} de tous les nombres entiers premiers
auquel on ajoute $\infty$.  On note $|\cdot|_\infty$ la norme usuelle
sur $\Q$, et pour chaque nombre premier $p$ on note $| \cdot |_p$ la
norme $p$-adique, normalis{\'e}e de telle sorte que $| p |_p = p^{-1}$.
Pour tout $\alpha\in \Q^*$, il existe au plus un nombre fini de $v \in
M_\Q$ tel que $|\alpha|_v \neq 1$, et on~a
$$
\prod_{v \in M_\Q} | \alpha |_v = 1~.
$$ Fixons maintenant une extension finie $K$ de $\Q$.  On d{\'e}signe
par $M_K$ la collection de toutes les normes sur $K$ qui {\'e}tendent
l'une des normes $| \cdot |_v$, avec $v \in M_\Q$.  Un {\'e}l{\'e}ment de $M_K$
est appel{\'e} \emph{place}. Pour toute place  $v$ de $K$, on notera encore $|
\cdot |_v$ la norme sur $K$ correspondante.  On dira que $v \in M_K$
est {\it infinie} lorsque $| \cdot |_v$ est une extension de la norme
$| \cdot |_\infty$, et que $v$ est {\it finie} sinon.  Pour toute place
$v \in M_\Q$, il n'existe qu'un nombre fini d'{\'e}l{\'e}ments de $M_K$
qui {\'e}tendent la norme $| \cdot |_v$.  En particulier, le nombre
d'{\'e}l{\'e}ments infinis de $M_K$ est fini.

Pour $v \in M_K$ fix{\'e}e, $K_v$ d{\'e}signera la compl{\'e}tion du corps
valu{\'e} $(K, | \cdot |_v)$, et $\Q_v$ le compl{\'e}t{\'e} de $\Q$ dans
$K_v$.  On pose alors $N_v = [K_v : \Q_v] / [K : \Q]$ et $ \| \cdot
\|_v = | \cdot |_v^{N_v}$.  On montre que pour tout $\alpha \in K^*$,
il n'existe qu'un nombre fini de $v \in M_K$ tel que $|\alpha|_v \neq
1$, et  on~a la {\it formule du produit}
$$
\prod_{v \in M_K} \| \alpha \|_v = 1~.
$$
Pour toute place $v \in M_K$, la norme $| \cdot |_v$ s'{\'e}tend de
mani{\`e}re unique {\`a} la cl{\^o}ture alg{\'e}brique $\ov{K}_v$ de $K_v$.  On
d{\'e}signera par $\C_v$ le compl{\'e}t{\'e} de $(\ov{K}_v, | \cdot |_v)$.
Lorsque $v$ est infinie, le corps $\C_v$ est isom{\'e}triquement
isomorphe au corps des nombres complexes $\C$.  Lorsque $v$ est finie,
la restriction de $| \cdot |_v$ {\`a} $\Q$ est une norme associ{\'e}e {\`a} un
nombre premier $p$, et $\C_v$ est isom{\'e}triquement isomorphe au corps
$\C_p$.

\medskip

Pour $v \in M_K$ finie, on d{\'e}signe par $\pv = \C_v \cup \{ \infty \}$
la droite projective correspondante et par $\berv$ l'espace des
semi-normes d{\'e}crit au~\S\ref{espace des semi-normes}.  Lorsque $v$
est infinie, on d{\'e}signera indiff{\'e}remment par $\P^1(\C_v)$ ou $\berv$ la
droite projective sur $\C_v$.

On note $\cM_v$ l'espace des mesures bor{\'e}liennes {\`a} support dans
$\berv$, et $(\cdot, \cdot)_v$ l'accouplement d{\'e}finie
au~\S\ref{S-energie} lorsque $v$ est infinie et au~\S\ref{capacite
p-adique} lorsque $v$ est finie.  Par commodit{\'e}, on pose
$$
\lpar \cdot, \cdot \rpar_v = N_v \, (\cdot, \cdot)_v~.
$$
Enfin, $\lambda_v$ d{\'e}signera la mesure de probabilit{\'e}
proportionnelle {\`a} la mesure de Lebesgue du cercle $S^1 \subset \C_v$
lorsque $v$ est infinie, et {\`a} la masse de Dirac situ{\'e}e au point
canonique $\cS_\can$ dans $\berv$ lorsque $v$ est finie.  Notons que
pour tout $v \in M_v$ et $\alpha \in \C_v$, on~a
$$
\log^+ |\alpha|_v = - ([\alpha], \lambda_v)_v
\ \mbox{ et } \
\log^+ \| \alpha \|_v = - \lpar [\alpha], \lambda_v\rpar_v~.
$$

%%%%%%%%%%%%%%%%%%%%%%%%%%%%%%%%%%%%%%%%%%%%%%%%%%%%%%%%%%%%

\subsection{Mesures ad{\'e}liques et hauteurs.}\label{S-comp}
On pose tout d'abord la d{\'e}finition suivante.
\begin{definition}
  Soit $K$ un corps de nombres. Une mesure ad{\'e}lique $\rho =\{
  \rho_v\}_{v \in M_K}$ est la donn{\'e}e pour chaque place $v\in M_K$
  d'une mesure de probabilit{\'e} $\rho_v$ d{\'e}finie sur $\berv$ et {\`a}
  potentiel continu, et telle que, hors un nombre fini de place, on ait
  $\rho_v = \lambda_v$.
\end{definition}
Etant donn{\'e}e une mesure ad{\'e}lique $\rho = \{ \rho_v \}_{v \in M_K}$, on
d{\'e}finit la fonction $h_\rho : \P^1(\ov{K}) \to \R$ de la mani{\`e}re
suivante.  Pour un sous ensemble fini $F$ de $\P^1(\ov{K})$ invariant
sous l'action de $\Gal(\ov{K}/K)$, on pose
\begin{equation}\label{hauteur arakelov}
h_\rho(F) \=
\frac12
\sum_{v \in M_K} \lpar [F] - \rho_v, [F] - \rho_v \rpar_v~,
\end{equation}
et pour $\a \in \P^1(\ov{K})$ on d{\'e}finit $h_\rho(\a) = h_\rho(F)$, o{\`u}
$F$ est l'orbite de $\a$ sous l'action du groupe de Galois
$\Gal(\ov{K}/K)$.  Notons que la somme \eqref{hauteur arakelov} est en
r{\'e}alit{\'e} une somme finie.  En effet, il est facile de voir que pour
toute place finie $v$, avec un nombre fini d'exceptions d{\'e}pendant de
$F$, on~a $([F] - \lambda_v, [F] - \lambda_v)_v = 0$.

Rappelons que la \emph{hauteur na{\"\i}ve} $\hn$ est d{\'e}finie par 
$$
\hn(F)
\=
|F|^{-1} \, \sum_{\alpha \in F} \sum_{v \in M_K}
\log^+\| x \|_v
=
- \sum_{v \in M_K} \lpar \lambda_v, [F] \rpar_v ~.
$$
Comme $\lpar \lambda_v, \lambda_v \rpar_v =0 $ en toutes les
places, la hauteur na{\"\i}ve $\hn$ co{\"\i}ncide avec la hauteur d{\'e}termin{\'e}e par
la mesure ad{\'e}lique $\{ \lambda_v \}_{v \in M_K}$.
\begin{proposition}\label{hauteur de Weil}
  Pour toute mesure ad{\'e}lique $\rho$, la fonction $h_\rho$ est une
  hauteur de Weil.  En d'autres termes, la diff{\'e}rence $h_\rho -\hn$
  est uniform{\'e}ment born{\'e}e sur $\ov{K}$.
\end{proposition}
Nous verrons plus tard que le minimum essentiel de $h_\rho$ est
non-n{\'e}gatif. On appelle $h_\rho$ \emph{la hauteur d{\'e}finie par} $\rho$.
Les hauteurs que nous avons ainsi construit correspondent pr{\'e}cis{\'e}ment
aux hauteurs \emph{ad{\'e}liques int{\'e}grables} consid{\'e}r{\'e}es en~\cite[\S
2.2]{ACL2} et intervenant dans l'{\'e}nonc{\'e} du Th{\'e}or{\`e}me~4.2, op. cit.

La d{\'e}monstration de la Proposition~\ref{hauteur de Weil} s'appuie sur
le lemme suivant.
\begin{lemme}\label{definition equivalente}
Pour chaque mesure ad{\'e}lique $\rho = \{ \rho_v \}_{v \in M_K}$ on~a
$$
h_\rho(\infty) = \frac{1}{2} \sum_{v \in M_K} \lpar \rho_v, \rho_v \rpar_v
\ \mbox{ et } \ 
h_\rho(F) = h_\rho(\infty) - \sum_{v \in M_K} \lpar \rho_v, [F] \rpar_v.
$$
De plus, ces propri{\'e}t{\'e}s d{\'e}terminent la hauteur $h_\rho$.
\end{lemme}
\begin{proof}
  La premi{\`e}re {\'e}galit{\'e} r{\'e}sulte du fait que pour chaque place $v$ on~a
  $([\infty] - \rho_v, [\infty] - \rho_v)_v = (\rho_v, \rho_v)_v$.
  Pour montrer la deuxi{\`e}me {\'e}galit{\'e} notons qu'on~a
\begin{equation}\label{E-1000}
([F] - \rho_v, [F] - \rho_v)_v = 
(\rho_v, \rho_v)_v - 2(\rho_v, [F])_v + ([F], [F])_v ~,
\end{equation}
et que $([F], [F])_v = - |F|^{-2} \log |\Delta_F|_v$, o{\`u} $\Delta_F =
\prod_{\a, \a' \in F, \a \neq \a'} (\a - \a')$.  Comme $F$ est
invariant sous l'action de $\Gal(\ov{K}/K)$ on~a $\Delta_F \in K$ et
alors la formule du produit implique qu'on a $\sum_{v \in M_K} \lpar
[F], [F] \rpar_v = 0$.  On obtient alors l'{\'e}galit{\'e} desir{\'e}e en
sommant~\eqref{E-1000} sur toutes les places.

Pour montrer la derni{\`e}re assertion du lemme, notons simplement que la
deu\-xi{\`e}me formule d{\'e}termine $h_\rho$ {\`a} une constante additive pr{\`e}s et
que la premi{\`e}re formule d{\'e}termine cette constante.
\end{proof}

\begin{proof}[D{\'e}monstration de la Proposition~\ref{hauteur de Weil}.]
  Il suffit de montrer que la difference $h_\rho - \hn$ est born{\'e}e.
  Soit $N$ le sous-ensemble fini de $M_K$ des places en lequelles
  $\rho_v \neq \lambda_v$.  Pour chaque $v \in N$ soit $g_v$ le
  potentiel de la mesure $\rho_v - \lambda_v$ normalis{\'e} de tel fa{\c c}on
  que $g_v(\infty) = 0$.  Notons que par hypoth{\`e}se chacune des
  fonctions $g_v$ est continue et donc born{\'e}e.
Par le Lemme~\ref{definition equivalente}, on~a
\begin{equation*}
  h_\rho(F) - h_\rho(\infty) - \hn(F) \ = \ - \sum_{v \in N} \lpar
  \rho_v - \lambda_v, [F] \rpar_v  = |F|^{-1} \sum_{v \in N}
  \sum_{\a \in F} g_v(\alpha)~,
\end{equation*}
d'o{\`u} $|h_\rho(F) - h_\rho(\infty) - \hn(F)| \le \sum_{v \in N} \sup \{
|g_v| \}.$
\end{proof}
%%%%%%%%%%%%%%%%%%%%%%%%%%%%%%%%%%%%%%%%%%%%%%%%%%%%%%%%%%%%%%%%%%%%%%%%%%%%%%%%%%%%
%
%
%     Tentative d explication du lien avec les hauteurs d Arakelov.
%
%
%
%
%Le fibr{\'e} $L_v\=O(1)$ sur $\pv$ est donn{\'e} dans les deux ouverts de
%trivialisation $[z_0:1]$ et $[1:z_1]$ par l'application de recollement
%$(z_1 , v_1 ) = (1/ z_0 , v_0 / z_0)$. On peut donc le munir de la
%m{\'e}trique continue $|\cdot|_\star$ d{\'e}finie dans la premi{\`e}re carte
%par $|v_0|_\star \= |v_0| \min \{ 1 , 1/|z_0| \}$.  Aux places
%infinies, la courbure de $|\cdot |_\star$ est la mesure de Lebesgue
%sur le cercle unit{\'e}, \ie $\lambda_v$.
%
%Si $\rho$ est une famille admissible de m{\'e}triques, on peut trouver
%des fonctions continues $g_v$ sur $\berv$ telles que $\rho_v -
%\lambda_v = \Delta g_v$ pour toute place. On peut alors d{\'e}finir une
%m{\'e}trique $|\cdot|_{\rho,v}$ sur chaque fibr{\'e} ample $L_v$ sur
%$\pv$ en posant $|\cdot|_{\rho,v} = |\cdot|_\star \times \exp ( -
%g_v)$. 
%
%
%
%Dans le cas o{\`u} $\rho_v=\lambda_v$ en toute place fini, la hauteur
%associ{\'e} $h_\rho$ est {\'e}gale {\`a} la hauteur d'Arakelov associ{\'e} {\`a}
%la famille de m{\'e}triques $|\cdot|_v$. 
%
%1. metrique formelle
%
%2. metrique algebrique
%
%3. metrique semi-positive
%
%4. metrique integrable
%
%
%
%
%%%%%%%%%%%%%%%%%%%%%%%%%%%%%%%%%%%%%%%%%%%%%%%%%%%%%%%%%%%%%%%%%%%%%%%%%%%%%%%%%%%%%

%%%%%%%%%%%%%%%%%%%%%%%%%%%%%%%%%%%%%%%%%%%%%%%%%%%%%%%%%%%%

\subsection{Formule de Mahler.}\label{S-formule de Mahler}
Ce paragraphe est d{\'e}di{\'e} {\`a} montrer la Proposition~\ref{P-standard}.

Soit $\rho = \{ \rho_v \}_{v \in M_K}$ une mesure ad{\'e}lique et fixons
$\alpha \in \P^1(\ov{K})$.  On d{\'e}signe par $F$ l'orbite de $\alpha$
sous l'action du groupe $\Gal(\ov{K} / K)$.  La preuve consiste {\`a}
r{\'e}{\'e}crire, pour chaque $v \in M_K$, l'{\'e}nergie $([F], \rho_v)_v$ en
ter\-me d'une int{\'e}grale sur $\berv$. On note $P$ le polynome minimal
et unitaire de $\a$ sur $K$.

Si $v$ est une place infinie, on~a $\prod_{\b \in F} (z - \b) =
P(z)$, et ~\eqref{e-cal} donne
\begin{eqnarray*}
\deg(\alpha) \, ( [F] , \rho_v)_v
& =& 
\sum_{\b\in F} ( \b, \rho_v)_v =
- \sum_{\b \in F} \int_{\berv}\log|z-\b|_v \, d\rho_v(z)
\\
&=& - \int_{\berv} \log|P(z)|_v\, d\rho_v(z)~.
\end{eqnarray*}
Supposons d{\'e}montr{\'e}e pour toute place $v$ finie, et tout
$\cS \in \berv$, l'{\'e}galit{\'e}
\begin{equation}\label{E-918}
 \sum_{\b \in F} \log \sup \{ \cS
, \b \} = \log |P(\cS) |_v~.
\end{equation}
De~\eqref{e-500}, on tire de m{\^e}me $\deg(\alpha) \, ( [F] , \rho_v)_v =
- \int_{\berv} \log|P(z)|_v \, d\rho_v$.  En sommant ceci sur toutes
les places (finies et infinies), en remarquant que $\rho_v$ ne charge
pas l'infini, et en appliquant le Lemme~\ref{definition equivalente},
on trouve le r{\'e}sultat escompt{\'e}.

Il nous reste donc {\`a} d{\'e}montrer~\eqref{E-918}. On remarque tout
d'abord que les deux membres ont m{\^e}me laplacien.  Celui-ci est en
effet {\'e}gal {\`a} $\sum_{\b\in F}( [\b] - [\infty])$
par~\eqref{e-330}. Deux fonctions sur $\berv$ dont les laplaciens
co{\"\i}ncident diff{\`e}rent d'une constante.  Mais il est clair que
les deux fonctions sont identiques sur $\pv$, ce qui
d{\'e}mon\-tre~\eqref{E-918}.

%%%%%%%%%%%%%%%%%%%%%%%%%%%%%%%%%%%%%%%%%%%%%%%%%%%%%%%%%%%%

\subsection{Equidistribution des points de petite hauteur.}
\label{S-equidistribution}
\begin{theoreme}\label{T-prec}
  Soit $\rho$ une mesure ad{\'e}lique et $\{ F_n \}_{n \ge 0}$ une suite
  d'ensembles finis distincts deux {\`a} deux et
  $\Gal(\ov{K}/K)$-invariants, et telle que
$$
\ov{\lim}_{n \to \infty} h_\rho(F_n) \le 0 ~.
$$
Alors on~a $\lim_{n \to \infty} h_\rho(F_n) = 0$ et pour toute
place $v$ de $M_K$ on~a $[F_n] \to \rho_v$ lorsque $n \to \infty$.
\end{theoreme}
Cet {\'e}nonc{\'e}, combin{\'e} {\`a} la Proposition~\ref{hauteur de Weil}, implique
imm{\'e}diatement les deux Th{\'e}or{\`e}mes~\ref{T-main1} et~\ref{T-main2} cit{\'e}s
dans l'introduction.  La d{\'e}monstration du Th{\'e}or{\`e}me~\ref{T-prec}
s'appuie sur le lemme suivant.
\begin{lemme}\label{bonne reduction}
  Pour toute place finie $v$ de $K$ et tout sous-ensemble fini $F$ de
  $\C_v$ on~a
$$
(\lambda_v - [F], \lambda_v - [F])_v \ge 0~.
$$
\end{lemme}
\begin{proof}
L'in{\'e}galit{\'e} ultram{\'e}trique $|z - z'|_v \le \max \{ |z|_v, |z'|_v
\}$ implique qu'on a $\log|z - z'|_v \le \log^+|z|_v +
\log^+|z'|_v$. On en d{\'e}duit
\begin{eqnarray*}
 (\lambda_v - [F], \lambda_v - [F])_v 
&=&
 \frac{2}{|F|} \sum_{z \in F}
\log^+|z|_v - \frac1{|F|^{2}} \sum_{z, z' \in F, z \neq z'} \log |z - z'|_v
\\
&\ge&
\frac2{|F|^2} \sum_{z \in F}
\log^+|z|_v \ge
0~,
\end{eqnarray*}
ce qui termine la preuve.
\end{proof}
\begin{proof}[D{\'e}monstration du Th{\'e}or{\`e}me~\ref{T-prec}]
  Soit $\{ F_n \}_{n \ge 0}$ une suite d'ensembles finis de la droite
  projective $\P^1(\ov{K})$, invariants sous l'action du groupe de
  Galois $\Gal(\ov{K} / K)$ et tels que $\ov{\lim}_{n \to \infty}
  h_\rho(F_n) \le 0$.  On suppose de plus que ces ensembles sont
  distincts deux {\`a} deux.
  
  Montrons tout d'abord que $|F_n| \to \infty$ lorsque $n \to \infty$.
  Supposons que ce ne soit pas le cas.  Quitte {\`a} prendre une
  sous-suite on suppose $M \= \sup |F_n| < \infty$.  Notons de m{\^e}me
  $H\= \sup \hn(F_n)$. C'est une quantit{\'e} finie, born{\'e}e par $\sup
  h_\rho(F_n) + \sup_{\ov{K}} | h_\rho -\hn| < + \infty$ car $h_\rho$
  est une hauteur de Weil. Tout {\'e}l{\'e}ment de $F_n$ est alors de degr{\'e} au
  plus $M$ et de hauteur (na{\"\i}ve) au plus $H$.  Mais la propri{\'e}t{\'e}
  de Northcott affirme que l'ensemble de ces points est fini, voir par
  exemple~\cite[Theorem~B.2.3]{hindry-silverman}. Ce qui contredit
  notre hypoth{\`e}se que les $F_n$ sont distincts deux {\`a} deux.
  
  Soit $B$ le sous ensemble de $M_K$ des places finies $v$ telles que
  $\rho_v = \lambda_v$.  Le compl{\'e}mentaire $N = M_K \setminus B$ est
  un ensemble fini qui contient toutes les places infinies.  Lorsque
  $v \in B$, on~a $\rho_v = \lambda_v$ et le Lemme~\ref{bonne
    reduction} implique que pour tout $n \ge 1$ on~a
$$
(\rho_v - [F_n], \rho_v - [F_n])_v \ge 0~.
$$
On a donc
$$
h_\rho(F_n)
\ge
\frac{1}{2} \sum_{v \in N} \lpar \rho_v - [F_n], \rho_v - [F_n] \rpar_v ~.
$$
Par ailleurs, pour chaque $v \in N$ la Proposition~\ref{P-130}
lorsque $v$ infinie et la Proposition~\ref{P-230} lorsque $v$ finie,
fournissent l'in{\'e}galit{\'e}
$$
\underline{\lim}_{n \to \infty}
(\rho_v - [F_n], \rho_v - [F_n])_v \ge 0~,
$$
d'o{\`u} on d{\'e}duit que $\lim_{n \to \infty} h_\rho (F_n) = 0$.

D'autre part, fixons $v_0 \in M_K$.  On peut maintenant appliquer la
Proposition~\ref{P-casdeg2} (lorsque $v_0$ est infinie) et la
Proposition~\ref{critere de convergence} (lorsque $v_0$ finie): pour
montrer que $[F_n] \to \rho_{v_0}$, il nous suffit de montrer que
\begin{equation}\label{e-304}
\ov{\lim}_{n \to \infty} \,  
(\rho_{v_0} - [F_n], \rho_{v_0} - [F_n])_{v_0} \le0 ~.
\end{equation}
On a,
\begin{eqnarray*}
\lefteqn{\ov{\lim}_{n \to \infty} \, 
\lpar \rho_{v_0} - [F_n], \rho_{v_0} - [F_n] \rpar_{v_0} \le }& &\\
& \le&  
\ov{\lim}_{n \to \infty} 
\lpar \rho_{v_0} - [F_n], \rho_{v_0} - [F_n]\rpar_{v_0}
+\\
&& \hspace*{1.5cm}
+ \sum_{v \in N \setminus \{ v_0 \}} \underline{\lim}_{n \to \infty}
 \lpar \rho_v - [F_n], \rho_v - [F_n] \rpar_v 
\\
& \le &
\ov{\lim}_{n \to \infty} \sum_{v \in N \cup \{ v_0 \}} 
\lpar \rho_v - [F_n], \rho_v - [F_n] \rpar_v \\
%&\le &
%\ov{\lim}_{n \to \infty} \sum_{v \in M_K}  
%\lpar \rho_v - [F_n], \rho_v - [F_n] \rpar_v \\
& \le  &  
\ov{\lim}_{n \to \infty} \, 2 h_\rho(F_n) \le 0.
\end{eqnarray*}
Ceci d{\'e}montre~\eqref{e-304} et termine la preuve du th{\'e}or{\`e}me.
\end{proof}

%%%%%%%%%%%%%%%%%%%%%%%%%%%%%%%%%%%%%%%%%%%%%%%%%%%%%%%%%%%%

\subsection{Version quantitative de l'{\'e}quidistribution.}\label{S-vitesse}
Nous nous pla{\c c}ons dans le m{\^e}\-me cadre que la section pr{\'e}c{\'e}dente. Afin
d'{\'e}noncer notre r{\'e}sultat de mani{\`e}re pr{\'e}ci\-se, nous aurons besoin
d'introduire un peu de terminologie.
\begin{definition}
  Soit $v$ une place finie. Une fonction continue $\varphi : \berv \to
  \R$ est dite \emph{de classe $\cC^k$} pour $k\ge 1$ si et seulement si elle est
  localement constante hors d'un sous-arbre \emph{fini et ferm{\'e}}
  $\cT\subset \H_v$, et si de plus $\cT$ est la r{\'e}union d'un nombre
  fini de segments ferm{\'e}s sur lesquels $\varphi$ est de classe
  $\cC^k$.
\end{definition}
Fixons $\cS_0\in \H_v$ un point base, et $\varphi$ une fonction de
classe $\cC^1$. On peut alors d{\'e}finir $\partial \varphi( \cS)$ comme
la d{\'e}riv{\'e}e {\`a} gauche en $\cS$ de l'application $\varphi$ restreint au
segment $[\cS_0,\cS]$ (param{\'e}tr{\'e} par la distance {\`a} $\cS_0$). Hors d'un
arbre fini de mesure $\lambda$-fini, $\partial \varphi$ est
identiquement nulle et on peut donc poser
$$
\langle \varphi, \varphi \rangle \= \int_\berv (\partial
\varphi)^2\, d\lambda~.
$$
Cette quantit{\'e} est ind{\'e}pendante du choix du point base $\cS_0$, car
$\varphi$ est d{\'e}termin{\'e}e au signe pr{\`e}s hors de $\H_v^\Q$ qui est
de mesure $\lambda$ nulle.

On montre que toute fonction $\varphi$ de classe $\cC^2$ d{\'e}finit un
potentiel au sens de \S\ref{S-poten}, et que la mesure $(\partial^2
\varphi)\times \lambda$ est une mesure sign{\'e}e, {\'e}gale {\`a} $\Delta
\varphi$. On a alors
$$
\langle \varphi , \varphi \rangle  = (\Delta \varphi , \Delta \varphi)~.
$$
Plus g{\'e}n{\'e}ralement, si $\rho$ est une mesure sign{\'e}e {\`a}
potentiel continu et telle que $\rho(\berv) = 0$, on~a $(\Delta
\varphi, \rho) = \int_\berv (\partial \varphi) \, f_\rho \cdot 
d\lambda$, o{\`u} $f_\rho(\cS) = \rho \{ \cS', \, [\cS_0, \cS] \cap [ \cS
, \cS'] = \{ \cS\}\}$. On en d{\'e}duit donc l'in{\'e}galit{\'e} de Cauchy-Schwartz
\begin{equation}\label{E-477}
\left| \int_\berv  \varphi\, d\rho \right| 
\le
 \langle \varphi , \varphi \rangle \,
 (\rho , \rho) ~. 
\end{equation}
Pour m{\'e}moire, l'in{\'e}galit{\'e} analogue lorsque la place $v$ est infinie
s'{\'e}crit
\begin{equation}\label{E-475}
\left| \int_\berv  \varphi \, d\rho \right| 
= \left| \int_\berv  d \varphi \wedge d^c g \right| 
\le \langle \varphi , \varphi \rangle \,
 (\rho , \rho) ~, 
\end{equation}
avec $\rho = dd^c g$ et $\langle \varphi , \varphi \rangle = \int
d\varphi \wedge d^c \varphi$. 

\medskip

Le contr{\^o}le de la vitesse de convergence n{\'e}cessite une borne sur
la constante de Lipschitz de la fonction consid{\'e}r{\'e}e. On fixe donc
sur $\mathbb{P}^1(\C_v)$ la m{\'e}trique sph{\'e}rique que l'on note $\mathsf{d}$,
et on pose
$$
\lip(\varphi) \= 
\sup\left\{ \frac{|\varphi(x) - \varphi(y)|}{\mathsf{d}(x,y)} , \, 
x \neq y \in \mathbb{P}^1(\C_v) \right\}~.
$$
Aux places infinies, il est clair que cette constante est finie si
$\varphi$ est de classe $\cC^1$, {\'e}gale au supremum de la norme de la
diff{\'e}rentielle de $\varphi$ (pour la m{\'e}trique $\mathsf{d}$) par le
th{\'e}or{\`e}me des accroissements finis. Aux places finies,
$\lip(\varphi)$ est aussi fini si $\varphi$ est de classe $\cC^1$, car
il existe $\e>0$ pour lequel $\mathsf{d}(x,y) < \e$ implique $\varphi(x) =
\varphi(y)$.  Notons cependant que $\lip(\varphi)$ n'est dans ce cas
pas born{\'e} par le supremum de $\partial \varphi$ sur $\berv$.

Nous pouvons maintenant {\'e}noncer le th{\'e}or{\`e}me principal de cette
section.  On dit qu'une mesure ad{\'e}lique $\rho = \{ \rho_v \}_{v \in
  M_K}$ est \emph{{\`a} potentiel H{\"o}lder} si pour toute place $v \in M_K$,
la mesure $\rho_v - \lambda_v$ a un potentiel H{\"o}lder. En d'autres
termes, si on peut {\'e}crire $\rho_v - \lambda_v = \Delta g_v$, avec
$g_v$ de telle sorte que $|g_v(z) - g_v(w)| \le C \, \mathsf{d}(z,w)^\kappa$
pour tout $z,w\in \pp$ et pour des constantes $C>0$, et $0<\kappa \le
1$. On dira alors que $g_v$ est $\kappa$-H{\"o}lder.
\begin{theoreme}\label{T-estimprecise}
  Soit $\rho = \{ \rho_v \}_{v \in M_K}$ une mesure ad{\'e}lique {\`a}
  potentiel $\kappa$-H{\"o}lder pour $\kappa \le 1$.  Alors il existe une
  constante $C > 0$ ne d{\'e}pendant que de $\rho$, telle que pour toute
  place $v$, pour toute fonction $\varphi$ de classe $\cC^1$ sur
  $\berv$, et pour tout ensemble fini $\Gal(\ov{K}/K)$-invariant
  $F\subset \P^1(\ov{K})$ on ait
\begin{equation}
\left| \frac1{|F|} 
\sum_{\a \in F} \varphi(\alpha) -  \int_{\berv} \varphi\, d\rho \right|
\le 
\frac{\lip (\varphi)}{|F|^{1/\kappa}}+
\left( h_\rho(F) + C\, \frac{\log |F|}{|F|}
\right) \,  \langle \varphi , \varphi \rangle~,
\end{equation}
o{\`u} $|F|$ d{\'e}signe la cardinalit{\'e} de $F$.
\end{theoreme}
\begin{remarque}
  Aux places infinies, on~a vu que le supremum de la diff{\'e}rentielle
  $d\varphi$ est {\'e}gal {\`a} $\lip(\varphi)$, ce qui implique $\langle
  \varphi , \varphi \rangle \le \lip(\varphi)$. En ajoutant {\`a} cel{\`a} le
  fait que $x^{-1/\kappa}$ est n{\'e}gligeable devant $\log x/x$ pour $x$
  grand, on en d{\'e}duit le Th{\'e}or{\`e}me~\ref{T-estim} cit{\'e} dans
  l'introduction.
\end{remarque}
\begin{proof}
  Soit $N$ la famille finie des places infinies et des places
  finies~$v$ en lesquelles $\rho_v\neq \lambda_v$ et posons $M = (|N|
  + 1)\max \{ N_v ; v \in M_K \}$.  Notons $\heta$ un module de
  continuit{\'e} commun {\`a} tous les potentiels de $\rho_v - \lambda_v$ pour
  $v \in N$ et soit $\eta(\e) \= \heta(\e) + \e$.
  
  Soit $\e > 0$ et soit $F$ un sous ensemble fini de $\P^1(\ov{K})$
  invariant sous l'action du groupe de Galois $\Gal(\ov{K}/K)$.
 % Notons que la Proposition~\ref{P-230} implique que pour toute place
 % finie $v$ on~a
%$$
%([F]_\e - \lambda_v, [F]_\e - \lambda_v)_v \le |F|^{-1} \log \e^{-1}.
%$$
Alors, pour chaque place $v_0 \in M_K$ l'application des
propositions~\ref{P-130} et~\ref{P-230} donne
\begin{multline*}
\sum_{v \in N \cup \{ v_0 \}} \lpar [F]_\e - \rho_v , [F]_\e - \rho_v \rpar_v
\le \\ \le
\sum_{v \in N \cup \{ v_0 \}} \lpar [F] - \rho_v , [F] - \rho_v \rpar_v
+ M(2 \eta(\e) + |F|^{-1} \, ( C + \log \e^{-1}) )~.
\end{multline*}
Pour chaque $v \in M_K \setminus N$ on~a $\rho_v = \lambda_v$ et le
Lemme~\ref{bonne reduction} implique qu'on a $([F] - \rho_v, [F] -
\rho_v)_v \ge 0$, dont on tire
$$
\lpar [F]_\e - \rho_{v_0} , [F]_\e - \rho_{v_0} \rpar_{v_0}
\le
h_\rho(F) + M \left( 2 \eta(\e) + |F|^{-1} ( C + \log \e^{-1}) \right)~.
$$
Par hypoth{\`e}se on peut supposer $\eta (r) \lesssim r^\kappa$, et en
prenant $\e^{-\kappa} = |F|$ on obtient finalement,
$$
\lpar [F]_\e - \rho_{v_0} , [F]_\e - \rho_{v_0} \rpar_{v_0}
\le h_\rho (F) + C ' \frac{\log |F| } {|F|} ~,
$$
pour une constante $C'>0$ ne d{\'e}pendant que de $\rho$.

On {\'e}crit alors $ [F]_\e - \rho_v = \Delta g_\e$ avec $g_\e$ continue.
Pour une fonction $\varphi$ de classe $\cC^1$, les in{\'e}galit{\'e}s de
Cauchy-Schwartz~\eqref{E-475} et~\eqref{E-477} impliquent
\begin{multline*}
\left| \int \varphi\, d([F]_\e - \rho_v) \right|
= 
\left| \int \varphi\, \Delta g_\e \right|
\le
\langle \varphi , \varphi \rangle
\, ([F]_\e - \rho_v, [F]_\e - \rho_v)
\le\\\le \langle \varphi , \varphi \rangle 
\left( h_\rho(F) + C' \, \frac{\log |F| } {|F|}
\right)~.
\end{multline*}
Il nous reste {\`a} estimer $\int\varphi\, d([F]_\e - [F])$.  Pour cel{\`a},
on remarque que la m{\'e}trique induite par $|\cdot|_v$ est plus grande ou
{\'e}gale que la m{\'e}trique sph{\'e}rique utilis{\'e}e pour d{\'e}finir $\lip(\varphi)$.
Pour toute boule $B(x,\e)$ pour la m{\'e}trique $|\cdot|_v$, on~a donc
l'estim{\'e}e $\sup_{B(x,\e)} |\varphi - \varphi(x)| \le \lip(\varphi)
\times \e$, ce qui implique
$$
\left| \int\varphi\, d([F]_\e -[F]) \right| \le
\lip(\varphi) \times \e = \lip(\varphi)|F|^{-1/\kappa}.
$$

Ceci conclut la preuve du Th{\'e}or{\`e}me~\ref{T-estimprecise}.
\end{proof}

%
%%%%%%%%%%%%%%%%%%%%%%%%%%%%%%%%%%%%%%%%%%%%%%%%%%%%%%%%%%%%%%%%%%
%

\section{Dynamique des fractions rationnelles.}\label{S-dyn}
Fixons un corps de nombres $K$ et une fraction rationnelle $R$ {\`a}
coefficients dans $K$ et de degr{\'e} $D \ge 2$.  Pour tout sous-ensemble
fini $F$ de $\P^1(\ov{K})$ invariant sous l'action du groupe de Galois
$\Gal(\ov{K} / K)$, l'ensemble $R(F)$ est aussi invariant sous
l'action du groupe de Galois.  On d{\'e}finit alors la {\it hauteur
  normalis{\'e}e} de $F$ associ{\'e}e {\`a} $R$ par
$$
h_R(F) \= \lim_{n \to \infty} D^{-n} \hn ( R^n(F) )~.
$$
On peut montrer que la limite ci-dessus existe directement (voir
par exemple~\cite[Theorem~B.4.1]{hindry-silverman}).  
%La preuve du Lemme~\ref{L-Mahler dynamique} que nous donnons plus bas
%contient implicitement une preuve de l'existence de cette limite.
Cette hauteur est caracteris{\'e}e comme l'unique hauteur de Weil, {\`a} une
constante multiplicative pr{\`e}s, satisfaisant la formule d'invariance
$h_R \circ R = D \cdot h_R$.

A chaque place $v$ de $K$ on associe {\`a} $R$ une mesure de probabilit{\'e}
$\rho_{R, v}$, appel{\'e}e \emph{mesure d'{\'e}quilibre}.  Cette mesure d{\'e}crit
la distribution des pr{\'e}images it{\'e}r{\'e}es de chaque point non exceptionnel
dans $\berv$, voir \S~\ref{S-mesure-inv} pour plus de pr{\'e}cisions.

Le th{\'e}or{\`e}me suivant est le r{\'e}sultat principal de cette section.
\begin{theoreme}\label{normalisee est arakelov}
  Soit $R$ une fraction rationnelle de degr{\'e} au moins~$2$ et {\`a}
  coefficients dans $K$.  Pour chaque place $v$ de $K$ soit $\rho_{R,
    v}$ la mesure d'{\'e}quilibre de $R$ correspondante.  Alors $\rho_R =
  \{ \rho_{R, v} \}_{v \in M_K}$ est une mesure ad{\'e}lique et la hauteur
  normalis{\'e}e $h_R$ co{\"\i}ncide avec la hauteur $h_{\rho_R}$ d{\'e}finie par
  la mesure ad{\'e}lique $\rho_R$.
\end{theoreme}

\begin{remarque}
  Comme nous l'avons mentionn{\'e} dans l'introduction, la mesure ad{\'e}lique
  $\rho_R$ est {\`a} potentiel H{\"o}lder et par cons{\'e}quent la hauteur
  normalis{\'e}e $h_R$ v{\'e}rifie toutes les conditions du
  Th{\'e}or{\`e}me~\ref{T-estimprecise}. Nous d{\'e}montrerons ce fait
  en~\S~\ref{S-holder}.
\end{remarque}

Le plan de cette section est le suivant.  Apr{\`e}s la construction de la
mesure d'{\'e}qui\-libre d'une fraction rationnelle {\`a} coefficients dans
$\C_v$ au~\S~\ref{S-mesure-inv}, nous parlerons succintement de la
notion de \emph{bonne r{\'e}duction} au~\S~\ref{S-bonne reduction}.  Nous
montrons au~\S~\ref{S-birapport} une formule de transformation du
birapport sous l'action d'une fraction rationnelle.  Le birapport est
{\'e}troitement li{\'e} {\`a} l'{\'e}nergie d{\'e}finie au~\S\ref{S-energie}
et~\S\ref{capacite p-adique}.  On donne alors la d{\'e}monstration du
Th{\'e}or{\`e}me~\ref{normalisee est arakelov} au~\S~\ref{s-normalisee est
  arakelov} et celle du Th{\'e}or{\`e}me~\ref{thm:mandel}
au~\S~\ref{S-mandel}. On montre ensuite en~\S~\ref{S-holder} la
continuit{\'e} H{\"o}lder des potentiels des mesure d'{\'e}quilibre.  Enfin, on
d{\'e}crit un exemple au~\S~\ref{S-exemple}, montrant que le proc{\'e}d{\'e} de
r{\'e}gularisation dans la preuve de l'{\'e}quidistribution des petits points
est en effet n{\'e}cessaire.

%%%%%%%%%%%%%%%%%%%%%%%%%%%%%%%%%%%%%%%%%%%%%%%%%%%%%%%%
\subsection{Mesure d'{\'e}quilibre.}\label{S-mesure-inv}
Fixons une place $v$ de $M_K$ et soit $R$ une fraction rationnelle {\`a}
coefficients dans $\C_v$ et de degr{\'e} $D \ge 2$.  Lorsque $v$ est
infinie il est clair que la fraction rationnelle $R$ agit contin{\^u}ment
sur $\P^1(\C_v) = \berv$.  Lorsque $v$ est finie, on~a vu
au~\S\ref{espace des semi-normes} que $R$ agissait contin{\^u}ment sur
$\berv$.

Dualement, on peut d{\'e}finir de mani{\`e}re unique une application
lin{\'e}aire continue $R^*$ agissant sur $\cM_v$, et v{\'e}rifiant les
propri{\'e}t{\'e}s suivantes:
\begin{itemize}
\item[$\bullet$]
$R^*[z] =
\sum_{R(w) = z} \deg_R(w)\, [w]$ \ pour tout $z\in \berv$; 
\item[$\bullet$]
$R^* \Delta g = \Delta (g \circ R)$ \ pour tout potentiel~$g$.
\end{itemize}
Lorsque $w\in\pv$, l'entier $\deg_R(w)$ est la multiplicit{\'e} locale
de $R$ en $w$. Lorsque $v$ est finie et $w\in\H_v$, nous renvoyons
{\`a}~\cite{R} pour une d{\'e}finition pr{\'e}cise de ce nombre, ainsi que
son interpr{\'e}tation g{\'e}om{\'e}trique.  Notons que la masse
de $R^*\rho$ est toujours {\'e}gale {\`a} $D$ fois la masse de $\rho$.

On peut trouver un potentiel \emph{continu} $g$ tel que
$D^{-1} R^*\lambda_v - \lambda_v = \Delta g$.
Par it{\'e}ration successive, on en d{\'e}duit pour tout $n\ge0$, que $D^{-n}
R^{n*}\lambda_v - \lambda_v = \Delta \left( \sum_{k \ge 0} D^{-k} g
\circ R^k\right)$. Cette derni{\`e}re s{\'e}rie d{\'e}finit une suite de
fonctions continues qui converge uniform{\'e}ment sur $\berv$.
Il s'ensuit que la suite de mesures de probabilit{\'e} $\{
D^{-n} (R^n)^*\lambda_v \}_{n \ge 1}$ converge vers une mesure que
l'on d{\'e}si\-gnera par $\rho_{R,v}$.  Cette mesure est invariante~: on~a
$R^* \rho_{R,v} = D \cdot \rho_{R,v}$ et $R_* \rho_{R,v} =
\rho_{R,v}$.  Notons que $\rho_{R,v}$ admet un potentiel d{\'e}fini et
continu sur $\berv$. En particulier pour tout $\alpha \in \C_v$,
l'accouplement $([\alpha], \rho_{R, v})$ est bien d{\'e}fini et fini
(voir~Propositions~\ref{L-integr} et~\ref{L-integr-Cp}).

La construction de la mesure $\rho_{R,v}$ ci-dessus est d{\^u}e
{\`a}~\cite{Brolin},~\cite{FLM}, et~\cite{Lyubich} dans le cas complexe,
et {\`a}~\cite{FR} dans le cas ultram{\'e}trique. Nous renvoyons {\`a} ces
articles pour plus de d{\'e}tails sur les propri{\'e}t{\'e}s dynamiques de
ces mesures.

%%%%%%%%%%%%%%%%%%%%%%%%%%%%%%%%%%%%%%%%%%%%%%%%%%%%%%%%

\subsection{Bonne r{\'e}duction.}\label{S-bonne reduction}
Soit $v$ une place finie de $K$.  On d{\'e}signe par $\cO_v = \{ z \in
\C_v ; \, |z|_v \le 1 \}$ l'anneau des entiers de $\C_v$ et par
$\widetilde{\C}_v$ le corps r{\'e}siduel de $\C_v$.  Pour tout
polyn{\^o}me $P \in \cO_v[z_0, z_1]$, on note $\widetilde{P}$ sa
projection dans $\widetilde{\C}_v$.

Consid{\'e}rons une fraction rationnelle $R$ {\`a} coefficients dans $\C_v$,
et {\'e}crivons $R$ en coordonn{\'e}es homog{\`e}nes $[P_0: P_1]$.  Quitte {\`a}
rempla{\c c}er $P_0$ et $P_1$ par $\lambda P_0$ et $\lambda P_1$
respectivement, on peut supposer que les coefficients de $P_0$ et de
$P_1$ sont de norme au plus {\'e}gale {\`a}~1 et qu'au moins l'un des
coefficients de $P_0$ ou de $P_1$ est de norme {\'e}gale {\`a}~1.  On dit
alors que $R$ a {\it bonne r{\'e}duction} lorsque les deux 
polyn{\^o}mes de degr{\'e} $\deg(R)$ et homog{\`e}nes $\widetilde{P}_0$,
$\widetilde{P}_1$ n'ont pas de racine commune autre que~$0$. Ceci
revient {\`a} dire que la fraction rationnelle
$[\widetilde{P}_0:\widetilde{P}_1]$ induite par $R$ sur
$\mathbb{P}^1(\widetilde{\C}_v)$ a le m{\^e}me degr{\'e} que $R$.

De fa{\c c}on {\'e}quivalente, $R$ a bonne r{\'e}duction si et seulement si
le point canonique $\cS_\can$ est compl{\`e}tement invariant par $R$~:
$R^{-1}\{\cS_\can\} = \{\cS_\can\}$.  Dans ce cas, il est facile de
voir que l'on a $\rho_R = \lambda_v$.

Lorsque $R$ est {\`a} coefficients dans une extension finie $K'$ de $\Q$,
l'ensemble des $v \in M_{K'}$ finies o{\`u} $R$ n'a pas bonne r{\'e}duction
est fini. 

%%%%%%%%%%%%%%%%%%%%%%%%%%%%%%%%%%%%%%%%%%%%%%%%%%%%%%%%

\subsection{Birapport.}\label{S-birapport}
Fixons $v \in M_K$.  Etant donn{\'e}s $z_0, z_1, w_0, w_1 \in \C_v$ deux
{\`a} deux distincts, on pose
$$ (z_0, z_1 , w_0, w_1)_v = \log \frac{|z_0 - w_0|_v \cdot |z_1 -
w_1|_v}{|z_0 - w_1|_v \cdot |z_1 - w_0|_v}~.
$$ 
On voit que cette fonction s'{\'e}tend en une fonction d{\'e}finie
et continue sur l'ensemble
$$
\cD = \{ z_0, z_1, w_0, w_1 \in 
\berv \ ; \ \text{si } z_i = w_j, \text{alors } z_i \notin \pv \}~,
$$ en prenant soin de remplacer $|z-w|_v$ par $\sup \{ z, w \}$ lorsque
$v $ est finie.  Notons qu'en particulier on~a
$$ (z_0, \infty , w_0, w_1)_v = (z_0,w_0)_v - (z_0, w_1)_v~.
$$
Lorsque $\mu_0, \mu_1, \nu_0, \nu_1$ sont des mesures dans $\cM_v$
telles que la fonction $(\cdot, \cdot , \cdot, \cdot)_v$ est int{\'e}grable
par rapport {\`a} $\mu_0 \otimes \mu_1 \otimes \nu_0 \otimes \nu_1$ sur
$\cD$, on pose
$$
(\mu_0, \mu_1 , \nu_0, \nu_1)_v 
= 
\int_\cD (z_0, z_1, w_0, w_1)_v \ d\mu_0(z_0)\otimes d\mu_1(z_1)\otimes 
d\nu_0(w_0)\otimes d\nu_1(w_1)~.
$$ Par commodit{\'e}, on notera aussi $\lpar \cdot, \cdot, \cdot, \cdot
\rpar_v = N_v \cdot (\cdot, \cdot, \cdot, \cdot)_v$.

\smallskip
Pour des mesures de \emph{probabilit{\'e}} et sous les conditions
d'int{\'e}grabilit{\'e} ad{\'e}quates pour que tous les termes existent et
soient finis, une simple int{\'e}gration montre
\begin{eqnarray}
(\mu_0, \mu_1, \nu_0, \nu_1)_v 
&=& 
(\mu_0, \nu_0)_v  + (\mu_1, \nu_1)_v - (\mu_0, \nu_1)_v - (\mu_1, \nu_0)_v~;
\label{E-123}
\\
(\mu_0, [\infty] , \nu_0, \nu_1)_v 
&=& 
(\mu_0, \nu_0)_v  - (\mu_0, \nu_1)_v~.
\label{E-124}
\end{eqnarray}
Les {\'e}quations pr{\'e}c{\'e}dentes sont v{\'e}rifi{\'e}es par exemple lorsque les
mesures de probabilit{\'e} $\mu_i, \nu_j$ sont atomiques, {\`a} support
disjoint et ne chargent pas $\{ \infty\}$, ou si l'on remplace
certaines d'entre elles par des mesures (toujours de probabilit{\'e}) {\`a}
potentiel continu.

\medskip

\noindent
{\bf Formule de transformation.}  {\it Soit $R$ une fonction
rationnelle de degr{\'e} $D \ge 1$ {\`a} coefficients dans $K$.  Soient
$\mu_0, \mu_1, \nu_0, \nu_1$ des mesures dans $\cM_v$ telles que
$$
(R_*\mu_0, R_*\mu_1, \nu_0, \nu_1)_v
\ \mbox{ et } \
(\mu_0, \mu_1, R^*\nu_0, R^*\nu_1 )_v
$$
soient bien d{\'e}finis.
Alors on~a
$$
(R_*\mu_0, R_*\mu_1, \nu_0, \nu_1)_v
=
D^{-1} \cdot (\mu_0, \mu_1, R^*\nu_0, R^*\nu_1)_v~.
$$}
\begin{proof}
Notons qu'on a
$$ (R_*\mu_0, R_*\mu_1, \nu_0, \nu_1)_v = \int_\cD (R(z_0),
R(z_1), w_0, w_1)_v \ d \mu_0 \otimes d \mu_1 \otimes d \nu_0 \otimes d \nu_1~.
$$ Soient $z_0, z_1, w_0, w_1 \in \berv$ tels que $(R(z_0), R(z_1),
w_0, w_1)_v \in \cD$.  Si l'on d{\'e}signe par $a_i^1, \ldots, a_i^D$ les
pr{\'e}images de $w_i$ par $R$ dans $\berv$, compt{\'e}es avec
multiplicit{\'e}, alors on~a
\begin{eqnarray}
(R(z_0), R(z_1), w_0, w_1)_v 
& = & 
\sum_{1 \le k \le D} 
(z_0, z_1, a_0^{k}, a_1^{k})_v
\notag
 \\ & = & 
\frac1D \sum_{1 \le k_0,k_1 \le D} 
(z_0, z_1, a_0^{k_0}, a_1^{k_1})_v~.
\label{E-334}
\end{eqnarray}
Ceci est clair lorsque $z_i, R(z_i), w_i$ et $a_i^k$ appartiennent {\`a}
$\C_v$, car dans ce cas
\begin{equation}\label{E-333}
\left| \frac{R(T) - w_0}{ R(T) - w_1} \right|_v
\cdot
\left| \frac{R(T') - w_1}{ R(T') - w_0} \right|_v
=
\prod_k \left| \frac{T - a_0^k}{T- a_1^k}\right|_v 
\cdot
\left| \frac{T' - a_1^k}{T' - a_0^k}\right|_v 
~.
\end{equation}
Par continuit{\'e}, on voit que~\eqref{E-334} reste vraie sans
restriction pour $z_i, R(z_i), w_i, a_i^k \in \pv$ et m{\^e}me dans
$\berv$.  Notons en effet, que la fonction $\sup\{\cdot, \cdot \}$ 
est continue (faiblement) en chacune des variables. 
On d{\'e}duit facilement de ceci 
$$
\int_\cD (R(z_0), R(z_1), w_0, w_1)_v \ 
d \mu_0 \otimes d \mu_1 \otimes d \nu_0\otimes d \nu_1
=
\frac1D  (\mu_0, \mu_1, R^*(\nu_0), R^*(\nu_1))_v~,
$$
ce qui termine la preuve. 
\end{proof}

%%%%%%%%%%%%%%%%%%%%%%%%%%%%%%%%%%%%%%%%%%%%%%%%%%%%%%%%

\subsection{D{\'e}monstration du Th{\'e}or{\`e}me~\ref{normalisee est arakelov}.}
\label{s-normalisee est arakelov}
Soit $R$ une fraction rationnelle de degr{\'e} $D \ge 2$ {\`a} coefficients
dans $K$.  Pour chaque place $v$ de $K$ la mesure d'{\'e}quilibre
$\rho_{R, v}$ est {\`a} potentiel continu.  De plus, pour toute place
finie $v$ en laquelle $R$ a bonne r{\'e}duction on~a $\rho_v = \lambda_v$.
Comme cette derni{\`e}re propri{\'e}t{\'e} est valable pour toute place, avec un
nombre fini d'exceptions, on conclut que $\rho_R = \{ \rho_{R, v}
\}_{v \in M_K}$ est une mesure ad{\'e}lique.
Il reste {\`a} montrer que les hauteurs $h_R$ et $h_{\rho_R}$ co{\"\i}ncident.
%On donnera deux d{\'e}monstrations de ce fait.

%\subsubsection*{D{\'e}monstration~1.}
La preuve se base sur le lemme g{\'e}n{\'e}ral suivant.
\begin{lemme}
  Soit $\rho = \{ \rho_v \}$ une mesure ad{\'e}lique et soient $F$, $F'$
  et $F''$ des sous ensembles finis de $\P^1(\ov{K})$ invariants sous
  l'action du groupe de Galois $\Gal(\ov{K}/K)$.  Alors on~a les
  propri{\'e}t{\'e}s suivantes.
\begin{enumerate}
\item[1.]
Si $F''$ est disjoint de $F \cup F'$, alors on~a
\begin{equation}\label{E-difference}
\sum_{v \in M_K} \lpar [F], [F'], [F''], \rho_v \rpar_v
=
h_\rho(F) - h_\rho(F').
\end{equation}
\item[2.]
Si les ensembles $F$ et $F'$ sont disjoints, alors on~a
\begin{equation}\label{E-somme}
\sum_{v \in M_K} \lpar \rho_v, [F], \rho_v, [F'] \rpar_v
=
h_\rho(F) + h_\rho(F').
\end{equation}
\end{enumerate}
\end{lemme}
\begin{proof}
  On ne montrera que~\eqref{E-difference},
  l'{\'e}quation~\eqref{E-somme} se d{\'e}montrant de mani{\`e}\-re similaire.
  Notons tout d'abord que pour chaque place $v$ de $K$ on~a
\begin{equation}\label{E-1001}
([F], [F'], [F''], \rho_v)_v
=
([F], [F''])_v + ([F'], \rho_v)_v - ([F], \rho_v)_v - ([F'], [F''])_v.
\end{equation}
On a $([F], [F''])_v = -\log|D(F, F'')|_v$, avec $D(F, F'') =
\prod_{\a \in F, \a' \in F''} (\a - \a')$.  Comme les ensembles $F$ et
$F''$ sont disjoints et invariants sous l'action du groupe de Galois
$\Gal(\ov{K}/K)$, on~a $D(F, F'') \in K^*$. La formule du produit
implique alors $\sum_{v \in M_K} \lpar [F], [F''] \rpar_v = 0$.  Un
raisonnement analogue donne enfin l'{\'e}galit{\'e} $\sum_{v \in M_K} \lpar
[F'], [F''] \rpar_v = 0$.  En sommant~\eqref{E-1001} sur toutes les
places, le Lemme~\ref{definition equivalente} donne la formule
d{\'e}sir{\'e}e.
\end{proof}

Pour montrer que les hauteurs $h_{\rho_R}$ et $h_R$ co{\"\i}ncident, on
utilisera le fait que $h_R$ est l'unique hauteur de Weil, {\`a} une
constante multiplicative pr{\`e}s, satisfaisant la formule d'invariance
$h_R \circ R= D \cdot h_R$.

Notons que pour chaque sous-ensemble fini $F$ de $\P^1(\ov{K})$
invariant sous l'action de $\Gal(\ov{K}/K)$, les ensembles $R(F)$ et
$R^{-1}(F)$ sont finis et invariants sous l'action de
$\Gal(\ov{K}/K)$.  On a de plus $R_*[F] = [R(F)]$ et si $F$ ne
contient aucune valeur critique de $R$, alors chaque {\'e}l{\'e}ment de $F$ a
exactement $D$ pr{\'e}images par $R$ et donc $D^{-1} R^*[F] =
[R^{-1}(F)]$.

Etant donn{\'e}s des sous ensembles finis $F$, $F'$ de $\P^1(\ov{K})$,
invariants sous l'action du groupe de Galois $\Gal(\ov{K}/K)$, soit
$F''$ un sous ensemble fini de $\P^1(\ov{K})$ invariant sous l'action
de $\Gal(\ov{K}/K)$ ne contenant aucune valeur critique de $R$ et tel
que $R^{-1}(F'')$ soit disjoint de $F \cup F'$.  Alors les birapports
$$
(R_*[F], R_*[F'], [F''], \rho_{R, v})_v
=
([R(F)], [R(F')], [F''], \rho_{R, v})_v
$$
et
$$
([F], [F'], D^{-1}R^*[F''], \rho_{R, v})_v
=
([F], [F'], [R^{-1}(F'')], \rho_{R, v})_v
$$
sont bien d{\'e}finis et la Formule de Transformation implique
$$
([R(F)], [R(F')], [F''], \rho_{R, v})_v
=
D \cdot ([F], [F'], [R^{-1}(F'')], \rho_{R, v})_v ~.
$$
Lorsqu'on somme sur toutes les places, la formule~\eqref{E-difference} donne
$$
h_{\rho_R}(R(F)) - h_{\rho_R}(R(F')) = D (h_{\rho_R}(F) - h_{\rho_R}(F')).
$$
Ceci implique que la fonction $h_{\rho_R} \circ R - D \cdot
h_{\rho_R}$ est constante.

Un raisonnement similaire avec la Formule de Transformation
et~\eqref{E-somme}, implique que si $F'$ ne contient aucune valeur
critique de $R$ et si $R^{-1}(F')$ est disjoint de $F$, alors
$$
h_{\rho_R} (R(F)) + h_{\rho_R} (F')
=
D \cdot ( h_{\rho_R} (F) + h_{\rho_R} (R^{-1}(F'))).
$$
On conclut que la fonction $h_{\rho_R} \circ R - D\,
h_{\rho_R}$ est identiquement nulle, ce qui termine la d{\'e}monstration
du Th{\'e}or{\`e}me~\ref{normalisee est arakelov}.
\subsection{Ensembles de Mandelbrot.}\label{S-mandel}
Notre but est de d{\'e}montrer le Th{\'e}or{\`e}me~\ref{thm:mandel}.  Rappelons le
contexte. On fixe un entier $D \ge 2$, et on note $P_c(z) = z^D + c$
pour $c \in \C$.  L'ensemble de Mandelbrot $\cM_D$ est l'ensemble des
param{\`e}tres complexes pour lesquels l'orbite du point critique~$z = 0$
de $P_c$ est born{\'e}e, \ie $\cM_D = \{ c \in \C ; |P^n_c(0)| = O(1) \}$.

On introduit alors la famille de mesures suivantes. Pour tout nombre
premier $p$, on pose $\rho_p = \lambda_p$. A la place infinie, on
d{\'e}finit $\rho_\infty$ comme {\'e}tant la mesure harmonique de $\cM_D$.
\begin{lemme}
La famille $\rho_{\cM_D} = \{\rho_v\}_{v\in M_\Q}$ est une mesure
ad{\'e}lique {\`a} potentiel H{\"o}lder. De plus, pour toute place $v$, on~a
$\rho_v = \lim_{n\to\infty} D^{-n} \Delta \log^+|P^n_c(0)|$.
\end{lemme}
\begin{proof}
  On montre tout d'abord la seconde assertion.  En une place $v$
  finie, on proc{\`e}de comme suit. Si $|c|\le1$, on~a $|P^n_c(0)| \le 1$
  pour tout $n$ donc $\lim_{n\to\infty} D^{-n} \Delta
  \log^+|P^n_c(0)|=0$. Lorsque $|c|>1$, on v{\'e}rifie par r{\'e}currence que
  $|P^n_c(0)| = |c|^{D^n}$. Donc $\rho_v = \Delta \log^+|z| =
  \lambda_p$. A la place infinie, la preuve que $G =\lim_{n\to\infty}
  D^{-n} \Delta \log^+|P^n_c(0)|$ est la mesure harmonique de $\cM_D$
  est compl{\`e}tement classique: il faut pour cela montrer tout
  d'abord l'existence de la limite et que $G$ est continue; puis que
  $G=0$ sur $\cM_D$; que $G$ est harmonique sur $\C \setminus \cM_D$;
  enfin que $\Delta G$ est une mesure de probabilit{\'e}.  Tout ceci est
  explicitement fait dans le cas $D=2$ dans~\cite[VIII.4]{CG}, et les
  preuves s'adaptent au cas $D\ge3$ sans modification.
  De~\cite[Theorem~3.2]{CG}, on tire enfin le caract{\`e}re H{\"o}ld{\'e}rien
  de~$G$.
\end{proof}
On peut donc introduire la hauteur $h_{\cM_D}$ associ{\'e}e {\`a} la mesure
ad{\'e}lique $\rho_{\cM_D}$ d{\'e}finie dans le lemme pr{\'e}c{\'e}dent.  Le
Th{\'e}or{\`e}me~\ref{thm:mandel} est maintenant une cons{\'e}quence imm{\'e}diate du
Th{\'e}or{\`e}me~\ref{T-estim} et du lemme suivant.  Rappelons qu'un param{\`e}tre
$c \in \C$ est critiquement fini, si on peut trouver des entiers distincts
$n$ et $m$ tels que $P_c^n(0) = P_c^m(0)$.
\begin{lemme}
  Pour tout $c\in \ov{\Q}$, on~a $h_{\cM_D}(c) = D^{-1} h_{P_c}(c)$.
  En particulier, pour un param{\`e}tre $c \in \C$ on~a $h_{\cM_D}(c) = 0$
  si et seulement si $c$ est critiquement fini.
\end{lemme}
\begin{proof}
  Soit $c \in \ov{\Q}$ et soit $F$ l'orbite de $c$ sous l'action du
  groupe de Galois $\Gal(\ov{\Q} / \Q)$.  La preuve r{\'e}sulte alors du
  calcul suivant.
\begin{multline*}
  h_{P_c}(c) = \lim_{n \to \infty} D^{-n} \hn (P^n_c(c)) = \lim_{n \to
    \infty} D^{-n} |F|^{-1} \sum_{\alpha \in F} \sum_{v \in M_\Q}
  \log^+|P_{\alpha}^n(\alpha)|_v =\\= D |F|^{-1} \sum_{\alpha \in F}
  \sum_{v \in M_\Q} \lim_{n \to \infty} D^{-(n + 1)}\log^+|P_\alpha^{n
    + 1}(0)|_v = D \cdot h_{\cM_D}(c)~.
\end{multline*}
La permutation des signes sommes et du limite se fait
sans souci car seul un nombre fini de termes de la somme est
non-nuls.  Enfin la derni{\`e}re {\'e}galit{\'e} est une cons{\'e}quence du
lemme pr{\'e}c{\'e}dent et des d{\'e}finitions.
\end{proof}

%%%%%%%%%%%%%%%%%%%%%%%%%%%%%%%%%%%%%%%%%%%%%%%%%%%%%%%%
\subsection{Continuit{\'e} H{\"o}lder.}\label{S-holder}
Dans cette partie, nous donnons une preuve du fait suivant.
\begin{proposition}
  Soit $R$ une fraction rationnelle d{\'e}finie sur $\C$ ou sur $\C_p$
  pour un $p$ premier. Alors la mesure d'{\'e}quilibre $\rho_R$ admet un
  potentiel H{\"o}lder.
\end{proposition}
Ce fait est classique sur le corps des complexes, et reste valable en
dimension sup{\'e}rieure. Nous renvoyons {\`a}~\cite[Section~1.7]{sibony}
pour une preuve dans un cadre assez g{\'e}n{\'e}ral. Nous donnons ici une
preuve valable aussi bien dans le cas complexe que dans le cas
$p$-adique.
\begin{proof}
  Pour fixer les notations, on se place dans le cas $p$-adique. Soit
  $D \ge 2$ le degr{\'e} de $R$ et soit $g$ un potentiel tel que
  $D^{-1}R^*\lambda_p - \lambda_p = \Delta g$.  La remarque
  essentielle est que $g$ induit sur $\berp$ une fonction
  lipschitzienne, et que l'on peut donc trouver une constante $C>0$
  telle que $|g(z) - g(w)|\le C\cdot d(z,w)$ pour tout $z,w\in\berp$.
  La fonction $R$ est aussi lipshitzienne pour une autre constante $M$
  que l'on peut supposer strictement plus grande que $D$. De $\rho_R
  =\lim_n D^{-n} R^{n*} \lambda_p$, on en d{\'e}duit que $\rho_R
  =\lambda_p + \Delta g_\infty$ avec $g_\infty = \sum_0^\infty D^{-k}
  g\circ R^k$.  On choisit $N$ entier, et on proc{\`e}de alors aux
  estimations suivantes.
\begin{multline*}
|g_\infty(z)-g_\infty(w)|\le \\
\le \sum_0^{N-1} D^{-k} |g\circ R^k(z) - g\circ R^k(w)| + \sup|g|\,
\sum_N^\infty  D^{-k}\le \\
\le C_1 \cdot d(z,w)\, (M/D)^N + C_2\cdot D^{-N} ~,
\end{multline*}
pour des constantes $C_1,C_2>0$.
On choisit alors $N$ de telle sorte que les deux termes de la
derni{\`e}re somme soient {\'e}gaux. On en d{\'e}duit que
$|g_\infty(z)-g_\infty(w)|\le C_3 \cdot d(z,w)^\kappa$ avec $\kappa =
\log D /\log M$.
\end{proof}

%%%%%%%%%%%%%%%%%%%%%%%%%%%%%%%%%%%%%%%%%%%%%%%%%%%%%%%%
\subsection{Un exemple.}\label{S-exemple}
Soit $K$ un corps de nombres et fixons une place finie $v$ de $K$ ne
divisant pas~$2$.  Etant donn{\'e}e une constante $C \in \ov{K}$ de norme
$|C|_v >1$, posons $P(T) \= T^2 + C$.  Pour chaque entier positif $n$
posons $F_n = \{ z \in \C_v;\, P^n(z) = z \}$. Notons que l'on a
$h_P(F_n) = 0$.  On va montrer que
\begin{equation}\label{e3827}
([F_n]-\rho_{P,v},[F_n]-\rho_{P,v})_v =
  -\frac{n}{2^{2n}} \log \sqrt{|C|}   < 0~.
\end{equation}
Comme ce terme local de la hauteur normalis{\'e}e $h_P$ est n{\'e}gatif, cet
exemple montre que le proc{\'e}d{\'e} de r{\'e}gularisation dans la preuve du
Th{\'e}or{\`e}me~\ref{T-prec} est en effet necessaire.

Le calcul de cette quantit{\'e} peut {\^e}tre r{\'e}alis{\'e} de la mani{\`e}re
suivante. Pour simplifer notons $\rho \= \rho_{P,v}$, et posons $g (z)
= \int_{\C_v} \log | z - w| \, d\rho(z)$.  On v{\'e}rifie que $g =
\lim_{n\to \infty} 2^{-n} \log^+ | P^n|$ et en particulier que
$g\equiv 0$ sur $F_n$. Le Lemme~\ref{L-calculCp} implique alors
$([F_n],\rho) = - \int g\, d[F_n] = 0$. Du fait que $P$ est unitaire,
pour tout $z,w\in \C_v$ on~a $P(z) -w = \prod_{P(w_i) = w} (z -
w_i)$, et on en d{\'e}duit que
\begin{multline*}
  2\, (\rho,\rho) = (\rho ,P^*\rho) = -\int_{\C_v} \sum_{P(w_i)=w}
  \log |z-w_i| \, d\rho \otimes d\rho =\\ -\int_{\C_v} \log |P(z) -w|
  \, d\rho \otimes d\rho = (P_*\rho , \rho).
\end{multline*}
On verra plus bas que le support de $\rho$ est inclus dans $\C_v$, ce
qui explique le fait que dans le cacul pr{\'e}c{\'e}dent l'on puisse
remplacer $\log \sup\{ \cS , \cS'\}$ sur l'espace de Berkovich par la
fonction $\log|z-w|$ sur $\C_v$.  Comme $P_*\rho = \rho$, il s'ensuit
$(\rho,\rho) =0$, et par suite $([F_n]-\rho_{P,v},[F_n]-\rho_{P,v})=
([F_n], [F_n])$. Pour calculer cette derni{\`e}re quantit{\'e} on proc{\`e}de
comme suit.

Posons $r \= \sqrt{ |C|}$ et
$$
K_P = \{ z \in \C_p; \; \{ |P^n(z)| \}_{n \ge 0} \mbox{ est born{\'e}e} \}.
$$
Il est facile de voir que $K_P$ est contenu dans la boule ferm{\'e}e
$B(0, r)$.  La pr{\'e}image de cette boule est {\'e}gale {\`a} la r{\'e}union des
boules $B_0$ et $B_1$, centr{\'e}es aux racines carr{\'e}es de $C$ dans $\C_v$
et de diam{\`e}tre~$1$.  Les boules $B_0$ et $B_1$ sont incluses dans
$\{|z| = r\}$ et pour tout couple $(z_0,z_1)\in B_0 \times B_1$ on~a
$|z_0 - z_1| = r$ (on utilise ici le fait que $v$ ne divise pas~$2$).
D'autre part, si $z, z'$ sont dans une m{\^e}me boule $B_i$ on v{\'e}rifie
ais{\'e}ment que $|P(z) - P(z')| = r |z - z'|$.

A tout point $z$ dans $K_P$ est associ{\'e} un codage $\e (z)= \{ \e_i
\}_{i\ge0} \in \{ 0,1 \}^{\N}$ de telle sorte que $P^i(z) \in
B_{\e_i}$ pour tout $i$.  Pour tout $z,z' \in K_P$, la distance
$|z-z'|$ ne d{\'e}pend donc que du codage de $z$ et $z'$.  Pour $\e, \e'
\in \{ 0, 1 \}^{\N}$ distincts, soit $k(\e,\e')$ l'entier positif tel
que $\e_i = \e_i$ pour tout $i < k$ et $\e_{k}(z) \neq \e_{k}(z')$.
Alors pour $z, z' \in K_P$ distincts, on~a $\log |z-z'|/\log r =
1-k(\e(z), \e(z'))$. On montre facilement que l'application $\e : K_P
\to \{ 0, 1 \}^{\N}$ est un hom{\'e}omorphisme qui conjugue $P$ au decalage
$\{ \e_i \}_{i\ge 0} \mapsto \{ \e_{i + 1} \}_{i \ge 0}$. On en d{\'e}duit
d'une part que $\rho$ est support{\'e}e sur $K_P$, et d'autre part que
pour chaque $n \ge 1$ l'ensemble $F_n$ est en bijection par $\e$ avec
les suites finies de $\{ 0,1 \}$ de longueur $n$.  On a donc ramen{\'e} le
calcul de $([F_n],[F_n])$ {\`a} un calcul purement combinatoire:
$$
2^{2n}([F_n],[F_n])
= 
- \sum_{z\neq z' \in F_n} \log |z-z'|
=
- \sum_{\e \neq \e' \in \{0,1\}^n}  (1-k(\e,\e'))\, \log r~. 
$$
On d{\'e}duit~\eqref{e3827} de la
formule pr{\'e}c{\'e}dente par un calcul direct (mais fastidieux).

%
%%%%%%%%%%%%%%%%%%%%%%%%%%%%%%%%%%%%%%%%%%%%%%%%%%%%%%%%%%%%%%%%%%
%

%
%%%%%%%%%%%%%%%%%%%%%%%%%%%%%%%%%%%%%%%%%%%%%%%%%%%%%%%%%%%%%%%%%%
%

\end{document}